\documentclass{elsarticle}
  
\usepackage[margin=1in]{geometry}

\usepackage{bbm}
\usepackage{graphicx}
\usepackage{pstool}
\usepackage{wrapfig}
\usepackage{caption}
\usepackage{subcaption}
\usepackage[section]{placeins} 

\usepackage{fancyhdr} 
\usepackage{lastpage} 
\usepackage{extramarks} 

\usepackage{xcolor} 
\usepackage{enumerate} 
\usepackage{paralist} 
\usepackage{amsmath, amsthm, amssymb, mathtools}
\usepackage{mathabx, pifont, stmaryrd} 
\usepackage[explicit]{titlesec} 
\usepackage{etoolbox} 
\usepackage{bibentry} 
\makeatletter\let\saved@bibitem\@bibitem\makeatother 
\usepackage[colorlinks, bookmarksopen, bookmarksnumbered,
            citecolor=red,urlcolor=red]{hyperref} 
\makeatletter\let\@bibitem\saved@bibitem\makeatother 
\usepackage{rotating}

\usepackage{algorithm}
\usepackage{algorithmic}

\newtheorem{remark}{Remark}

\theoremstyle{definition}

\usepackage{bm} 
\usepackage{amsfonts} 

\DeclareMathOperator*{\argmin}{arg\,min}

\newcommand{\ds}[1]{\ensuremath{\displaystyle{#1}}}
\newcommand{\func}[3]{\ensuremath{#1 : #2 \rightarrow #3}}

\newcommand{\norm}[1]{\ensuremath{\left\| #1 \right\|}}

\newcommand{\suchthat}{\mathrel{}\middle|\mathrel{}}

\newcommand{\pder}[2]{\ensuremath{\frac{\partial #1}{\partial #2}}}

\newcommand{\pderH}[3]{\ensuremath{\frac{\partial^{#3} #1}{\partial #2^{#3}}}}

\newcommand{\Ecal}{\ensuremath{\mathcal{E}}}

\newcommand{\Gcal}{\ensuremath{\mathcal{G}}}
\newcommand{\Hcal}{\ensuremath{\mathcal{H}}}

\newcommand{\Jcal}{\ensuremath{\mathcal{J}}}

\newcommand{\Ocal}{\ensuremath{\mathcal{O}}}
\newcommand{\Pcal}{\ensuremath{\mathcal{P}}}
\newcommand{\Qcal}{\ensuremath{\mathcal{Q}}}

\newcommand{\Scal}{\ensuremath{\mathcal{S}}}
\newcommand{\Tcal}{\ensuremath{\mathcal{T}}}

\newcommand{\Vcal}{\ensuremath{\mathcal{V}}}
\newcommand{\Wcal}{\ensuremath{\mathcal{W}}}

\newcommand{\Gbb}{\ensuremath{\mathbb{G}}}

\newcommand{\Nbb}{\ensuremath{\mathbb{N} }}

\newcommand{\Rbb}{\ensuremath{\mathbb{R} }}

\newcommand\Bbm{{\ensuremath{\bm{B}}}}

\newcommand\Dbm{{\ensuremath{\bm{D}}}}

\newcommand\Fbm{{\ensuremath{\bm{F}}}}

\newcommand\Jbm{{\ensuremath{\bm{J}}}}

\newcommand\Pbm{{\ensuremath{\bm{P}}}}
\newcommand\Qbm{{\ensuremath{\bm{Q}}}}
\newcommand\Rbm{{\ensuremath{\bm{R}}}}

\newcommand\gbm{{\ensuremath{\bm{g}}}}

\newcommand\rbm{{\ensuremath{\bm{r}}}}

\newcommand\ubm{{\ensuremath{\bm{u}}}}

\newcommand\xbm{{\ensuremath{\bm{x}}}}
\newcommand\ybm{{\ensuremath{\bm{y}}}}
\newcommand\zbm{{\ensuremath{\bm{z}}}}

\newcommand\phibold{{\ensuremath{\boldsymbol{\phi}}}}

\newcommand\Psibold{{\ensuremath{\boldsymbol{\Psi}}}}

\newcommand\zerobold{\ensuremath{\mathbf{0}}}

\usepackage{tikz}
\usepackage{pgfplots}
\usepackage{pgfplotstable, filecontents, booktabs}
\pgfplotsset{compat=1.9}

\usetikzlibrary{pgfplots.groupplots}
\usepgfplotslibrary{fillbetween}
\usetikzlibrary{calc,fit,matrix,arrows,automata,positioning,shapes}
\usetikzlibrary{arrows.meta}

\pgfplotsset{select coords between index/.style 2 args={
    x filter/.code={
        \ifnum\coordindex<#1\fi
        \ifnum\coordindex>#2\fi
    }
}}

\tikzset{
 invisible/.style={opacity=0},
 visible on/.style={alt={#1{}{invisible}}},
 alt/.code args={<#1>#2#3}{%
   \alt<#1>{\pgfkeysalso{#2}}{\pgfkeysalso{#3}}
 },
}

\newcommand{\colorbarMatlabParula}[5]{
\begin{tikzpicture}
\begin{axis}[
   hide axis, scale only axis,
   height=0pt, width=0pt,
   colormap={parula}{rgb255=(62,38,168) rgb255=(62,39,172) rgb255=(63,40,175) rgb255=(63,41,178) rgb255=(64,42,180) rgb255=(64,43,183) rgb255=(65,44,186) rgb255=(65,45,189) rgb255=(66,46,191) rgb255=(66,47,194) rgb255=(67,48,197) rgb255=(67,49,200) rgb255=(67,50,202) rgb255=(68,51,205) rgb255=(68,52,208) rgb255=(69,53,210) rgb255=(69,55,213) rgb255=(69,56,215) rgb255=(70,57,217) rgb255=(70,58,220) rgb255=(70,59,222) rgb255=(70,61,224) rgb255=(71,62,225) rgb255=(71,63,227) rgb255=(71,65,229) rgb255=(71,66,230) rgb255=(71,68,232) rgb255=(71,69,233) rgb255=(71,70,235) rgb255=(72,72,236) rgb255=(72,73,237) rgb255=(72,75,238) rgb255=(72,76,240) rgb255=(72,78,241) rgb255=(72,79,242) rgb255=(72,80,243) rgb255=(72,82,244) rgb255=(72,83,245) rgb255=(72,84,246) rgb255=(71,86,247) rgb255=(71,87,247) rgb255=(71,89,248) rgb255=(71,90,249) rgb255=(71,91,250) rgb255=(71,93,250) rgb255=(70,94,251) rgb255=(70,96,251) rgb255=(70,97,252) rgb255=(69,98,252) rgb255=(69,100,253) rgb255=(68,101,253) rgb255=(67,103,253) rgb255=(67,104,254) rgb255=(66,106,254) rgb255=(65,107,254) rgb255=(64,109,254) rgb255=(63,110,255) rgb255=(62,112,255) rgb255=(60,113,255) rgb255=(59,115,255) rgb255=(57,116,255) rgb255=(56,118,254) rgb255=(54,119,254) rgb255=(53,121,253) rgb255=(51,122,253) rgb255=(50,124,252) rgb255=(49,125,252) rgb255=(48,127,251) rgb255=(47,128,250) rgb255=(47,130,250) rgb255=(46,131,249) rgb255=(46,132,248) rgb255=(46,134,248) rgb255=(46,135,247) rgb255=(45,136,246) rgb255=(45,138,245) rgb255=(45,139,244) rgb255=(45,140,243) rgb255=(45,142,242) rgb255=(44,143,241) rgb255=(44,144,240) rgb255=(43,145,239) rgb255=(42,147,238) rgb255=(41,148,237) rgb255=(40,149,236) rgb255=(39,151,235) rgb255=(39,152,234) rgb255=(38,153,233) rgb255=(38,154,232) rgb255=(37,155,232) rgb255=(37,156,231) rgb255=(36,158,230) rgb255=(36,159,229) rgb255=(35,160,229) rgb255=(35,161,228) rgb255=(34,162,228) rgb255=(33,163,227) rgb255=(32,165,227) rgb255=(31,166,226) rgb255=(30,167,225) rgb255=(29,168,225) rgb255=(29,169,224) rgb255=(28,170,223) rgb255=(27,171,222) rgb255=(26,172,221) rgb255=(25,173,220) rgb255=(23,174,218) rgb255=(22,175,217) rgb255=(20,176,216) rgb255=(18,177,214) rgb255=(16,178,213) rgb255=(14,179,212) rgb255=(11,179,210) rgb255=(8,180,209) rgb255=(6,181,207) rgb255=(4,182,206) rgb255=(2,183,204) rgb255=(1,183,202) rgb255=(0,184,201) rgb255=(0,185,199) rgb255=(0,186,198) rgb255=(1,186,196) rgb255=(2,187,194) rgb255=(4,187,193) rgb255=(6,188,191) rgb255=(9,189,189) rgb255=(13,189,188) rgb255=(16,190,186) rgb255=(20,190,184) rgb255=(23,191,182) rgb255=(26,192,181) rgb255=(29,192,179) rgb255=(32,193,177) rgb255=(35,193,175) rgb255=(37,194,174) rgb255=(39,194,172) rgb255=(41,195,170) rgb255=(43,195,168) rgb255=(44,196,166) rgb255=(46,196,165) rgb255=(47,197,163) rgb255=(49,197,161) rgb255=(50,198,159) rgb255=(51,199,157) rgb255=(53,199,155) rgb255=(54,200,153) rgb255=(56,200,150) rgb255=(57,201,148) rgb255=(59,201,146) rgb255=(61,202,144) rgb255=(64,202,141) rgb255=(66,202,139) rgb255=(69,203,137) rgb255=(72,203,134) rgb255=(75,203,132) rgb255=(78,204,129) rgb255=(81,204,127) rgb255=(84,204,124) rgb255=(87,204,122) rgb255=(90,204,119) rgb255=(94,205,116) rgb255=(97,205,114) rgb255=(100,205,111) rgb255=(103,205,108) rgb255=(107,205,105) rgb255=(110,205,102) rgb255=(114,205,100) rgb255=(118,204,97) rgb255=(121,204,94) rgb255=(125,204,91) rgb255=(129,204,89) rgb255=(132,204,86) rgb255=(136,203,83) rgb255=(139,203,81) rgb255=(143,203,78) rgb255=(147,202,75) rgb255=(150,202,72) rgb255=(154,201,70) rgb255=(157,201,67) rgb255=(161,200,64) rgb255=(164,200,62) rgb255=(167,199,59) rgb255=(171,199,57) rgb255=(174,198,55) rgb255=(178,198,53) rgb255=(181,197,51) rgb255=(184,196,49) rgb255=(187,196,47) rgb255=(190,195,45) rgb255=(194,195,44) rgb255=(197,194,42) rgb255=(200,193,41) rgb255=(203,193,40) rgb255=(206,192,39) rgb255=(208,191,39) rgb255=(211,191,39) rgb255=(214,190,39) rgb255=(217,190,40) rgb255=(219,189,40) rgb255=(222,188,41) rgb255=(225,188,42) rgb255=(227,188,43) rgb255=(230,187,45) rgb255=(232,187,46) rgb255=(234,186,48) rgb255=(236,186,50) rgb255=(239,186,53) rgb255=(241,186,55) rgb255=(243,186,57) rgb255=(245,186,59) rgb255=(247,186,61) rgb255=(249,186,62) rgb255=(251,187,62) rgb255=(252,188,62) rgb255=(254,189,61) rgb255=(254,190,60) rgb255=(254,192,59) rgb255=(254,193,58) rgb255=(254,194,57) rgb255=(254,196,56) rgb255=(254,197,55) rgb255=(254,199,53) rgb255=(254,200,52) rgb255=(254,202,51) rgb255=(253,203,50) rgb255=(253,205,49) rgb255=(253,206,49) rgb255=(252,208,48) rgb255=(251,210,47) rgb255=(251,211,46) rgb255=(250,213,46) rgb255=(249,214,45) rgb255=(249,216,44) rgb255=(248,217,43) rgb255=(247,219,42) rgb255=(247,221,42) rgb255=(246,222,41) rgb255=(246,224,40) rgb255=(245,225,40) rgb255=(245,227,39) rgb255=(245,229,38) rgb255=(245,230,38) rgb255=(245,232,37) rgb255=(245,233,36) rgb255=(245,235,35) rgb255=(245,236,34) rgb255=(245,238,33) rgb255=(246,239,32) rgb255=(246,241,31) rgb255=(246,242,30) rgb255=(247,244,28) rgb255=(247,245,27) rgb255=(248,247,26) rgb255=(248,248,24) rgb255=(249,249,22) rgb255=(249,251,21) },
   colorbar horizontal,
   point meta min=#1, point meta max=#5,
   colorbar style={width=10cm, xtick={#1,#2,#3,#4,#5}}
]
\addplot [draw=none] coordinates {(0,0)};
\end{axis}
\end{tikzpicture}
}

\newcommand{\colorbarMatlabParulandscalse}[7]{
\begin{tikzpicture}
\begin{axis}[
   hide axis, scale only axis,
   height=0pt, width=0pt,
   colormap={parula}{rgb255=(62,38,168) rgb255=(62,39,172) rgb255=(63,40,175) rgb255=(63,41,178) rgb255=(64,42,180) rgb255=(64,43,183) rgb255=(65,44,186) rgb255=(65,45,189) rgb255=(66,46,191) rgb255=(66,47,194) rgb255=(67,48,197) rgb255=(67,49,200) rgb255=(67,50,202) rgb255=(68,51,205) rgb255=(68,52,208) rgb255=(69,53,210) rgb255=(69,55,213) rgb255=(69,56,215) rgb255=(70,57,217) rgb255=(70,58,220) rgb255=(70,59,222) rgb255=(70,61,224) rgb255=(71,62,225) rgb255=(71,63,227) rgb255=(71,65,229) rgb255=(71,66,230) rgb255=(71,68,232) rgb255=(71,69,233) rgb255=(71,70,235) rgb255=(72,72,236) rgb255=(72,73,237) rgb255=(72,75,238) rgb255=(72,76,240) rgb255=(72,78,241) rgb255=(72,79,242) rgb255=(72,80,243) rgb255=(72,82,244) rgb255=(72,83,245) rgb255=(72,84,246) rgb255=(71,86,247) rgb255=(71,87,247) rgb255=(71,89,248) rgb255=(71,90,249) rgb255=(71,91,250) rgb255=(71,93,250) rgb255=(70,94,251) rgb255=(70,96,251) rgb255=(70,97,252) rgb255=(69,98,252) rgb255=(69,100,253) rgb255=(68,101,253) rgb255=(67,103,253) rgb255=(67,104,254) rgb255=(66,106,254) rgb255=(65,107,254) rgb255=(64,109,254) rgb255=(63,110,255) rgb255=(62,112,255) rgb255=(60,113,255) rgb255=(59,115,255) rgb255=(57,116,255) rgb255=(56,118,254) rgb255=(54,119,254) rgb255=(53,121,253) rgb255=(51,122,253) rgb255=(50,124,252) rgb255=(49,125,252) rgb255=(48,127,251) rgb255=(47,128,250) rgb255=(47,130,250) rgb255=(46,131,249) rgb255=(46,132,248) rgb255=(46,134,248) rgb255=(46,135,247) rgb255=(45,136,246) rgb255=(45,138,245) rgb255=(45,139,244) rgb255=(45,140,243) rgb255=(45,142,242) rgb255=(44,143,241) rgb255=(44,144,240) rgb255=(43,145,239) rgb255=(42,147,238) rgb255=(41,148,237) rgb255=(40,149,236) rgb255=(39,151,235) rgb255=(39,152,234) rgb255=(38,153,233) rgb255=(38,154,232) rgb255=(37,155,232) rgb255=(37,156,231) rgb255=(36,158,230) rgb255=(36,159,229) rgb255=(35,160,229) rgb255=(35,161,228) rgb255=(34,162,228) rgb255=(33,163,227) rgb255=(32,165,227) rgb255=(31,166,226) rgb255=(30,167,225) rgb255=(29,168,225) rgb255=(29,169,224) rgb255=(28,170,223) rgb255=(27,171,222) rgb255=(26,172,221) rgb255=(25,173,220) rgb255=(23,174,218) rgb255=(22,175,217) rgb255=(20,176,216) rgb255=(18,177,214) rgb255=(16,178,213) rgb255=(14,179,212) rgb255=(11,179,210) rgb255=(8,180,209) rgb255=(6,181,207) rgb255=(4,182,206) rgb255=(2,183,204) rgb255=(1,183,202) rgb255=(0,184,201) rgb255=(0,185,199) rgb255=(0,186,198) rgb255=(1,186,196) rgb255=(2,187,194) rgb255=(4,187,193) rgb255=(6,188,191) rgb255=(9,189,189) rgb255=(13,189,188) rgb255=(16,190,186) rgb255=(20,190,184) rgb255=(23,191,182) rgb255=(26,192,181) rgb255=(29,192,179) rgb255=(32,193,177) rgb255=(35,193,175) rgb255=(37,194,174) rgb255=(39,194,172) rgb255=(41,195,170) rgb255=(43,195,168) rgb255=(44,196,166) rgb255=(46,196,165) rgb255=(47,197,163) rgb255=(49,197,161) rgb255=(50,198,159) rgb255=(51,199,157) rgb255=(53,199,155) rgb255=(54,200,153) rgb255=(56,200,150) rgb255=(57,201,148) rgb255=(59,201,146) rgb255=(61,202,144) rgb255=(64,202,141) rgb255=(66,202,139) rgb255=(69,203,137) rgb255=(72,203,134) rgb255=(75,203,132) rgb255=(78,204,129) rgb255=(81,204,127) rgb255=(84,204,124) rgb255=(87,204,122) rgb255=(90,204,119) rgb255=(94,205,116) rgb255=(97,205,114) rgb255=(100,205,111) rgb255=(103,205,108) rgb255=(107,205,105) rgb255=(110,205,102) rgb255=(114,205,100) rgb255=(118,204,97) rgb255=(121,204,94) rgb255=(125,204,91) rgb255=(129,204,89) rgb255=(132,204,86) rgb255=(136,203,83) rgb255=(139,203,81) rgb255=(143,203,78) rgb255=(147,202,75) rgb255=(150,202,72) rgb255=(154,201,70) rgb255=(157,201,67) rgb255=(161,200,64) rgb255=(164,200,62) rgb255=(167,199,59) rgb255=(171,199,57) rgb255=(174,198,55) rgb255=(178,198,53) rgb255=(181,197,51) rgb255=(184,196,49) rgb255=(187,196,47) rgb255=(190,195,45) rgb255=(194,195,44) rgb255=(197,194,42) rgb255=(200,193,41) rgb255=(203,193,40) rgb255=(206,192,39) rgb255=(208,191,39) rgb255=(211,191,39) rgb255=(214,190,39) rgb255=(217,190,40) rgb255=(219,189,40) rgb255=(222,188,41) rgb255=(225,188,42) rgb255=(227,188,43) rgb255=(230,187,45) rgb255=(232,187,46) rgb255=(234,186,48) rgb255=(236,186,50) rgb255=(239,186,53) rgb255=(241,186,55) rgb255=(243,186,57) rgb255=(245,186,59) rgb255=(247,186,61) rgb255=(249,186,62) rgb255=(251,187,62) rgb255=(252,188,62) rgb255=(254,189,61) rgb255=(254,190,60) rgb255=(254,192,59) rgb255=(254,193,58) rgb255=(254,194,57) rgb255=(254,196,56) rgb255=(254,197,55) rgb255=(254,199,53) rgb255=(254,200,52) rgb255=(254,202,51) rgb255=(253,203,50) rgb255=(253,205,49) rgb255=(253,206,49) rgb255=(252,208,48) rgb255=(251,210,47) rgb255=(251,211,46) rgb255=(250,213,46) rgb255=(249,214,45) rgb255=(249,216,44) rgb255=(248,217,43) rgb255=(247,219,42) rgb255=(247,221,42) rgb255=(246,222,41) rgb255=(246,224,40) rgb255=(245,225,40) rgb255=(245,227,39) rgb255=(245,229,38) rgb255=(245,230,38) rgb255=(245,232,37) rgb255=(245,233,36) rgb255=(245,235,35) rgb255=(245,236,34) rgb255=(245,238,33) rgb255=(246,239,32) rgb255=(246,241,31) rgb255=(246,242,30) rgb255=(247,244,28) rgb255=(247,245,27) rgb255=(248,247,26) rgb255=(248,248,24) rgb255=(249,249,22) rgb255=(249,251,21) },
   colorbar horizontal,
   point meta min=#1, point meta max=#6,
   colorbar style={width=#7, xtick={#1,#2,#3,#4, #5,#6}}
]
\addplot [draw=none] coordinates {(0,0)};
\end{axis}
\end{tikzpicture}
}

\usepackage{setspace}

\newbool{fastcompile}
\setbool{fastcompile}{false}

\begin{document}
\title{An $rp$-adaptive method for accurate resolution of shock-dominated viscous flow based on implicit shock tracking}

\author[rvt1]{Huijing Dong\fnref{fn1}}
\ead{hdong2@nd.edu}

\author[rvt2]{Masayuki Yano\fnref{fn2}}
\ead{masa.yano@utoronto.ca}

\author[rvt1]{Tianci Huang\fnref{fn3}}
\ead{thuang5@alumni.nd.edu}

\author[rvt1]{Matthew J. Zahr\fnref{fn4}\corref{cor1}}
\ead{mzahr@nd.edu}

\address[rvt1]{Department of Aerospace and Mechanical Engineering, University
               of Notre Dame, Notre Dame, IN 46556, United States}
\address[rvt2]{Institute for Aerospace Studies, University
	      of Toronto, Toronto, M3H 5T6, Ontario, Canada}
\cortext[cor1]{Corresponding author}

\fntext[fn1]{Graduate Student, Department of Aerospace and Mechanical
             Engineering, University of Notre Dame}
\fntext[fn2]{Associate Professor, Institute for Aerospace Studies, University of Toronto}
\fntext[fn3]{Research Associate, Department of Aerospace and Mechanical
             Engineering, University of Notre Dame}
\fntext[fn4]{Assistant Professor, Department of Aerospace and Mechanical
             Engineering, University of Notre Dame}


\begin{abstract}
	This work introduces an optimization-based $rp$-adaptive numerical method to
	approximate solutions of viscous, shock-dominated flows using implicit shock
	tracking and a high-order discontinuous Galerkin discretization on traditionally
	coarse grids without nonlinear stabilization (e.g., artificial viscosity or
	limiting). The proposed method adapts implicit shock tracking methods,
	originally developed to align mesh faces with solution discontinuities, to
	compress elements into viscous shocks and boundary layers, functioning as
	a novel approach to aggressive $r$-adaptation. This form of $r$-adaptation
	is achieved naturally as the minimizer of the enriched residual with respect
	to the discrete flow variables and coordinates of the nodes of the grid.
	Several innovations to the shock tracking optimization solver are
	proposed to ensure sufficient mesh compression at viscous features
	to render stabilization unnecessary, including residual weighting,
	step constraints and modifications, and viscosity-based continuation.
	Finally, $p$-adaptivity is used to locally increase the polynomial
	degree with three clear benefits: (1) lessens the mesh compression
	requirements near shock waves and boundary layers, (2) reduces the
	error in regions where $r$-adaptivity is not sufficient with the
	given grid topology, and (3) reduces computational cost by performing
	a majority of the $r$-adaptivity iterations on the coarsest discretization.
	A series of numerical experiments show the proposed method effectively resolves
	viscous, shock-dominated flows, including accurate prediction of heat flux
	profiles produced by hypersonic flow over a cylinder, and compares
	favorably in terms of accuracy per degree of freedom to $h$-adaptation
	with a high-order discretization.
\end{abstract}

\maketitle

\section{Introduction}
\label{sec:intro}
Computational fluid dynamics (CFD) has become an essential tool for gaining a
deeper understanding of flow physics and enabling the prediction and control
of complex fluid flows. However, high-speed flows, particularly in the hypersonic
regime, present significant challenges even for modern CFD methods. One major
difficulty is their sensitivity to the computational mesh as tight grid spacing
near and precise alignment with shock waves and boundary layers are
needed \cite{lee1999spurious, candler2009current, gnoffo2009multi, candler2015next}.
Properly aligning the grid with shocks is a user-intensive process, often requiring
iterating between grid generation and flow simulation \cite{candler2009current}.
Additionally, nonlinear instabilities become increasingly problematic as the
Mach number and vehicle bluntness increase \cite{candler2015advances}, which
further complicates the stability and accuracy of numerical methods.

Shock-dominated flows are typically handled using shock capturing methods, which
employ various techniques to stabilize the numerical solution on a fixed computational
grid while preserving shock structures. One common approach involves the use of limiters,
which are widely applied in second-order finite volume methods to prevent spurious
oscillations \cite{van1979towards}. Another widely used method is Weighted Essentially
Non-Oscillatory (WENO) \cite{harten1997uniformly,jiang1996efficient,liu1994weighted}
schemes, which reconstruct high-resolution solutions using different stencils to mitigate
numerical oscillations. However, WENO is most effective on Cartesian grids, thus making it
less suitable for flows around complex geometries. Artificial viscosity is perhaps the most
popular approach in the high-order community. These approaches selectively add artificial
viscosity to elements where oscillatory behavior is detected using a shock
sensor \cite{persson2006sub,fernandez2018physics,barter2010shock,ching2019shock}.
For high Mach flows, greater levels of artificial viscosity are required, which introduces
the challenge of balancing stabilization with solution accuracy. Excessive artificial
viscosity can alter shock structure, affecting boundary layer characteristics and
potentially leading to inaccurate predictions of aerodynamic
heating in the hypersonic regime \cite{candler2015advances,lee1999spurious}.

On the other hand, shock-fitting methods align the computational mesh with solution
discontinuities to represent them perfectly without additional stabilization.
Such techniques can generate highly accurate solutions on coarse meshes, especially
when combined with high-order numerical methods. Traditional shock-fitting methods,
use the Rankine--Hugoniot conditions to compute the shock speed and the downstream
state, and subsequently determine the shock motion \cite{moretti2002thirty,salas2009shock,ciallella2020extrapolated,bonfiglioli2016unsteady,geisenhofer2020extended}.
A major challenge of this approach is the need to generate a fitted mesh that
conforms to the shock surface, which becomes increasingly difficult as the
complexity of the shock topologies increase, making it unsuitable for some complex
flows \cite{johnsen2010assessment}. Efforts to extend shock-fitting methods
to viscous flows, where shock waves are rapid transitions\footnote{Throughout this work, we define a \textit{transition} as a continuous flow feature bridging distinct states. The strength of a transition is proportional to the magnitude difference between the two states and inversely proportional to the distance over which the transition occurs.} over a finite thickness,
approximate shocks as discontinuities, thereby neglecting their thickness and internal
structure \cite{assonitis2021numerical,moretti1970numerical}.
Moretti \cite{moretti1970numerical} argues that this is the only viable approach to
handle high Reynolds flows; however, the internal shock structure is governed by
complex thermochemical effects and neglecting them could lead to inaccurate predictions
of shock-boundary-layer interactions, aerodynamic heating, and post-shock thermodynamic
states, particularly for low Reynolds number flows.

A new class of high-order methods, implicit shock tracking, that includes
moving discontinuous Galerkin method with interface condition enforcement
(MDG-ICE) \cite{corrigan2019moving,kercher2021moving,kercher2021least,ching2024moving}
and high-order implicit shock tracking (HOIST) method \cite{zahr2018shktrk,zahr2020implicit,huang2022robust,shi2022implicit,naudet2024space,huang2023high} have
overcome the key limitation of traditional shock-fitting methods by discretizing
the conservation law on a shock-agnostic mesh and simultaneously computing
the discrete flow field and coordinates of the mesh nodes as the solution
of an optimization problem. This key innovation has led to general approaches
that is equation- and problem-independent, capable of aligning grids with
complex shock structures, and recovering the design-order of the numerical
discretization for non-smooth flows. The MDG-ICE method, originally developed
for inviscid flows where shocks are discontinuities, has been extended to
an aggressive $r$-adaptive method capable of resolving the internal structure
of viscous shocks and boundary layers by compressing elements into the
transition \cite{kercher2021moving,kercher2021least,ching2024moving}.
By using a sufficiently large polynomial degree for both the solution
and mesh, these features are accurately resolved with a single, highly
stretched element. The viscous variant of the MDG-ICE method uses a
novel variational formulation that write conservation, the jump condition,
and the viscous constitutive law in their strong, first-order (flux)
form. They introduce independent test variables for each equation, which
leads to four residuals that are simultaneously minimized. This approach
is elegant, directly connected to the equations governing shock-dominated
viscous flow, and has been shown to accurately resolve challenging viscous
flows in two- and three-dimensions. However, it suffers from two drawbacks.
First, the flux formulation requires simultaneous computation of the flow
variables and their gradient, which increases the number of discrete unknowns
by a factor of $d+1$ in $d$ dimensions, which can make the linearized system
prohibitively large. Second, the viscous MDG-ICE method requires high-degree 
(usually $p \geq 4$) polynomials to
accurately approximate the internal structure of shocks and boundary layers
with a single element, but high-degree polynomials are unnecessary in the regions 
away from these viscous transitions where 
flow field is constant or slowly varying.

In this work, we introduce an $rp$-adaptive discretization as an extension
of the HOIST method that addresses the aforementioned limitations of the
viscous MDG-ICE method. First, our method uses the primal form of a
standard interior penalty discontinuous Galerkin (IPDG) \cite{arnold1982interior,arnold2002unified,hartmann2008optimal,hartmann2014higher}
to eliminate auxiliary variables prior to discretization, whiche
amounts to a savings factor of $d+1$ in terms of degrees of freedom.
The implicit shock tracking problem is formulated as an optimization
problem over the discrete flow variables and mesh nodal coordinates.
The optimization problem minimizes the IPDG residual with an enriched
test space and it is constrained by the IPDG equations to ensure the
properties of IPDG (stability, high-order accuracy, consistency, adjoint
consistency) \cite{arnold2002unified,hartmann2008optimal} are invherited.
To ensure both the shock wave and boundary layer are given similar priority
by the minimum-residual method, despite the shock wave being a stronger feature
(in terms of the magnitude of the transition and transition width), the residual
contribution of all elements touching a viscous wall are amplified.
The optimization problem is embedded in a viscosity continuation loop,
similar to that in \cite{kercher2021moving},
which proved critical for robust element compression
into viscous transitions and to avoid carbuncles.
On its own, the proposed framework functions as an aggressive
$r$-adaptive method in primal form to accurately resolve
viscous shock-dominated flows, where nodes are adjusted
during the nonlinear solve (as opposed to the more common
approach of moving nodes \textit{a posteriori} after the
solve is complete). To avoid the need for excessive
mesh compression or a large constant polynomial degree across the
domain, the entire framework is further embedded in a $p$-adaptivity
setting to iteratively increase the polynomial degree in underresolved
regions. Experimentation with several error indicators revealed
the enriched residual, which defines the shock tracking objective
function, is the most effective because it rapidly increased the
polynomial degree inside the viscous features so they can be
resolved with one or two elements, then focused the $p$-refinement
to other underresolved features. Finally, the underlying optimization
problem is solved using the sequential quadratic programming method
proposed in \cite{zahr_implicit_2020, huang2022robust} with several
robustness enhancements required to compress elements to a small
fraction of their original size. These robustness measures include
step length constraints, step modifications, and a new regularization
strategy for the Levenberg--Marquardt Hessian approximation.

The remainder of this paper is organized as follows. Section~\ref{sec:govern} introduces
the governing system of compressible viscous conservation law, its reformulation on a
fixed reference domain, and its discretization using an IPDG method; the section also
reviews solution-based error estimation and
$p$-adaptivity. Section~\ref{sec:ist} extends the implicit shock tracking
formulation originally proposed in \cite{zahr2018shktrk,zahr_implicit_2020,huang2022robust}
to incorporate continuation and $p$-adaptivity. It also introduces several robustness
enhancements to the optimization solver in \cite{huang2022robust} necessary to robustly
compress mesh elements into thin shock waves and boundary layers. Section~\ref{sec:numexp}
presents a series of numerical experiments that demonstrate the capability of the proposed
$rp$-adaptive to accurately resolve shock-dominated viscous flows, including hypersonic
flow over a cylinder with accurate heat flux profile predictions. It also systematically
studies several of the solver innovations proposed in this work to directly demonstrate
their impact. Section~\ref{sec:conclude} offers conclusions.

\section{Governing equations and high-order discretization}
\label{sec:govern}
In this section, we introduce the governing equations (steady or space-time, viscous conservation law) (Section \ref{subsec:claw}) and recast them over a fixed reference domain such that domain deformations appear explicitly (Section \ref{subsec:transfdom}). The transformed conservation law is discretized using a high-order IPDG method
(Section \ref{subsec:ipdg}). We close the section with a review of standard
\textit{a posteriori} error indicators (Section \ref{subsec:errindcator}).

\subsection{System of viscous conservation laws}
\label{subsec:claw}
A general system of $m$ steady viscous conservation laws, defined in spatial domain
$\Omega \subset \Rbb^d$ takes the form 
\begin{equation} \label{eqn:vclaw-phys0}
 \nabla \cdot F(U, \nabla U) = S(U, \nabla U) \quad \text{in}~~ \Omega,
\end{equation}
where $\func{U}{\Omega}{\Rbb^m}$, $U : x \mapsto U(x)$ is the solution of the system of
conservation laws, $\func{F}{\Rbb^m \times \Rbb^{m \times d}}{\Rbb^{m\times d}}$
is the flux function,  and $\func{S}{\Rbb^m\times\Rbb^{m\times d}}{\Rbb^{m}}$ is
the source term, and $\ds{\nabla \coloneqq (\partial_{x_1},\dots,\partial_{x_d})}$
is the gradient operator in the physical domain. The formulation of the conservation
law in (\ref{eqn:vclaw-phys0}) is sufficiently general to encapsulate
steady conservation laws in a $d$-dimensional spatial domain or
time-dependent conservation laws in a $(d-1)$-dimensional domain, i.e.,
a $d$-dimensional space-time domain.

We assume the flux function can be split into an inviscid and viscous term as
\begin{equation}
	F(W, Q) =  F_{\text{inv}}(W) - F_{\text{visc}}(W, Q),
\end{equation}
for any $W\in\Rbb^m$ and $Q\in\Rbb^{m\times d}$,
where $\func{F_\text{inv}}{\Rbb^m}{\Rbb^{m\times d}}$,
$F_\text{inv} : W \mapsto F_\text{inv}(W)$ is the inviscid flux function
and $\func{F_\text{visc}}{\Rbb^m\times\Rbb^{m\times d}}{\Rbb^{m\times d}}$,
$F_\text{visc} : (W, Q) \mapsto F_\text{visc}(W, Q)$ is
the viscous flux function, which reduces \eqref{eqn:vclaw-phys0} to 
\begin{equation}\label{eqn:vclaw-phys1}
\nabla \cdot F_{\text{inv}}(U) - \nabla \cdot F_{\text{visc}}(U,\nabla U) = S(U,\nabla U). 
\end{equation}
Additionally, we assume the viscous term is linear in its second argument
\begin{equation}
	F_{\text{visc}}(W, Q) = D(W) : Q,
\end{equation}
for any $W\in\Rbb^m$ and $Q\in\Rbb^{m\times d}$, where
$\func{D}{\Rbb^m}{\Rbb^{m \times d \times m \times d}}$ is the viscous tensor.

\subsection{Transformed conservation laws on a fixed reference domain}
\label{subsec:transfdom}
Before introducing a discretization of \eqref{eqn:vclaw-phys1}, it is convenient to
transform domain of the conservation law $\Omega$ to a fixed reference domain
$\Omega_0\subset\Rbb^d$. Let $\Gbb$ be the collection of diffeomorphisms from
the reference domain to the physical domain; i.e., for any
$\Gcal\in\Gbb$, $\func{\Gcal}{\Omega_0}{\Omega}$.
Then the conservation law on the physical domain is transformed to
a conservation law on the reference domain as
\begin{equation} \label{eqn:vclaw-ref}
\bar\nabla\cdot \bar{F}_\text{inv}(\bar{U};G) - \bar\nabla\cdot \bar{F}_\text{visc}(\bar{U}, \bar\nabla \bar{U};G) = \bar{S}(\bar{U}, \bar\nabla \bar{U};G) \quad \text{in}~~ \Omega_0,
\end{equation}
where $\func{\bar{U}}{\Omega_0}{\Rbb^m}$, $\bar{U} : X \mapsto \bar{U}(X)$ is the solution
of the transformed conservation law,
$\func{\bar{F}_\text{inv}}{\Rbb^m\times\Rbb^{d \times d}}{\Rbb^{m\times d}}$
is the  transformed inviscid flux function, $\func{\bar{F}_\text{visc}}{\Rbb^m \times \Rbb^{m \times d} \times \Rbb^{d \times d}}{\Rbb^{m\times d}}$ is the transformed viscous flux
function,
$\func{\bar{S}}{\Rbb^m \times \Rbb^{m \times d} \times \Rbb^{d\times d}}{\Rbb^{m}}$ 
is the transformed source term, and
$\ds{\bar\nabla \coloneqq (\partial_{X_1},\dots,\partial_{X_d})}$
is the gradient operator in the reference domain. The deformation gradient
$\func{G}{\Omega_0}{\Rbb^{d \times d}}$ is defined as
\begin{equation}
	G : X \mapsto \bar\nabla \Gcal(X).
\end{equation}
The reference domain solution is related to the solution over the physical domain as
\begin{equation}
	\bar{U}(X) = U(\Gcal(X)), \quad
	\bar{\nabla} \bar{U}(X) = \nabla U(\Gcal(X)) G(X).
\end{equation}
for any $X \in \Omega_0$. The fluxes and source terms transform as 
\begin{equation} \label{eqn:transf}
	\begin{aligned}
		\bar{F}_\text{inv} : (\bar{W}; \Theta) &\mapsto
		\det(\Theta) F_\text{inv}(\bar{W}) \Theta^{-T} \\
		\bar{F}_\text{visc} : (\bar{W}, \bar{Q}; \Theta) &=
		\det(\Theta) F_\text{visc}(\bar{W}, \bar{Q}\Theta^{-1}) \Theta^{-T} \\
		\bar{S} : (\bar{W}, \bar{Q}; \Theta) &=
		\det(\Theta) S(\bar{W}, \bar{Q}\Theta^{-1}).
	\end{aligned}
\end{equation}
The transformed viscous flux maintains linearity with respect to the state gradient; i.e.,
\begin{equation} \label{eqn:vflx-ref}
	\bar{F}_\text{visc}(\bar{W}, \bar{Q}; \Theta) = \bar{D}(\bar{W}; \Theta) : \bar{Q},
\end{equation}
for any $\bar{W}\in\Rbb^m$, $\bar{Q}\in\Rbb^{m\times d}$, $\Theta\in\Rbb^{d\times d}$,
where the transformed viscous tensor is
\begin{equation}
	\bar{D}_{ijps} : (\bar{W}; \Theta) \mapsto
	\det(\Theta) D(\bar{W})_{ikpl} \Theta_{sl}^{-1} \Theta_{jk}^{-1}.
\end{equation}

\subsection{Interior penalty discontinuous Galerkin discretization of transformed conservation law}
\label{subsec:ipdg}
We use an IPDG
method \cite{arnold1982interior,arnold2002unified,hartmann2008optimal,hartmann2014higher}
to discretize the transformed conservation law (\ref{eqn:vclaw-ref}).
Let $\Ecal_{h,q}$ represent a discretization of the reference domain $\Omega_0$
into non-overlapping computational elements, where $h$ is the mesh size parameter
and each element $K\in\Ecal_{h,q}$ is generated by a $q$-degree polynomial mapping
applied to a parent element $\Omega_\square\subset\Rbb^d$; i.e.,
$K = \Qcal_q^K(\Omega_\square)$ where $\Qcal_q^K \in [\Pcal_q(\Omega_\square)]^d$
and $\Pcal_a(\Scal)$ is the space of polynomial functions of degree at most $a$
over the domain $\Scal\subset\Rbb^d$.  To establish the finite-dimensional IPDG
formulation, we first convert the second-order transformed conservation law into
a system of first-order equations
\begin{equation} \label{eqn:vclaw-ref-sys}
\bar{\nabla} \cdot \bar{F}_\text{inv} (\bar{U}; G) - \bar{\nabla} \cdot \bar{\sigma} = \bar{S}(\bar{U}, \bar{\nabla}\bar{U}; G), \qquad 
    \bar{\sigma} = \bar{D}(\bar{U}; G) : \bar{\nabla} \bar{U},
\end{equation}
where $\bar\sigma : \Omega_0 \rightarrow \Rbb^{m\times d}$,
$\bar\sigma : X \mapsto \bar\sigma(X)$ is an auxiliary variable that
represents the transformed viscous flux. Next, we introduce two
finite-dimensional discontinuous piecewise polynomial approximation spaces
\begin{equation} \label{eqn:trial}
	\begin{aligned}
		\Vcal_{h,q,s} & = \left\{v \in [L^2(\Omega_0)]^m \suchthat
		 \left.v\right|_K \circ [\Qcal_q^K]^{-1} \in
		 [\Pcal_{p(K)+s}(\Omega_\square)]^m,~\forall K\in\Ecal_{h,q}\right\}, \\
		\Tcal_{h,q,s} &= \left\{v \in [L^2(\Omega_0)]^{m\times d} \suchthat
		 \left.v\right|_K \circ [\Qcal_q^K]^{-1} \in
		 [\Pcal_{p(K)+s}(\Omega_\square)]^{m\times d},
		 ~\forall K\in\Ecal_{h,q}\right\},
	\end{aligned}
\end{equation}
where $\Vcal_{h,q,s}$ is the trial space for the flow solution ($\bar{U}$),
$\Tcal_{h,q,s}$ is the trial space for the viscous flux ($\bar\sigma$),
and $p : \Ecal_{h,q} \rightarrow \Nbb$ is the distribution
of polynomial approximation degree over the mesh.
Next, we define the space of globally continuous piecewise polynomials
of degree $q$ associated with the mesh $\Ecal_{h,q}$ as
\begin{equation} \label{eqn:gfcnsp}
	\Wcal_{h,q} = \left\{v \in H^1(\Omega_0) \suchthat
	 \left.v\right|_K \circ [\Qcal_q^K]^{-1} \in
	 \Pcal_q(\Omega_\square),~\forall K\in\Ecal_{h,q}\right\},
\end{equation}
and discretize the domain mapping $\Gcal$ with the corresponding vector-valued
space $[\Wcal_{h,q}]^d$. 
From these definitions, the primal weighted residual of the IPDG method is
$\func{r_{h,q,\delta,\Delta}(\,\cdot\,,\,\cdot\,;\,\cdot\,,\rho)}{\Vcal_{h,q,\Delta} \times \Vcal_{h,q,\delta} \times [\Wcal_{h,q}]^d}{\Rbb}$ with
\begin{equation}
        r_{h,q,\delta,\Delta} :
        (\bar\psi_\mathrm{h}, \bar{W}_\mathrm{h}; \Gcal_\mathrm{h}, \rho) \mapsto
        \sum_{K\in\Ecal_{h,q}} \rho(K)  r_{h,q,\delta,\Delta}^K(\bar\psi_\mathrm{h},\bar{W}_\mathrm{h},\Gcal_\mathrm{h}),
\end{equation}
given the trial enrichment $\delta \geq 0$, test enrichment $\Delta \geq 0$,
trial function $\bar{W}_\mathrm{h}\in\Vcal_{h,q,\delta}$,
test function $\bar\psi_\mathrm{h}\in\Vcal_{h,q,\Delta}$,
mapping $\Gcal_\mathrm{h}\in\Wcal_{h,q}$, and elemental scaling function
$\rho : \Ecal_{h,q} \rightarrow \Rbb_{\ge 0}$ that will be used in
Section~\ref{subsubsec:resscale} to emphasize boundary layer features.
The elemental residual,
$\func{r_{h,q,\delta,\Delta}^K}{\Vcal_{h,q,\Delta} \times \Vcal_{h,q,\delta} \times [\Wcal_{h,q}]^d}{\Rbb}$, is given by
\begin{equation} \label{eqn:ipdg}
\begin{aligned}
	r_{h,q,\delta,\Delta}^K : (\bar\psi_\mathrm{h}, \bar{W}_\mathrm{h},\Gcal_\mathrm{h}) \mapsto
	& \int_{\partial K} \bar\psi_\mathrm{h}^+ \cdot \left(\bar\Hcal(\bar{W}_\mathrm{h}^+,\bar{W}_\mathrm{h}^-,N_\mathrm{h};\bar\nabla\Gcal_\mathrm{h}) - \hat{\sigma}(\bar{W}_\mathrm{h}^+,\bar{W}_\mathrm{h}^-, \bar\nabla \bar{W}_\mathrm{h}^+, \bar\nabla \bar{W}_\mathrm{h}^-; \bar\nabla\Gcal_\mathrm{h} )N_\mathrm{h}\right) \, dS \\
	& - \frac{1}{2}\int_{\partial K}  \left(\nabla\bar\psi_\mathrm{h}^+ : \bar{D}(\bar{W}_\mathrm{h}^+;\bar\nabla\Gcal_\mathrm{h})\right)  :  \left[(\bar{W}_\mathrm{h}^+ - \bar{W}_\mathrm{h}^-)\otimes N_\mathrm{h}\right] \, dS \\
    & - \int_{K} \bar\nabla \bar\psi_\mathrm{h} : \left ( \bar{F}_\text{inv}(\bar{W}_\mathrm{h}; \bar\nabla\Gcal_\mathrm{h}) - \bar{D}(\bar{W}_\mathrm{h};\bar\nabla\Gcal_\mathrm{h}) : \bar\nabla\bar{W}_\mathrm{h} \right ) \, dV \\
    & - \int_K \bar\psi_\mathrm{h} \cdot \bar{S}(\bar{W}_\mathrm{h}, \bar\nabla\bar{W}_\mathrm{h}; \bar\nabla\Gcal_\mathrm{h}) \, dV.
\end{aligned}
\end{equation}
where $\func{N_\mathrm{h}}{\partial K}{\Rbb^d}$ is the unit outward normal to
element $K \in \Ecal_{h,q}$, and $\bar{W}_\mathrm{h}^+(\bar{W}_\mathrm{h}^-)$
denotes the interior (exterior) trace of $\bar{W}_\mathrm{h} \in \Vcal_{h,q,\delta}$
to the element, and
$\func{\bar\Hcal}{\Rbb^m \times \Rbb^m \times \Rbb^d \times \Rbb^{d \times d}}{\Rbb^m}$
is the transformed inviscid numerical flux function \cite{zahr_implicit_2020}.
The transformed viscous numerical flux function,
$\hat\sigma : \Rbb^m \times \Rbb^m \times \Rbb^{m\times d} \times \Rbb^{m\times d} \times \Rbb^{d\times d} \rightarrow \Rbb^{m\times d}$,
is defined using the approach in \cite{hartmann2008optimal}
\begin{equation}
	\hat\sigma : (\bar{W}_\mathrm{h}^+,\bar{W}_\mathrm{h}^-,\bar\nabla\bar{W}_\mathrm{h}^+,\bar\nabla\bar{W}_\mathrm{h}^-; \bar\nabla\Gcal_\mathrm{h}) \rightarrow \left\{\!\!\left\{\bar{D}(\bar{W}_\mathrm{h}; \bar\nabla\Gcal_\mathrm{h})\bar\nabla\bar{W}_\mathrm{h})\right\}\!\!\right\} - \delta_\mathrm{IP}(\bar{W}_\mathrm{h}^+, \bar{W}_\mathrm{h}^-; \bar\nabla\Gcal_\mathrm{h}),
\end{equation}
where $\{\!\!\{v\}\!\!\} = (v^+ + v^-)/2$ is the average operator and
$\delta_\mathrm{IP} : \Rbb^m \times \Rbb^m \times \Rbb^{d\times d} \rightarrow \Rbb^{m\times d}$ is
the penalty term
\begin{equation}
	\left.\delta_\mathrm{IP}\right|_{\partial K \cap \partial K'} : (\bar{W}_\mathrm{h}^+, \bar{W}_\mathrm{h}^-; \bar\nabla\Gcal_\mathrm{h}) \mapsto C_\mathrm{IP}\frac{\max\{p(K), p(K')\}+\delta}{h} \{\!\!\{\bar{D}(\bar{W}_\mathrm{h}; \bar\nabla\Gcal_\mathrm{h})\}\!\!\}: \left[(\bar{W}_h^+ - \bar{W}_h^-)\otimes N_h\right],
\end{equation}
with $h = \min\{|K|, |K'|\} / |\partial K \cap \partial K'|$ being the linear mesh size
associated with the face $\overline K \cap \overline K'$ and $N_h$ is the outward unit normal to element $K$,
and $C_\mathrm{IP} \in \Rbb_{> 0}$ is a user-defined stabilization parameter.
On boundary faces, i.e., $\partial K \cap \partial\Omega_0$, the inviscid and viscous
numerical fluxes are modified to enforce the appropriate boundary condition.

Finally, we introduce a basis for the test space ($\Vcal_{h,q,\Delta}$), trial space
($\Vcal_{h,q,\delta}$), and domain mapping space ($[\Wcal_{h,q}]^d$) to reduce the weak
formulation in residual form to a system of nonlinear algebraic equations in
residual form. In the case where $\delta=\Delta = 0$ and $\rho = \mathbbm{1}$
is the one function, $\mathbbm{1} : \Ecal_{h,q} \rightarrow \{1\}$ with
$\mathbbm{1} : K \mapsto 1$, we denote the algebraic residual
\begin{equation} \label{eqn:ipdg-alg}
	\rbm : \Rbb^{N_\ubm}\times\Rbb^{N_\xbm} \rightarrow \Rbb^{N_\ubm}, \qquad
	\rbm : (\ubm,\xbm) \mapsto \rbm(\ubm,\xbm),
\end{equation}
where $N_\ubm = \dim\Vcal_{h,q,0}$ and $N_\xbm = \dim([\Wcal_{h,q}]^d)$,
$\ubm\in\Rbb^{N_\ubm}$ are the coefficients of the flow solution ($\bar{U}_\mathrm{h}$),
and $\xbm\in\Rbb^{N_\xbm}$ are the coefficients of the domain mapping ($\Gcal_\mathrm{h}$).
Assuming a nodal basis is used for the continuous space $\Wcal_{h,q}$, the
coefficients $\xbm$ are equal to nodal coordinates in the physical domain.
In this notation, a standard IPDG discretization (algebraic form) reads:
given $\xbm\in\Rbb^{N_\xbm}$ (from mesh generation), find $\ubm\in\Rbb^{N_\ubm}$
such that $\rbm(\ubm,\xbm) = \zerobold$. In the case where $\delta=\Delta=1$
and $\rho = \mathbbm{1}$, we denote the algebraic residual
\begin{equation} \label{eqn:ipdg-plus1-alg}
	\tilde\rbm:\Rbb^{\tilde{N}_\ubm}\times\Rbb^{N_\xbm}\rightarrow\Rbb^{\tilde{N}_\ubm},
	\qquad
        \tilde\rbm : (\tilde\ubm,\xbm) \mapsto \tilde\rbm(\tilde\ubm,\xbm),
\end{equation}
where $\tilde{N}_\ubm = \dim\Vcal_{h,q,1}$ and $\tilde\ubm\in\Rbb^{\tilde{N}_\ubm}$
are the coefficients of the flow solution in the finer trial space ($\delta = 1$).
The fine-space residual $\tilde\rbm$ is the result of a single level of global
$p$-refinement applied to the coarse-space residual $\rbm$, and will be used
to drive dual-weighted residual-based $p$-adaptivity. Finally, in the case where
$\delta = 0$ and $\Delta > 0$, we denote the algebraic residual
\begin{equation} \label{eqn:ipdg-enr-alg}
	\Rbm_\rho : \Rbb^{N_\ubm}\times\Rbb^{N_\xbm} \rightarrow \Rbb^{\hat{N}_\ubm},
	\qquad
        \Rbm_\rho : (\ubm,\xbm) \mapsto \Rbm_\rho(\ubm,\xbm),
\end{equation}
where $\hat{N}_\ubm = \dim\Vcal_{h,q,\Delta}$. For brevity,
we denote the special case $\rho = \mathbbm{1}$, as
\begin{equation} \label{eqn:ipdg-enr-alg2}
	\Rbm : \Rbb^{N_\ubm}\times\Rbb^{N_\xbm} \rightarrow \Rbb^{\hat{N}_\ubm},
	\qquad
	\Rbm : (\ubm,\xbm) \mapsto \Rbm_\mathbbm{1}(\ubm,\xbm).
\end{equation}
This enriched residual has a larger test space than trial space, and will be used to
formulate the implicit shock tracking method (Section~\ref{sec:ist}).

\begin{remark} \label{rem:elemscale}
	Due to the discontinuous test function, the IPDG method requires each
	element residual be zero. Therefore, the elemental scaling function
	$\rho$ does not alter the solution provided it is non-zero in all elements.
	As such, typically it is taken to be unity, i.e.,  $\rho = \mathbbm{1}$.
	However, $\rho$ changes the relative scaling of the element residuals,
	which alters local minima of the enriched residual. Therefore, $\rho$
	plays an important role in implicit shock tracking, which minimizes
	the enriched residual to compress elements into steep features
	(Section~\ref{sec:ist}). We introduce a simple strategy to define $\rho$
	to emphasize boundary layer features (Section~\ref{subsubsec:resscale})
	and demonstrate its importance (Section~\ref{subsubsec:bow}).
\end{remark}

\begin{remark} \label{rem:homtens}
	The transformed viscous tensor $\bar{D}$ is required by
	the IPDG method in (\ref{eqn:ipdg}). However, the viscous tensor only
	appears in contraction operations with either a trial or test function.
	By appealing to (\ref{eqn:vflx-ref}), any contraction involving the transformed
	viscous tensor can be written as an evaluation of the transformed
	viscous flux ($\bar{F}_\text{visc}$), which can in turn be written
	as an evaluation of the physical viscous flux ($F_\text{visc}$) using
	(\ref{eqn:transf}). This approach is convenient as it only requires a
	single, efficient implementation of the viscous flux in (\ref{eqn:vclaw-phys0})
	and transformations in (\ref{eqn:transf}) to implement the IPDG scheme.
\end{remark}

\begin{remark} \label{rem:liftgrad}
	There are two obvious ways to compute quantities of interest that depend
	on the viscous flux $F_\mathrm{visc}(\bar{U}_\mathrm{h}, \nabla\bar{U}_\mathrm{h}; \bar\nabla\Gcal_\mathrm{h})$, where
	$\bar{U}_\mathrm{h}\in\Vcal_{h,q,\delta}$ satisfies
	$r_{h,q,0,0}(\bar\psi_\mathrm{h}, \bar{U}_\mathrm{h}; \Gcal_\mathrm{h}, \rho) = 0$
	for all $\bar\psi_\mathrm{h}\in\Vcal_{h,q,0}$ given $\Gcal_\mathrm{h}$ and $\rho$.
	The first option is, for any $X \in \Omega_0$,
	to evaluate $F_\mathrm{visc}(\bar{U}_\mathrm{h}(X), \nabla\bar{U}_\mathrm{h}(X); \bar\nabla\Gcal_\mathrm{h}(X))$,
	where $\nabla\bar{U}_\mathrm{h}(X)$ is computed by direct differentiation
	of the polynomial basis used to represent $\Vcal_{h,q,\delta}$ over each element
	(using the same coefficients defining $\bar{U}$). The other option appeals to
	the first-order form of the conservation law (\ref{eqn:vclaw-ref-sys}) and uses
	the viscous flux that is consistent with the IPDG discretization, i.e., given the
	solution $\bar{U}_\mathrm{h}\in\Vcal_{h,q,\delta}$, find
	$\bar\sigma_\mathrm{h} \in \Tcal_{h,q,\delta}$ such that
	\begin{equation}
		\int_K \bar\tau_\mathrm{h} : (\bar\sigma_\mathrm{h} - \bar{D}(\bar{U}_\mathrm{h}) : \bar\nabla \bar{U}_\mathrm{h} )\, dV + \int_{\partial K} (\bar\tau^+ : \bar{D}(\bar{U}_\mathrm{h}^+)) : [(\bar{U}_\mathrm{h}^+ - \bar{U}_\mathrm{h}^-)\otimes N_h] \, dS = 0
	\end{equation}
	holds for all $\bar\tau_\mathrm{h} \in \Tcal_{h,q,\delta}$ and $K\in\Ecal_{h,q}$.
	Then, the viscous flux at $X \in \Omega_0$ is approximated as
	$\bar\sigma_\mathrm{h}(X)$. In Section~\ref{subsubsec:bow} we show the
	second option, i.e., the lifted viscous flux, provides more accurate and
	less oscillatory quantities of interest, which is consistent with previous
	work showing this choice leads to dual-consistent evaluation of output
	functionals and hence to superconvergence of the output
	\cite{lu2005posteriori,Hartmann_2007_Adjoint_Consistency}.
\end{remark}

\subsection{Error estimation and $p$-adaptivity}
\label{subsec:errindcator}
We close this section with a review of elemental error estimation and $p$-adaptivity,
which will be used in Section~\ref{sec:ist} to formulate the proposed $rp$-adaptive
implicit shock tracking method. First, we introduce a sequence of polynomial degree
distributions, $\{p_j\}_{j=0}^n$ with $p_j : \Ecal_{h,q} \rightarrow \Nbb$.
The function spaces in (\ref{eqn:trial}) based on this sequence
of polynomial distributions are denoted $\Vcal_{h,q,s}^{(j)}$ and
$\Tcal_{h,q,s}^{(j)}$, respectively.
Furthermore, the algebraic residuals in (\ref{eqn:ipdg-alg})--(\ref{eqn:ipdg-enr-alg2})
based on this sequence of polynomial distributions are similarly denoted
$\rbm^{(j)} : \Rbb^{N_\ubm^{(j)}}\times\Rbb^{N_\xbm} \rightarrow \Rbb^{N_\ubm^{(j)}}$,
$\tilde\rbm^{(j)} : \Rbb^{\tilde{N}_\ubm^{(j)}}\times\Rbb^{N_\xbm} \rightarrow \Rbb^{\tilde{N}_\ubm^{(j)}}$,
$\Rbm_\rho^{(j)} : \Rbb^{N_\ubm^{(j)}}\times\Rbb^{N_\xbm} \rightarrow \Rbb^{\hat{N}_\ubm^{(j)}}$, and
$\Rbm^{(j)} : \Rbb^{N_\ubm^{(j)}}\times\Rbb^{N_\xbm} \rightarrow \Rbb^{\hat{N}_\ubm^{(j)}}$,
respectively, where $N_\ubm^{(j)} = \dim\Vcal_{h,q,0}^{(j)}$,
$\tilde{N}_\ubm^{(j)} = \dim\tilde\Vcal_{h,q,1}^{(j)}$,
and
$\hat{N}_\ubm^{(j)} = \dim\Vcal_{h,q,\Delta}^{(j)}$.

The initial polynomial degree distribution, $p_0$, is usually chosen to be constant
across the mesh. Each $p$-adaptive iteration locally increases the polynomial
degree from $p_j$ to $p_{j+1}$ based on the distribution of an elemental
error indicator $s_j:\Ecal_{h,q}\rightarrow \Rbb_{\geq 0}$ as
\begin{equation} \label{eqn:pref}
	p_{j+1} : K \mapsto
        \begin{dcases}
		p_j(K) + 1, & s_j(K) > \tau\max_{K'\in\Ecal_{h,q}} s_j(K') \\
		p_j(K),   & \text{otherwise},
        \end{dcases}
\end{equation}
where $\tau > 0$ is a threshold used to identify which elements to refine.
The error indicator is a constant scalar over each element of the mesh $\Ecal_{h,q}$
that approximates some measure of the error between the discrete and true ($\bar{U}$)
attributed to the element. It is used to identify the elements with the largest
contribution to the overall error for $p$-refinement. In the following sections,
we review three popular error indicators studied in the context of the proposed
$rp$-adaptive method.

\subsubsection{Enriched residual indicator}
\label{subsubsec:enrres}
First, we consider the unweighted elemental enriched residual,
$\hat{s}_j^\mathrm{uwr} : \Ecal_{h,q} \times \Rbb^{N_\ubm^{(j)}} \times \Rbb^{N_\xbm} \rightarrow \Rbb_{\geq 0}$, defined as
\begin{equation}
        \hat{s}_j^\mathrm{uwr} :
	(K, \ubm, \xbm) \mapsto \norm{(\hat\Pbm_K^{(j)})^T\Rbm^{(j)}(\ubm,\xbm)}_2,
\end{equation}
where 
$\hat\Pbm_K^{(j)}\in\Rbb^{\hat{N}_\ubm^{(j)}\times \hat{N}_\ubm^{(j),K}}$
is the assembly operator that maps the local degrees of
freedom associated with element $K\in\Ecal_{h,q}$ to
the corresponding global degrees of freedom for the
test space $\hat\Vcal_K^{(j)}$, and the transpose
restricts global degrees of freedom to the degrees
of freedom associated with element $K$.
This error indicator is less popular than dual-weighted variants \cite{Becker_1996_Feed_Error_Control,Becker_2001_DWR,fidkowski2011review}
(Section~\ref{subsubsec:dwr});
however, it is natural to use in the context of implicit shock tracking because it
is used to define the shock tracking objective function (Section~\ref{subsec:optimform}).

\subsubsection{Dual-weighted residual indicator}
\label{subsubsec:dwr}
Next, we consider dual-weighted elemental residual indicator \cite{Becker_1996_Feed_Error_Control,Becker_2001_DWR,fidkowski2011review},
$\hat{s}_j^\mathrm{dwr} : \Ecal_{h,q} \times \Rbb^{N_\ubm^{(j)}} \times \Rbb^{N_\xbm} \rightarrow \Rbb_{\geq 0}$, defined as
\begin{equation}
	\hat{s}_j^\mathrm{dwr} :
	(K,\ubm,\xbm)\mapsto
	\left|\left[(\hat\Qbm_K^{(j)})^T\Psibold^{(j)}(\ubm,\xbm)\right] \cdot (\hat\Qbm_K^{(j)})^T\Rbm^{(j)}(\ubm,\xbm)\right|
\end{equation}
where the dual variable
$\Psibold^{(j)} : \Rbb^{N_\ubm^{(j)}}\times\Rbb^{N_\xbm} \rightarrow \Rbb^{\hat{N}_\ubm^{(j)}}$
is implicitly defined as the solution of the fine-space adjoint equations, i.e.,
\begin{equation}
	\Psibold^{(j)} : (\ubm, \xbm) \mapsto
	\left[\pder{\tilde\rbm^{(j)}}{\tilde\ubm}(\tilde\Qbm^{(j)}\ubm,\xbm)\right]^{-T}
	\pder{\tilde{J}^{(j)}}{\tilde\ubm}(\tilde\Qbm^{(j)}\ubm,\xbm)^T
\end{equation}
and $\tilde\Qbm^{(j)}\in\Rbb^{\tilde{N}_\ubm^{(j)}\times N_\ubm^{(j)}}$ is a discrete
prolongation operator that maps the coarse-space ($\delta=0$) coefficients
$\ubm\in\Rbb^{N_\ubm^{(j)}}$ to the fine-space ($\delta=1$) coefficients
$\tilde\ubm\in\Rbb^{\tilde{N}_\ubm^{(j)}}$ that
reproduces the same function over $\Omega_0$. Furthermore,
$\tilde{J}^{(j)} : \Rbb^{\tilde{N}_\ubm^{(j)}}\times\Rbb^{N_\xbm}\rightarrow \Rbb$ with
$\tilde{J}^{(j)} : (\tilde\ubm,\xbm) \mapsto \tilde{J}^{(j)}(\tilde\ubm,\xbm)$ is the
algebraic form (fine-space, $\delta = 1$) of the continuous-level quantity of interest
\begin{equation}
	\Jcal : (\bar{U}, \Gcal) \mapsto
	\int_{\Gcal(\Omega_0)} q_\mathrm{v}(\bar{U}; \Gcal) \, dV +
	\int_{\Gcal(\partial\Omega_0)} q_\mathrm{b}(\bar{U}; \Gcal) \, dS,
\end{equation}
based on the trial space $\Vcal_{h,q,1}^{(j)}$,
where $q_\mathrm{v}(\bar{U},\Gcal)$ is the volumetric quantity of interest and
$q_\mathrm{b}(\bar{U},\Gcal)$ is the boundary quantity of interest.
In this work, we take $q_\mathrm{v} = 0$
and $q_\mathrm{b}$ is the drag (Section~\ref{subsubsec:fp}) or
heat flux (Section~\ref{subsubsec:bow_snsr}) over a relevant boundary.
This is a popular error indicator because the sum
of $\hat{s}_j^\mathrm{dwr}$ across all elements in the mesh provides an estimate
of the error between the true quantity of interest $\Jcal(\bar{U},\Gcal)$ and the
discrete approximate $\tilde{J}^{(j)}(\tilde\ubm,\xbm)$
\cite{venditti2003anisotropic,venditti2002grid,fidkowski2011review}.

\subsubsection{Feature-based indicator}
\label{subsubsec:fbs}
Lastly, we consider a feature-based indicator,
$\hat{s}_j^\mathrm{fbs} : \Ecal_{h,q} \times \Rbb^{N_\ubm^{(j)}} \times \Rbb^{N_\xbm} \rightarrow \Rbb_{\geq 0}$, defined as
\begin{equation}
        \hat{s}_j^\mathrm{fbs} :
	(K, \ubm, \xbm) \mapsto \int_{\Gcal(K)} \vartheta(\bar{U},\Gcal) \, dV,
\end{equation}
where $\vartheta : \Vcal_{h,q,0}^{(j)} \times [\Wcal_{h,q}]^d \rightarrow \Rbb$
is a relevant volumetric sensor. In this work, we use the gradient of the solution
vector $\vartheta(\bar{U},\Gcal) = \norm{\bar\nabla\bar{U} \cdot (\bar\nabla \Gcal)^{-1}}$.
Such sensors detect regions of significant physical change in the state and can
selectively identify shocks and boundary layers. However, the sensor does not
diminish as resolution increases, which can cause refinement to be incorrectly
focused and does not admit a natural termination criteria.

\section{A $p$-adaptive high-order implicit shock tracking method for viscous conservation laws}
\label{sec:ist}

In this section, we embed the HOIST method \cite{zahr2018shktrk,zahr2020implicit,huang2022robust} in a $p$-adaptivity setting (Section~\ref{subsec:optimform}) and recall the sequential quadratic programming solver \cite{huang2022robust} (Section~\ref{subsec:sqp}). Finally, we introduce modifications to the HOIST formulation and solver necessary to fully resolve shock and boundary layers in viscous conservation law solutions (Section~\ref{subsec:slvrimprov}) and summarize the complete algorithm (Section~\ref{subsec:summary}).

\subsection{Optimization formulation}
\label{subsec:optimform}
The HOIST method treats both the discrete solution of
the conservation law ($\ubm$) and the nodal coordinates of the mesh ($\xbm$)
as unknowns, which are determined as the solution of a PDE-constrained
optimization problem. In the viscous setting, we seek nodal coordinates
that cause mesh elements to compress into steep features such as shocks
and boundary layers, and flow variables that satisfy the IPDG equations
with $\delta = \Delta = 0$ on the tailored mesh. On the $j$th $p$-adaptivity
iteration, the optimization problem is formulated as
\begin{equation} \label{eqn:pde-opt}
 \left(\ubm_{\star}^{(j)}(\Xi), \ybm_{\star}^{(j)}(\Xi)\right) \coloneqq
 \argmin_{\ubm\in\Rbb^{N_\ubm},\ybm\in\Rbb^{N_\ybm}} f^{(j)}\left(\ubm, \ybm; \Xi \right) \quad \text{subject to:} \quad \rbm^{(j)}\left(\ubm,\phibold(\ybm); \Xi \right) = \zerobold,
\end{equation}
where $\func{f^{(j)}}{\Rbb^{N_\ubm}\times\Rbb^{N_\ybm}}{\Rbb}$ is the objective
function, $\ybm\in\Rbb^{N_\ybm}$ are the unconstrained mesh degrees of freedom,
$\phibold : \Rbb^{N_\ybm} \rightarrow \Rbb^{N_\xbm}$ is a parametrization
of the nodal coordinates that ensures $\xbm = \phibold(\ybm)$ represents a
physical mesh where all nodes remain on their original boundaries for any
$\ybm\in\Rbb^{N_\ybm}$ \cite{zahr_implicit_2020,huang2022robust,perez2025istdomain},
and $\Xi \in \Rbb$ is a continuation parameter. Continuation
(Section~\ref{subsec:cont}~and~\ref{subsubsec:viscont}) will be used as an
``outer loop'' around the optimization problem in (\ref{eqn:pde-opt})
to facilitate robust initialization and convergence. 
The $r$-adapted mesh coordinates are reconstructed from the unconstrained
degrees of freedom as $\xbm_\star^{(j)} = \phibold(\ybm_\star^{(j)}(\Xi))$.
The objective function is composed of two terms as
\begin{equation}\label{eqn:obj0}
 f^{(j)} : \left(\ubm, \ybm; \Xi\right) \mapsto f_\text{err}^{(j)}\left(\ubm, \ybm; \Xi \right) + \kappa^2 f_\text{msh}(\ybm),
\end{equation}
which balances alignment of the mesh with compression of the mesh into steep features
and the quality of mesh. The mesh alignment term,
$\func{f_\text{err}^{(j)}}{\Rbb^{N_\ubm}\times\Rbb^{N_\ybm}}{\Rbb_{>0}}$,
is taken to be the norm of enriched DG residual
\begin{equation}\label{eqn:obj1}
 f_\text{err}^{(j)} : \left(\ubm, \ybm; \Xi \right) \mapsto \frac{1}{2}\norm{\Rbm_\rho^{(j)}\left(\ubm,\phibold(\ybm); \Xi\right)}_2^2,
\end{equation}
where $\rho : \Omega_0 \rightarrow \Rbb_{\ge 0}$ is an element scaling defined
in Section~\ref{subsec:ipdg} that will be used to
emphasize the boundary layer (Section~\ref{subsubsec:resscale}). The mesh distortion term,
$f_\text{msh} : \Rbb^{N_\ybm} \rightarrow \Rbb_{>0}$, is defined as
\begin{equation}\label{eqn:obj2}
f_{\text{msh}} : \ybm \mapsto \frac{1}{2}\norm{\Rbm_{\text{msh}}(\phibold(\ybm))}_2^2,
\end{equation}
where $\Rbm_\text{msh} : \Rbb^{N_\xbm} \rightarrow \Rbb_{>0}^{|\Ecal_{h,q}|}$ is the
elemental mesh distortion defined in \cite{zahr2020implicit}, and $\kappa\in\Rbb_{>0}$
is the mesh quality penalty parameter (Section~\ref{subsubsec:mshpnlty}).
The objective function can be written as the two-norm of a residual function as
\begin{equation}
f^{(j)}(\ubm, \ybm; \Xi) = \frac{1}{2}\Fbm^{(j)}(\ubm, \ybm;\Xi)^T\Fbm^{(j)}(\ubm, \ybm;\Xi), \quad
\Fbm^{(j)} : (\ubm, \ybm;\Xi) \mapsto 
\begin{bmatrix}
  \Rbm_\rho^{(j)}\left(\ubm,\phibold(\ybm); \Xi\right) \\ \Rbm_{\text{msh}}^{(j)}(\phibold(\ybm))
\end{bmatrix},
\end{equation}
where $\Fbm^{(j)} : \Rbb^{N_\ubm^{(j)}}\times\Rbb^{N_\ybm} \rightarrow \Rbb^{\hat{N}_\ubm^{(j)}+|\Ecal_{h,q}|}$ is the residual, which allows for a Levenberg--Marquardt approximation to its Hessian (Section~\ref{subsec:sqp}).

\subsection{Sequential quadratic programming method}
\label{subsec:sqp}
To solve the system in (\ref{eqn:pde-opt}), we adapt the full-space sequential
quadratic programming (SQP) method developed in \cite{zahr_implicit_2020,huang2022robust}
that simultaneously drives
$\ubm$ and $\ybm$ to their optimal values. The SQP solver defines a sequence of
iterations $\{\zbm_0^{(j)}(\Xi), \zbm_1^{(j)}(\Xi), \dots\}$, where each iterate
contains both the flow and mesh degrees of freedom, i.e.,
$\zbm_k^{(j)}(\Xi) = (\ubm_k^{(j)}(\Xi), \ybm_k^{(j)}(\Xi))$ that converges to
the solution of (\ref{eqn:pde-opt}): $\zbm_k^{(j)}(\Xi) \rightarrow (\ubm_\star^{(j)}(\Xi), \ybm_\star^{(j)}(\Xi))$. The starting point
$\zbm_0^{(j)}(\Xi) = (\ubm_0^{(j)}(\Xi), \ybm_0^{(j)}(\Xi))$ is provided by
the user, continutation, or $p$-adaptation strategy, and the SQP update is
\begin{equation} \label{eqn:stpupd}
	\zbm_{k+1}^{(j)}(\Xi) = \chi_k^{(j)}(\tilde\zbm_{k+1}^{(j)}(\Xi);\Xi), \qquad
	\tilde\zbm_{k+1}^{(j)}(\Xi) = \zbm_{k}^{(j)}(\Xi) + \alpha_k^{(j)}(\Xi) \Delta \zbm_k^{(j)}(\Xi),
\end{equation}
where $\chi_k^{(j)}(\,\cdot\,;\Xi) : \Rbb^{N_\ubm^{(j)} + N_\ybm} \rightarrow \Rbb^{N_\ubm^{(j)} + N_\ybm}$
is a step modification function used to enhance the robustness of the solver
(Section~\ref{subsubsec:stepmod}).
The search direction, $\Delta \zbm_k^{(j)}(\Xi) \in \Rbb^{N_\ubm^{(j)} + N_\ybm}$, is
defined as the solution of the following quadratic program:
\begin{equation} \label{eqn:quadprog}
\begin{aligned}
	&\underset{\Delta \zbm \in \Rbb^{N_\ubm^{(j)} + N_\ybm}}{\text{minimize}} \quad \gbm_k^{(j)}(\Xi) \Delta\zbm + \frac{1}{2} \Delta\zbm^T \Bbm_k^{(j)}(\Xi) \Delta\zbm \\
	&\text{subject to} \quad \rbm_k^{(j)}(\Xi) + \Jbm_k^{(j)}(\Xi)\Delta\zbm=\zerobold,
\end{aligned}
\end{equation}
where $\rbm_k^{(j)}(\Xi) = \rbm^{(j)}(\ubm_k^{(j)}(\Xi), \phibold(\ybm_k^{(j)}(\Xi)))$
is the constraint function evaluated at $\zbm_k^{(j)}$, the first-order terms are defined as
\begin{equation}
	\begin{aligned}
		\gbm_k^{(j)}(\Xi) &= \begin{bmatrix} \partial_\ubm f^{(j)}(\ubm_k^{(j)}(\Xi),\ybm_k^{(j)}(\Xi); \Xi)^T \\  \partial_\ybm f^{(j)}(\ubm_k^{(j)}(\Xi),\ybm_k^{(j)}(\Xi); \Xi)^T \end{bmatrix}, \\
		\Jbm_k^{(j)}(\Xi) &= \begin{bmatrix} \partial_\ubm\rbm^{(j)}(\ubm_k^{(j)}(\Xi),\ybm_k^{(j)}(\Xi); \Xi) &  \partial_\ybm\rbm^{(j)}(\ubm_k^{(j)}(\Xi),\ybm_k^{(j)}(\Xi); \Xi) \end{bmatrix}.
	\end{aligned}
\end{equation}
The second-order term uses the following Levenberg--Marquardt Hessian approximation
\begin{equation}
	\Bbm_k^{(j)}(\Xi) =
	\begin{bmatrix}
		\Bbm_{\ubm\ubm,k}^{(j)}(\Xi) & \Bbm_{\ubm\ybm,k}^{(j)}(\Xi) \\
		\Bbm_{\ubm\ybm,k}^{(j)}(\Xi)^T & \Bbm_{\ybm\ybm,k}^{(j)}(\Xi)
	\end{bmatrix},
\end{equation}
where the individual terms are
\begin{equation}
	\begin{aligned}
		\Bbm_{\ubm\ubm,k}^{(j)}(\Xi) &= \partial_\ubm\Fbm^{(j)}(\ubm_k^{(j)}(\Xi),\ybm_k^{(j)}(\Xi); \Xi)^T\partial_\ubm\Fbm^{(j)}(\ubm_k^{(j)}(\Xi),\ybm_k^{(j)}(\Xi); \Xi), \\ 
		\Bbm_{\ubm\ybm,k}^{(j)}(\Xi) &= \partial_\ubm\Fbm^{(j)}(\ubm_k^{(j)}(\Xi),\ybm_k^{(j)}(\Xi); \Xi)^T\partial_\ybm\Fbm^{(j)}(\ubm_k^{(j)}(\Xi),\ybm_k^{(j)}(\Xi); \Xi), \\
		\Bbm_{\ybm\ybm,k}^{(j)}(\Xi) &= \partial_\ybm\Fbm^{(j)}(\ubm_k^{(j)}(\Xi),\ybm_k^{(j)}(\Xi); \Xi)^T\partial_\ybm\Fbm^{(j)}(\ubm_k^{(j)}(\Xi),\ybm_k^{(j)}(\Xi); \Xi) + \gamma_k^{(j)}(\Xi) \Dbm_k^{(j)}(\Xi),
	\end{aligned}
\end{equation}
and $\Dbm_k^{(j)}(\Xi) \in \Rbb^{N_\ybm \times N_\ybm}$ is a regularization matrix
and $\gamma_k^{(j)}(\Xi) \in \Rbb_{\geq 0}$ is a regularization parameter
(Section~\ref{subsubsec:hessreg}) applied
only to the mesh component of the Hessian approximation.
Finally, the step length, $\alpha_k^{(j)}(\Xi) \in (0, \hat\alpha_k^{(j)}]$, is obtained
from a backtracking line search \cite{nocedal2006numerical} on the $\ell_1$ merit
function of (\ref{eqn:pde-opt}) with penalty parameter chosen according
to \cite{huang2022robust}.
The maximum step length $\hat\alpha_k^{(j)}(\Xi) \in (0, 1]$ is chosen to ensure the next
iteration remains physical (Section~\ref{subsubsec:stepcon}). The quadratic program in
(\ref{eqn:quadprog}) can be recast as a linear system of equations by appealing
to the first-order optimality conditions; for additional detail see \cite{huang2022robust}.

\subsection{Continuation}
\label{subsec:cont}
The continuation parameter, $\Xi$, is used to improve the robustness of the shock
tracking solver (Section~\ref{subsec:sqp}) and avoid spurious features (such as carbuncles
\cite{peery1988blunt,robinet2000shock}) by defining a sequence of optimization problems,
each one with an reasonable starting point, ending with the optimization problem of
interest ($\Xi = \Xi_\star$).
Continuation is only required \textit{prior} to adaption of the polynomial degree ($j=0$),
i.e., when the flow and grid are far from convergence.
To this end, we introduce a sequence $\{\Xi_1,\Xi_2,\dots,\Xi_c = \Xi_\star\}$, where
$c$ is the number of continuation steps, and initialize the optimization problem
in (\ref{eqn:pde-opt}) for $\Xi = \Xi_i$ with
\begin{equation}
\ubm_0^{(0)}(\Xi_i) \coloneqq \ubm_\star^{(0)}(\Xi_{i-1}), \qquad
\ybm_0^{(0)}(\Xi_i) \coloneqq \ybm_\star^{(0)}(\Xi_{i-1})
\end{equation}
for $i=2,\dots,c$.
The initial guess for (\ref{eqn:pde-opt}) with $\Xi = \Xi_1$ comes from a
shock capturing simulation (initialization of $\ubm$) on a shock-agnostic mesh
obtained from mesh generation (initialization of $\ybm$). Once the continuation process
is complete, the solution $(\ubm_\star^{(0)}(\Xi_\star), \ybm_\star^{(0)}(\Xi_\star))$
is available and represents a purely $r$-adapted solution.

\subsection{$p$-Adaptivity}
\label{subsec:istpadapt}
Next, we introduce $p$-adaptivity to further improve the $r$-adapted approximation.
Given the solution at the final continuation parameter $\Xi = \Xi_\star$ of the $j$th
$p$-adaption iteration, i.e., $(\ubm_\star^{(j)}(\Xi_\star), \ybm_\star^{(j)}(\Xi_\star))$,
we locally increase the polynomial degree according to (\ref{eqn:pref}) using
the indicator
\begin{equation}
	s_j : K \mapsto
	\hat{s}_j(K, \ubm_\star^{(j)}(\Xi_\star), \xbm_\star^{(j)}(\Xi_\star)),
\end{equation}
where $\hat{s}_j$ is any of the indicators introduced in
Sections~\ref{subsubsec:enrres}-\ref{subsubsec:fbs}, to obtain the
new polynomial degree distribution $p_{j+1}$. The optimization
problem at the new refinement level, i.e., $j \leftarrow j + 1$
and $\Xi = \Xi_\star$, is initialized from prolongation of the
previous solution $\ubm_\star^{(j)}(\Xi_\star)$ into the finer
space to define $\ubm_0^{(j+1)}(\Xi_\star)$, and the mesh is directly
transferred, i.e., $\ybm_0^{(j+1)}(\Xi_\star) = \ybm_\star^{(j)}(\Xi_\star)$.
From this starting point, the solver in Section~\ref{subsec:sqp} is applied
to obtain $\ubm_\star^{(j+1)}(\Xi_\star)$ and $\ybm_\star^{(j+1)}(\Xi_\star)$.
The $p$-adaption iterations are terminated at iteration $j = J$, where one of
the following conditions are met
\begin{equation} \label{eqn:padapt_conv}
	\norm{\Rbm^{(J)}(\ubm_\star^{(J)}(\Xi_\star),\phibold(\ybm_\star^{(J)}(\Xi_\star)); \Xi_\star)} \leq \hat\epsilon \qquad \text{or} \qquad
	J = J_\text{max},
\end{equation}
where $\hat\epsilon > 0$ and $J_\mathrm{max} \in \Nbb$ are user-defined parameters.
The final $rp$-adapted flow solution is $\ubm_\star^{(J)}(\Xi_\star)$
on the transformed mesh defined by the nodal coordinates
$\phibold(\ybm_\star^{(J)}(\Xi_\star))$.

\subsection{Solver details tailored to viscous flows}
\label{subsec:slvrimprov}
The exposition in Section~\ref{subsec:optimform}-\ref{subsec:istpadapt}
provides an overview of the HOIST method initially introduced
in \cite{zahr_implicit_2020,huang2022robust}
enriched with continutation in \cite{huang2023high} and $p$-adaptivity
in \cite{dong2025high}. The remainder of this section completes the missing algorithmic
details necessary for the HOIST method to act as an aggressive
$r$-adaptive method to compress elements into steep viscous flow
features.

\subsubsection{Enrichment degree}
\label{subsubsec:enrdeg}
The enriched residual is used to defined the implicit shock tracking objective
function in (\ref{eqn:obj1}); however, the enrichment degree $\Delta$, defined
in (\ref{eqn:ipdg-enr-alg})-(\ref{eqn:ipdg-enr-alg2}), has been unspecified
to this point. The choice $\Delta = 1$ has been shown to lead to
accurate and robust shock tracking for inviscid flows across a number of
applications \cite{zahr_implicit_2020,zahr2020react,huang2022robust,huang2023high,naudet2024space}.
However, for viscous problems, $\Delta = 1$ leads to
suboptimal meshes and insufficient boundary layer resolution
(Section~\ref{subsubsec:bow_enrch}), especially when the solution
is severely underresolved. On the other hand, $\Delta = 2$
has led to reliable and aggressive $r$-adaptivity to both
viscous shocks and boundary layers for all problems considered
in this work. 

\subsubsection{Choice of mesh penalty parameter}
\label{subsubsec:mshpnlty}
For inviscid problems, $\kappa$ is adapted throughout the iterations to ensure
the tracking and mesh quality terms remain in balance \cite{huang2022robust}.
However, in the viscous
setting, our goal is to highly distort the elements to compress them into steep features,
e.g., viscous shocks and boundary layers, which means $f_\text{msh}(\ybm)$ will inevitably
become large. Therefore, we choose $\kappa$ to be a fixed, small, non-zero number (we take
$\kappa = 10^{-6}$ in this work). This ensures that $\kappa^2 f_\text{msh}(\ybm)$ will
blow up if an inverted mesh is encountered and cause the line search to decrease the step
length. However, for non-inverted meshes, the tracking term $f_\text{err}$ will dominate.

\subsubsection{Residual scaling to emphasize boundary layers}
\label{subsubsec:resscale}
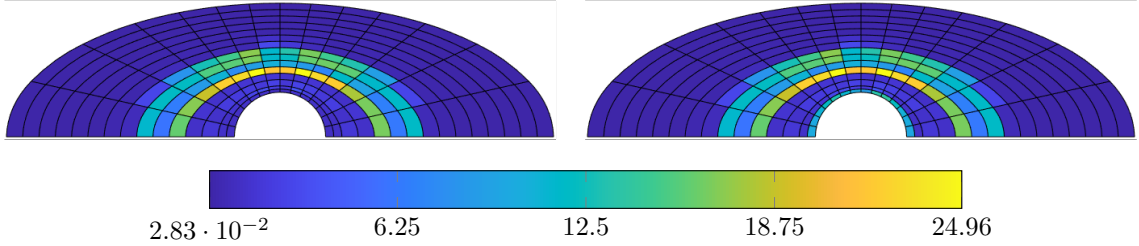
\begin{figure}
\centering
\begin{tikzpicture}[scale=0.8]
\begin{groupplot}[
  group style={
      group size=2 by 1,
      horizontal sep=0.5cm,
      vertical sep=0.5cm
  },
  width=0.65\textwidth,
  axis equal image,
  xlabel={$x_1$},
  ylabel={$x_2$},
  xtick = {-6.0, 0.0, 6.0},
  xticklabels={-6.0, 0.0, 6.0},
  ytick = {0.0, 3.0},
  xmin=-6.0, xmax=6.0,
  ymin=0, ymax=3.0
]

\nextgroupplot[xtick=\empty, xlabel={}, ytick=\empty, ylabel={}]
\addplot graphics [xmin=-6.0, xmax=6.0, ymin=0, ymax=3.0] {{_img/cyl_M5Re1e3_Enres_woscale_refdom}.png};

\nextgroupplot[xtick=\empty, xlabel={}, ytick=\empty, ylabel={}]
\addplot graphics [xmin=-6.0, xmax=6.0, ymin=0, ymax=3.0] {{_img/cyl_M5Re1e3_Enres_wscale_refdom}.png};

\end{groupplot}
\node[anchor=north] at ($(group c1r1.south)!0.5!(group c2r1.south)$) {\colorbarMatlabParula{0.02827900185167261}{6.25}{12.5}{18.75}{24.95653237639528}};
\end{tikzpicture}
	\caption{
		Elementwise magnitude of the enriched DG residual ($\Delta = 2$) for
		$M_\infty=5$, $\mathrm{Re}=10^3$ hypersonic flow over a cylinder
		(isothermal wall) without (\textit{left}) and with (\textit{right})
		boundary layer scaling ($\lambda$ = 10).
	}
 \label{fig:cyl_M5Re1e3_Enres_comp}
\end{figure}

Implicit shock tracking methods prioritize stronger
transitions over weaker ones, in the sense that many more SQP iterations
are required to track weak transitions than strong ones, owing to the residual-based
objective function \cite{huang2022robust}. For viscous flows, particularly in
the hypersonic regime, shock waves are much stronger transitions than boundary
layers. This is not typically a problem for standard DG discretizations that seek
a root of $\rbm(\cdot,\xbm)$ on a given mesh $\xbm\in\Rbb^{N_\xbm}$ because the
residual on each element is driven to zero regardless of its initial magnitude;
however, implicit shock tracking methods are minimum-residual methods whose solutions
are highly dependent on the scaling of the residual. Because the enriched-residual
objective function (\ref{eqn:obj1}) is less sensitive to weak transitions such
as boundary layers than strong transitions, the unweighted enriched residual
(\ref{eqn:ipdg-enr-alg2}) is not an appropriate objective function, despite its
success for inviscid problems \cite{zahr_implicit_2020,zahr2020react,huang2022robust,huang2023high,naudet2024space}, because it will effectively ignore the boundary
layer. Instead, we choose $\rho$ in the objective function (\ref{eqn:obj1}) to be
\begin{equation} \label{eqn:weight-scaling}
	\rho : K \mapsto
	\begin{cases}
		\lambda, & K \in \Ecal_{h,q}^\text{wall} \\
		1, & \text{otherwise},
	\end{cases}
\end{equation}
where $\lambda \geq 1$ is a boundary layer amplification factor and
$\Ecal_{h,q}^\text{wall} \subset \Ecal_{h,q}$ is the subset of the
elements of the mesh with at least one boundary touching a viscous wall.
The impact of the boundary layer scaling on the enriched residual is
demonstrated in Figure~\ref{fig:cyl_M5Re1e3_Enres_comp}. Boundary layer
scaling is critical in obtaining accurate solutions to shock-dominated
viscous flows using implicit shock tracking, particularly with regard to
prediction of heat flux profiles (Section~\ref{subsubsec:bow_bndwght}).

\subsubsection{Step length constraint}
\label{subsubsec:stepcon}
Non-physical states and oscillations are often a byproduct of the flow variables
changing too rapidly when far from the solution. To avoid this instabilities, we
employ an increment-limiting strategy that restricts the maximum step length for
the line search \cite{modisette2011automated,yano2011adaptation}. To this end, let
$\phi_k^{(j)}(\,\cdot\,; \Xi) : \Omega_0 \rightarrow \Rbb_{>0}$
denote a positive, scalar-valued function over domain defined from the SQP
iterate $\ubm_k^{(j)}(\Xi)$ and let
$\Delta\phi_k^{(j)}(\,\cdot\,; \Xi) : \Omega_0 \rightarrow \Rbb_{>0}$
be its proposed modification defined from $\Delta\zbm_k^{(j)}(\Xi)$.
With the goal of limiting the change in $\phi_k^{(j)}(\Xi)$, we define the function
\begin{equation} \label{eqn:iota}
	\iota_k^{(j)} : (X; \Xi) \mapsto
	\begin{cases}
		1, & d_k^{(j)}(X; \Xi) \leq \theta_\mathrm{l} \\
		0, & d_k^{(j)}(X; \Xi) \geq \theta_\mathrm{u} \\
		\theta_\mathrm{l}/d_k^{(j)}(X;\Xi), & d_k^{(j)}(X; \Xi) \in (\theta_\mathrm{l}, \theta_\mathrm{u}),
	\end{cases}
\end{equation}
where $d_k^{(j)}(\,\cdot\,; \Xi) : \Omega_0 \rightarrow \Rbb$ is the relative change with
$d_k^{(j)} : (X; \Xi) \mapsto |\Delta\phi_k^{(j)}(X;\Xi)| / \phi_k^{(j)}(X;\Xi)$,
$\theta_\mathrm{l}\in\Rbb_{>0}$ is the largest relative change that will be allowed
a full Newton step (unit step length),
and $\theta_\mathrm{u}\in\Rbb_{>0}$ is the smallest relative change that will be rejected
(zero step length). The minimum of the pointwise function (\ref{eqn:iota}) over the
quadrature nodes of the mesh $\Ecal_{h,q}$ defines the maximum step length, i.e.,
\begin{equation} \label{eqn:stplencon}
	\hat\alpha_k^{(j)}(\Xi) = \min_{X\in\Lambda_{h,q}}~\iota_k^{(j)}(X;\Xi),
\end{equation}
where $\Lambda_{h,q}\subset\Omega_0$ is the collection of quadrature nodes associated
with the mesh $\Ecal_{h,q}$. In this work, we apply this procedure to both the density
and pressure (using whichever step length is smaller) and take $\theta_\mathrm{l}=0.1$
and $\theta_\mathrm{u} = 10$. Assuming the current iterate $\ubm_k^{(j)}(\Xi)$ is free of
negative densities and pressures, this procedure ensures they will not be introduced
during the update from $\zbm_k^{(j)}(\Xi)$ to $\tilde\zbm_{k+1}^{(j)}(\Xi)$.

\subsubsection{Step modification}
\label{subsubsec:stepmod}
Robustness of the SQP solver (Section~\ref{subsec:sqp}) is further improved by
occasionally freezing the mesh $\xbm_k^{(j)}$ and updating only the flow variables
to satisfy the constraint of (\ref{eqn:pde-opt}). This flexibility is afforded by the step
modification function $\chi_k^{(j)}(\,\cdot\,;\Xi)$ that maps the candidate step
$\tilde\zbm_{k+1}^{(j)}(\Xi)$ to the actual step $\zbm_{k+1}^{(j)}(\Xi)$.
In practice, the constraint cannot
be satisfied until the mesh optimization is complete because we do not use
any artificial stabilization. Instead, we only require the magnitude of the
constraint residual reduces by a preset amount. To this end, let
$\check\ubm^{(j)}(\,\cdot\,,\,\cdot\,;\Xi) : \Rbb^{N_\ubm^{(j)}}\times\Rbb^{N_\xbm} \rightarrow \Rbb^{N_\ubm^{(j)}}$ be any function that satisfies
\begin{equation} \label{eqn:stpmod_cond}
	\norm{\rbm^{(j)}(\check\ubm^{(j)}(\ubm,\xbm;\Xi),\xbm;\Xi)} < \epsilon \norm{\rbm^{(j)}(\ubm,\xbm;\Xi)}
\end{equation}
for any $\ubm\in\Rbb^{N_\ubm^{(j)}}$ and $\xbm\in\Rbb^{N_\xbm}$, where
$\epsilon\in(0, 1)$ is a user-defined parameter (in this work, $\epsilon=0.9$).
Then, we define the step modification function as
\begin{equation} \label{eqn:stpmod}
	\chi_k^{(j)} : ((\ubm, \ybm); \Xi) \mapsto
	\begin{cases}
		(\check\ubm^{(j)}(\ubm,\phibold(\ybm);\Xi), \ybm), &\lfloor k/\omega_\mathrm{mod}\rfloor = k/\omega_\mathrm{mod} \\ 
		(\ubm, \ybm), &\text{otherwise},
	\end{cases}
\end{equation}
where the modification is only applied every $\omega_\mathrm{mod}$ iterations. 
We use pseudo-transient continuation
(PTC) \cite{kelley1998convergence,modisette2011automated,yano2011adaptation} applied to the nonlinear
system $\rbm^{(j)}(\ubm,\xbm;\Xi)$ to evaluate the function $\check\ubm$.
To safeguard this process, which is particularly important at early iterations
before the mesh has been substantially compressed, we take $n_\mathrm{ptc}$ PTC
steps (in this work, $n_\mathrm{ptc} = 5$) and only accept the modification if
(\ref{eqn:stpmod_cond}) is satisfied.

\subsubsection{Hessian regularization}
\label{subsubsec:hessreg}
Following the work in \cite{huang2022robust}, we take the regularization matrix to be
\begin{equation} \label{eqn:regmat}
	\Dbm_k^{(j)}(\Xi) =
	\partial_\ybm \phibold(\ybm_j^{(k)}(\Xi))^T \hat\Dbm_k^{(j)}(\Xi) \cdot \partial_\ybm \phibold(\ybm_j^{(k)}(\Xi)),
\end{equation}
where $\hat\Dbm_k^{(j)}(\Xi)\in\Rbb^{N_\xbm\times N_\xbm}$ is the stiffness matrix of a
continuous Galerkin finite element discretization of the linear elasticity equations. 
To avoid bulging of high-order elements, we set the Poisson's ratio to $\nu = 0$.
Unlike \cite{huang2022robust}, which was tailored to inviscid problems, selection of the
Young's modulus field is delicate for viscous problems. First, we observed that a
discontinuous Young's modulus field leads to mesh degradation, which was particularly
troublesome when compressing elements into thin features. Therefore, we choose the
Young's modulus to lie in the piecewise linear, continuous space
$E(\,\cdot\;\Xi) \in \Wcal_{h,1}$. Next, we aim to stiffen the
mesh in regions where the objective function is small and allow larger deformations
other regions. The enriched residual can vary by several orders of magnitude in
adjacent elements, which would be poorly represented by the $\Wcal_{h,1}$ space,
so we first introduce $\hat{C}(\,\cdot\,;\Xi) \in \Wcal_{h,1}$, where $\hat{C}$
is the $L^2$ projection of the logarithm of the normalized elemental enriched residual,
i.e., $C(\,\cdot\,;\Xi) : \Omega_0 \rightarrow \Rbb_{\geq 0}$ with
\begin{equation}
	C : (X, \Xi) \mapsto \log_{10}
	\left(
		\dfrac{\norm{(\hat\Pbm_{K(X)}^{(j)})^T\Rbm_\rho^{(j)}(\ubm_k^{(j)}(\Xi),\phibold(\ybm_k^{(j)}(\Xi)); \Xi)}}{\max_{K'\in\Ecal_{h,q}} \norm{(\hat\Pbm_{K'}^{(j)})^T\Rbm_\rho^{(j)}(\ubm_k^{(j)}(\Xi),\phibold(\ybm_k^{(j)}(\Xi)); \Xi)}}
	\right),
\end{equation}
where $K(X) \in \Ecal_{h,q}$ is the element in which $X\in\Omega_0$ lies.
The logarithm transforms a variable that varies over many orders of magnitude to one
that varies over a more limited range. Finally, we define the Young's modulus field
by converting back to the original scale, applying a range transformation, and
converting the compliance ($C$) into the stiffness, we have
\begin{equation} \label{eqn:youngmod}
	E : (X; \Xi) \mapsto g(10^{\hat{C}(X;\Xi)}; \Xi)^{-1},
\end{equation}
where $g(\,\cdot\,;\Xi) : [0,10^{\hat{C}_\mathrm{max}(\Xi)}] \rightarrow [C_\mathrm{l},C_\mathrm{u}]$ with
$g : (s; \Xi) \mapsto C_\mathrm{u} \theta(s;\Xi) + C_\mathrm{l} (1-\theta(s;\Xi))$
is the range transformation, $\theta(\,\cdot\,;\Xi) : [0,10^{\hat{C}_\mathrm{max}(\Xi)}]\rightarrow[0, 1]$
is a smooth cutoff function
\begin{equation}
	\theta: (s; \Xi) \mapsto \frac{1}{2} \left(\sin\left(\pi\frac{s}{10^{\hat{C}_\mathrm{max}(\Xi)}}-\frac{\pi}{2}\right) + 1\right),
\end{equation}
$(C_\mathrm{l}, C_\mathrm{u})$ is the range of $g$
($C_\mathrm{l} = 1$ and $C_\mathrm{u} = 10$ in this work),
and $\hat{C}_\mathrm{max}(\Xi) = \sup_{X\in\Omega_0} \hat{C}(X;\Xi)$.
A comparison between the Young's modulus distribution in \cite{zahr_implicit_2020}
and the proposed enriched residual-based distribution is illustrated
in Figure~\ref{fig:ns1_YoungE}.

\ifbool{fastcompile}{}{
	\begin{figure}[!htbp]
		\centering
		\begin{tikzpicture}[scale=0.8]
\begin{groupplot}[
  group style={
      group size=2 by 2,
      horizontal sep=1cm,
      vertical sep=0.5cm
  },
  width=0.7\textwidth,
  axis equal image,
  xlabel={$x_1$},
  ylabel={$x_2$},
  xtick = {-6.0, 0.0, 6.0},
  xticklabels={-6.0, 0.0, 6.0},
  ytick = {0.0, 3.0},
  xmin=-6.0, xmax=6.0,
  ymin=0, ymax=3.0
]

\nextgroupplot[xlabel={}, ylabel={}, xtick=\empty, ytick=\empty]
\addplot graphics [xmin=-6.0, xmax=6.0, ymin=0, ymax=3.0] {{_img/cyl_M5Re1e3_YoungE_invvol_phydom}.png};

\nextgroupplot[xlabel={}, ylabel={}, xtick=\empty, ytick=\empty]
\addplot graphics [xmin=-6.0, xmax=6.0, ymin=0, ymax=3.0] {{_img/cyl_M5Re1e3_YoungE_enres_phydom}.png};

\nextgroupplot[xlabel={}, ylabel={}, xtick=\empty, ytick=\empty]
\addplot graphics [xmin=-6.0, xmax=6.0, ymin=0, ymax=3.0] {{_img/cyl_M5Re1e3_YoungE_invvol_refdom}.png};

\nextgroupplot[xlabel={}, ylabel={}, xtick=\empty, ytick=\empty]
\addplot graphics [xmin=-6.0, xmax=6.0, ymin=0, ymax=3.0] {{_img/cyl_M5Re1e3_YoungE_enres_refdom}.png};

\end{groupplot}
\node[anchor=north] at ($(group c1r2.south)!0.5!(group c2r2.south)$) {\colorbarMatlabParulandscalse{0.03}{0.224}{0.418}{0.612}{0.806}{1.0}{10cm}};
\end{tikzpicture}
		\caption{
			Distribution of the Young's modulus field $E$ based on the
			inverse volume of reference mesh \cite{zahr_implicit_2020}
			(\textit{left}) and the enriched residual (\textit{right})
			for $M_\infty=5$, $\mathrm{Re}=10^3$ hypersonic flow over a
			cylinder (prior to $p$-adaptation). The Young's modulus are
			provided on the physical domain (\textit{top}) and reference
			domain (\textit{bottom}) for clarity.
		}
		\label{fig:ns1_YoungE}
	\end{figure}
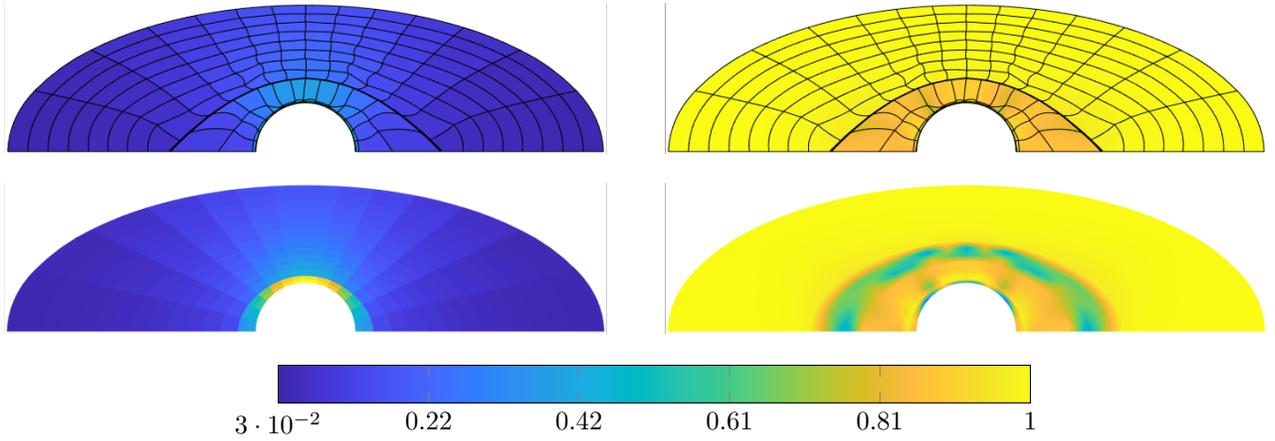
}

The regularization matrix $\Dbm_k^{(j)}(\Xi)$ controls the behavior of the regularization
across the domain $\Omega_0$. On the other hand, the regularization parameter
$\gamma_k^{(j)}(\Xi)$ controls the magnitude of the regularization, which must be
carefully balanced. Too little regularization leads to aggressive steps and
poor mesh quality, whereas too much regularization leads to minuscule step
sizes. Furthermore, it is important to have more regularization at early
iterations when the mesh and solution are far from convergence. To this end,
we define the regularization parameter as
\begin{equation} \label{eqn:regpar}
	\gamma_k^{(j)}(\Xi) = \frac{\hat\gamma}{\Xi^{\eta_1} k^{\eta_2}},
\end{equation}
where $\hat\gamma, \eta_1, \eta_2 \in\Rbb_{\ge 0}$ are user-defined parameters
($\eta_1 > 0$ if $\{\Xi_1,\dots,\Xi_c\}$ is an increasing sequence, e.g.,
$\Xi$ represents Reynolds number, and $\eta_2 < 0$ if $\{\Xi_1,\dots,\Xi_c\}$
is a decreasing sequence, e.g., $\Xi$ represents viscosity).
The explicit dependence on the optimization iteration $k$ and continutation parameter
$\Xi$ drives the regularization parameter toward zero as the SQP and continutation
iterations progress.

\begin{remark}
	For domains with planar boundaries, the parametrization $\phibold$ is affine
	in $\ybm$, which means its derivative is independent of $\ybm$. In these cases,
	$\Dbm_k^{(j)}(\Xi)$ does not depend on the state $\ybm_k^{(j)}(\Xi)$.
\end{remark}

\subsubsection{Viscosity continuation}
\label{subsubsec:viscont}
In the context of implicit shock tracking, continuation has been used to solve high-speed,
inviscid flows by using continuation with respect to Mach number \cite{huang2023high}.
For viscous flows, continuation with respect to viscosity is necessary to avoid carbuncles
and compress the mesh into thin shock and boundary layers. By starting with a large
viscosity, these features are thicker, which can easily be detected on
coarse, shock-agnostic meshes and require less extreme mesh compression
to adequately resolve. Once a reasonable mesh is obtained for the large,
initial viscosity ($\Xi_0$), the viscosity is decreased (to $\Xi_1$) and
the previous solution (mesh and flow variables) are used to initialize the SQP solver.
This process is repeated until the target viscosity is obtained ($\Xi_\star$).
Viscosity-based continuation was first introduced in the context of implicit shock
tracking in \cite{kercher2021moving}, and we demonstrate its importance in avoiding
carbuncles in Section~\ref{subsubsec:bow_recont}. Finally, we found that, provided
continuation is sufficiently gradual to avoid carbuncles, the final solution is
relatively insensitive to the continuation strategy.

\begin{remark} \label{rem:visccont}
	As is always the case with continuation, it is unnecessary to deeply converge
	the optimization problems for intermediate viscosities, i.e., $\Xi \neq \Xi_\star$.
	In practice, a predefined number of SQP iterations are performed for all intermediate
	viscosities.
\end{remark}

\subsection{Summary}
\label{subsec:summary}
To close this section, the complete HOIST method for viscous conservation laws is
summarized in Algorithms~\ref{alg:ist}-\ref{alg:ist-cont-padapt}. Algorithm~\ref{alg:ist}
summarizes the HOIST method for a fixed continuation parameter and $p$-adaptation level
to define the $\mathtt{HOIST}$ function ($r$-adaptation only). Algorithm~\ref{alg:ist-cont}
incorporates continuation into $\mathtt{HOIST}$ to define $\mathtt{HOISTC}$ for a fixed
$p$-adaptation level ($r$-adaptation only). Finally, Algorithm~\ref{alg:ist-cont-padapt} 
summarizes the complete $p$-adaptive HOIST method with continuation ($rp$-adaptation).

\begin{algorithm}
	\caption{HOIST method for viscous problems ($\mathtt{HOIST}$)}
	\label{alg:ist}
        \begin{algorithmic}[1]
		\REQUIRE Continuation parameter $\Xi$, $p$-adaptivity iteration $j$, initial state $\ubm_0^{(j)}(\Xi)$ and mesh $\ybm_0^{(j)}(\Xi)$, user-defined parameters ($\Delta$, $\kappa$, $\theta_\mathrm{l}$, $\theta_\mathrm{u}$, $\epsilon$, $\omega_\mathrm{mod}$, $C_\mathrm{l}$, $C_\mathrm{u}$, $\hat\gamma$, $\eta_1$, $\eta_2$)
		\ENSURE Converged state $\ubm_\star^{(j)}(\Xi)$ and mesh $\ybm_\star^{(j)}(\Xi)$
		\STATE \textbf{Residual scaling:} Compute $\rho$ from (\ref{eqn:weight-scaling}) to emphasize boundary layers
                \FOR{$k = 1, 2, \dots$}
			\STATE \textbf{Convergence check:} Terminate iterations if
			first-order optimality conditions of (\ref{eqn:pde-opt}) are
			satisfied \cite{zahr_implicit_2020,huang2022robust}
			\STATE \textbf{SQP search direction:} Compute $\Delta\zbm_k^{(j)}(\Xi)$ from (\ref{eqn:quadprog}) with regularization matrix from (\ref{eqn:regmat})--(\ref{eqn:regpar})
			\STATE \textbf{Step length constraint:} Compute $\hat\alpha_k^{(j)}(\Xi)$ from (\ref{eqn:stplencon})
			\STATE \textbf{Line search:} Compute step size $\alpha_k^{(j)}(\Xi) \in (0, \hat\alpha_k^{(j)}(\Xi)]$ such that $\ell_1$ merit function \cite{huang2022robust} satisfies sufficient decrease condition \cite{nocedal2006numerical}
			\STATE \textbf{SQP update:} Compute $\tilde\zbm_{k+1}^{(j)}(\Xi)$ from (\ref{eqn:stpupd})
			\STATE \textbf{State modification:} Compute $\zbm_{k+1}^{(j)}(\Xi)$ from (\ref{eqn:stpupd}) with $\chi_k^{(j)}$ defined in (\ref{eqn:stpmod})
                \ENDFOR
		\STATE \textbf{Finalize solution:} $\ubm_\star^{(j)}(\Xi) \leftarrow \ubm_k^{(j)}(\Xi)$ and $\ybm_\star^{(j)}(\Xi) \leftarrow \ybm_k^{(j)}(\Xi)$
        \end{algorithmic}
\end{algorithm}

\begin{algorithm}
	\caption{HOIST method for viscous problems with continuation ($\mathtt{HOISTC}$)}
	\label{alg:ist-cont}
	\begin{algorithmic}[1]
		\REQUIRE Continuation sequence $\{\Xi_1,\dots,\Xi_c=\Xi_\star\}$, $p$-adaptivity iteration $j$, initial state $\ubm_0^{(j)}(\Xi_1)$ and mesh $\ybm_0^{(j)}(\Xi_1)$, user-defined parameters $\Upsilon \coloneqq (\Delta, \kappa, \theta_\mathrm{l}, \theta_\mathrm{u}, \epsilon, \omega_\mathrm{mod}, C_\mathrm{l}, C_\mathrm{u}, \hat\gamma, \eta_1, \eta_2)$
		\ENSURE Converged state $\ubm_\star^{(j)}(\Xi_\star)$ and mesh $\ybm_\star^{(j)}(\Xi_\star)$
		\FOR{$i = 1, \dots, c$}
			\STATE \textbf{IST solve:} $(\ubm_\star^{(j)}(\Xi_i), \ybm_\star^{(j)}(\Xi_i)) = \mathtt{HOIST}(\Xi_i, j, \ubm_0^{(j)}(\Xi_i), \ybm_0^{(j)}(\Xi_i), \Upsilon)$
			\STATE \textbf{Continuation update:} $\ubm_0^{(j)}(\Xi_{i+1}) = \ubm_\star^{(j)}(\Xi_i)$ and $\ybm_0^{(j)}(\Xi_{i+1}) = \ybm_\star^{(j)}(\Xi_i)$
		\ENDFOR
	\end{algorithmic}
\end{algorithm}

\begin{algorithm}
	\caption{HOIST method for viscous problems with continuation and $p$-adaptivity ($\mathtt{pHOISTC}$)}
	\label{alg:ist-cont-padapt}
	\begin{algorithmic}[1]
		\REQUIRE Continuation sequence $\{\Xi_1,\dots,\Xi_c=\Xi_\star\}$, initial state $\ubm_0^{(0)}(\Xi_1)$ and mesh $\ybm_0^{(0)}(\Xi_1)$, user-defined parameters $\Upsilon \coloneqq (\Delta, \kappa, \theta_\mathrm{l}, \theta_\mathrm{u}, \epsilon, \omega_\mathrm{mod}, C_\mathrm{l}, C_\mathrm{u}, \hat\gamma, \eta_1, \eta_2)$, $p$-adaptivity convergence parameters $\hat\epsilon$ and $J_\mathrm{max}$
		\ENSURE Converged state $\ubm_\star^{(J)}(\Xi_\star)$ and mesh $\ybm_\star^{(J)}(\Xi_\star)$
		\STATE \textbf{IST continuation:} $(\ubm_\star^{(0)}(\Xi_\star), \ybm_\star^{(0)}(\Xi_\star)) = \mathtt{HOISTC}(\{\Xi_1,\dots,\Xi_c\}, 0, \ubm_0^{(0)}(\Xi_1), \ybm_0^{(0)}(\Xi_1), \Upsilon)$
		\FOR{$j = 1, \dots, J_\mathrm{max}$}
			\STATE \textbf{Convergence check:} Terminate iterations if conditions in (\ref{eqn:padapt_conv}) satisfied with $J = j$
			\STATE \textbf{IST solve:} $(\ubm_\star^{(j)}(\Xi_\star), \ybm_\star^{(j)}(\Xi_\star)) = \mathtt{HOIST}(\Xi_\star, j, \ubm_\star^{(j-1)}(\Xi_\star), \ybm_\star^{(j-1)}(\Xi_\star), \Upsilon)$
		\ENDFOR
	\end{algorithmic}
\end{algorithm}

\section{Numerical experiments}
\label{sec:numexp}
In this section, we apply the proposed $rp$-adaptive HOIST method to three
viscous Burgers' equation examples (Section~\ref{subsec:vburg}) and two
steady Navier-Stokes examples (Section~\ref{subsec:ns}) to demonstrate
its ability to accurately resolve viscous shocks (both straight-sided
and curved) and boundary layers, and accurately predict quantities of
interest (both integrated and profiles). We use the most interesting
example, hypersonic flow over a cylinder (Section~\ref{subsubsec:bow}),
to demonstrate the importance of the various solver enhancements proposed
in this work.

\subsection{Viscous Burgers' equation}
\label{subsec:vburg}
The time-dependent, viscous Burgers' equation governs nonlinear advection of a
scalar quantity $\func{\phi}{\Omega}{\Rbb}$ through a one-dimensional domain
$\Omega\subset\Rbb$
\begin{equation} \label{eqn:vburg}
  \pder{}{t}\phi(x,t) + \phi(x,t)\pder{}{x}\phi(x,t) = \nu \pderH{}{x}{2}\phi(x,t)
\end{equation}
for all $x\in\Omega$ and $t\in\Tcal$ with diffusion coefficient
$\nu\in\Rbb_{\ge 0}$. This unsteady problem is written as a conservation
law of the form (\ref{eqn:vclaw-phys0}) using a space-time formulation
\cite{kercher2021moving,zahr2020implicit}.

\subsubsection{Straight-sided space-time shock: pure $r$-adaptivity}
\label{subsubsec:vburg0}
First, we consider the example from \cite{kercher2021moving}, where
$\Omega = (0, 1)$ and $\Tcal = (0, 1)$ with initial condition
$\phi(x, 0) = \sin(2\pi x)/\pi + 0.2$ and constant boundary conditions
that are consistent with the initial condition. The sine wave steepens
into a viscous shock with thickness proportional to the viscosity, which
then propagates. We use this simple problem to demonstrate the capability
of the viscous HOIST method to resolve thin shocks with pure $r$-adaptation
(no $p$-adaptation).

We discretize the space-time domain using a structured mesh of $384$ straight-sided
($q=1$) triangles with quadratic solution approximation ($p : K \mapsto 2$). The
HOIST solver is initialized from a first-order finite volume solution of
(\ref{eqn:vburg}) with $\nu = 10^{-3}$ on the structured mesh
(Figure~\ref{fig:vburg0_sptm_ic}).
Viscosity continuation is performed in 11 stages with
$\Xi_i = (1 - 9(i-1)/100)\cdot 10^{-3}$ for $i=1,2,\dots,11$
with each continuation stage allotted $100$ SQP iterations.
For this problem, we do not use
residual weighting (Section~\ref{subsubsec:resscale}),
step length modifications (Section~\ref{subsubsec:stepcon})
or step modifications (Section~\ref{subsubsec:stepmod}) as
they proved unnecessary, and we used the Hessian regularization
matrix of \cite{zahr_implicit_2020}. The enrichment degree is $\Delta = 2$.
Due to the low polynomial degree used ($p : K \mapsto 2$), three
layers of elements are compressed into the shock thickness, which
decreases with viscosity as $\Ocal(\nu)$, to accurately resolve the steep gradient. 
At the final viscosity $\nu = 10^{-4}$, the shock is extremely thin.
The HOIST method also compresses the
mesh near the evolving peaks of the sine wave (gradient changes sign)
and near the shock formation singularity (Figure~\ref{fig:vburg0_sptm_physdom}).
In the reference domain $\Omega_0$, where the mesh is decompressed, the
internal shock structure can be clearly seen (Figure~\ref{fig:vburg0_sptm_refdom}).
Both of these figures demonstrate
that the space-time solution is well-resolved by the $r$-adaptation approach.
Temporal slices of the spatial solution field further show the solution
is well-resolved and oscillation-free, and the viscous shock is represented
as a continuous feature (Figure~\ref{fig:vburg0_slice}). Furthermore, temporal
slices of the HOIST solution match a reference solution computed with a
second-order finite volume method with superbee limiter on a grid with
2000 cells that is integrated with fourth-order Runge--Kutta (explicit)
with fixed time step $\Delta t = 0.0005$.

\ifbool{fastcompile}{}{
	\begin{figure}[!htbp]
		\centering
		\includegraphics[width=0.4\textwidth]{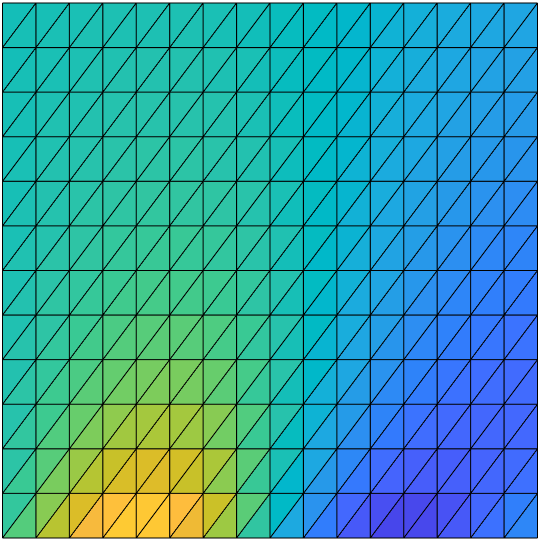}
		\colorbarMatlabParula{-0.118309886}{0.04}{0.20000000000000004}{0.35915494300000006}{0.518309886}
		\caption{
			Mesh and solution (first-order finite volume with $\nu=10^{-3}$)
			used to initialize HOIST method (pure $r$-adaptation) for viscous
			Burgers' problem with straight shock.
		}
		\label{fig:vburg0_sptm_ic}
	\end{figure}

	\begin{figure}[!htbp]
		\centering
		\begin{tikzpicture}
\begin{groupplot}[
  group style={
      group size=3 by 2,
      horizontal sep=0.3cm,
      vertical sep=0.3cm
  },
  width=0.45\textwidth,
  axis equal image,
  xlabel={$x$},
  ylabel={$t$},
  xtick = {0.0, 1.0},
  xticklabels={0.0, 1.0},
  ytick = {0.0, 0.5, 1.0},
  xmin=0.0, xmax=1.0,
  ymin=0, ymax=1.0
]

\nextgroupplot[xtick=\empty, xlabel={}, ytick=\empty, ylabel={}]
\addplot graphics [xmin=-0.0, xmax=1.0, ymin=0, ymax=1.0] {{_img/vburg_sptm_nu1e-3}.png};

\nextgroupplot[xtick=\empty, xlabel={}, ytick=\empty, ylabel={}]
\addplot graphics [xmin=-0.0, xmax=1.0, ymin=0, ymax=1.0] {{_img/vburg_sptm_nu5e-4}.png};

\nextgroupplot[xtick=\empty, xlabel={}, ytick=\empty, ylabel={}]
\addplot graphics [xmin=-0.0, xmax=1.0, ymin=0, ymax=1.0] {{_img/vburg_sptm_nu1e-4}.png};

\nextgroupplot[xtick=\empty, xlabel={}, ytick=\empty, ylabel={}]
\addplot graphics [xmin=-0.0, xmax=1.0, ymin=0, ymax=1.0] {{_img/vburg_sptm_nu1e-3_phydom_womsh}.png};

\nextgroupplot[xtick=\empty, xlabel={}, ytick=\empty, ylabel={}]
\addplot graphics [xmin=-0.0, xmax=1.0, ymin=0, ymax=1.0] {{_img/vburg_sptm_nu5e-4_phydom_womsh}.png};

\nextgroupplot[xtick=\empty, xlabel={}, ytick=\empty, ylabel={}]
\addplot graphics [xmin=-0.0, xmax=1.0, ymin=0, ymax=1.0] {{_img/vburg_sptm_nu1e-4_phydom_womsh}.png};

\end{groupplot}
\end{tikzpicture} \\
		\caption{
			HOIST solution (pure $r$-adaptation) \textit{in the physical
			domain $\Gcal(\Omega_0)$} to the viscous Burgers' problem with
			straight shock at $\nu = 10^{-3}, 5\times 10^{-4}, 10^{-4}$
			(\textit{left-to-right}) with and without mesh edges shown.
			Colorbar in Figure~\ref{fig:vburg0_sptm_ic}.
		}
		\label{fig:vburg0_sptm_physdom}
	\end{figure}

	\begin{figure}[!htbp]
                \centering
                \begin{tikzpicture}
\begin{groupplot}[
  group style={
      group size=3 by 2,
      horizontal sep=0.2cm,
      vertical sep=0.2cm
  },
  width=0.45\textwidth,
  axis equal image,
  xlabel={$x$},
  ylabel={$t$},
  xtick = {0.0, 1.0},
  xticklabels={0.0, 1.0},
  ytick = {0.0, 0.5, 1.0},
  xmin=0.0, xmax=1.0,
  ymin=0, ymax=1.0
]

\nextgroupplot[xtick=\empty, xlabel={}, ytick=\empty, ylabel={}]
\addplot graphics [xmin=0.0, xmax=1.0, ymin=0, ymax=1.0] {{_img/vburg_sptm_nu1e-3_refdom_wmsh}.png};

\nextgroupplot[xtick=\empty, xlabel={}, ytick=\empty, ylabel={}]
\addplot graphics [xmin=-0.0, xmax=1.0, ymin=0, ymax=1.0] {{_img/vburg_sptm_nu5e-4_refdom_wmsh}.png};

\nextgroupplot[xtick=\empty, xlabel={}, ytick=\empty, ylabel={}]
\addplot graphics [xmin=0.0, xmax=1.0, ymin=0, ymax=1.0] {{_img/vburg_sptm_nu1e-4_refdom_wmsh}.png};

\nextgroupplot[xtick=\empty, xlabel={}, ytick=\empty, ylabel={}]
\addplot graphics [xmin=-0.0, xmax=1.0, ymin=0, ymax=1.0] {{_img/vburg_sptm_nu1e-3_refdom_womsh}.png};

\nextgroupplot[xtick=\empty, xlabel={}, ytick=\empty, ylabel={}]
\addplot graphics [xmin=0.0, xmax=1.0, ymin=0, ymax=1.0] {{_img/vburg_sptm_nu5e-4_refdom_womsh}.png};

\nextgroupplot[xtick=\empty, xlabel={}, ytick=\empty, ylabel={}]
\addplot graphics [xmin=-0.0, xmax=1.0, ymin=0, ymax=1.0] {{_img/vburg_sptm_nu1e-4_refdom_womsh}.png};

\end{groupplot}
\end{tikzpicture} \\
                \caption{
                        HOIST solution (pure $r$-adaptation) \textit{in the reference
                        domain $\Omega_0$} to the viscous Burgers' problem with straight
                        shock at $\nu = 10^{-3}, 5\times 10^{-4}, 10^{-4}$
                        (\textit{left-to-right}) with and without mesh edges shown.
                        Colorbar in Figure~\ref{fig:vburg0_sptm_ic}.
                }
                \label{fig:vburg0_sptm_refdom}
        \end{figure}
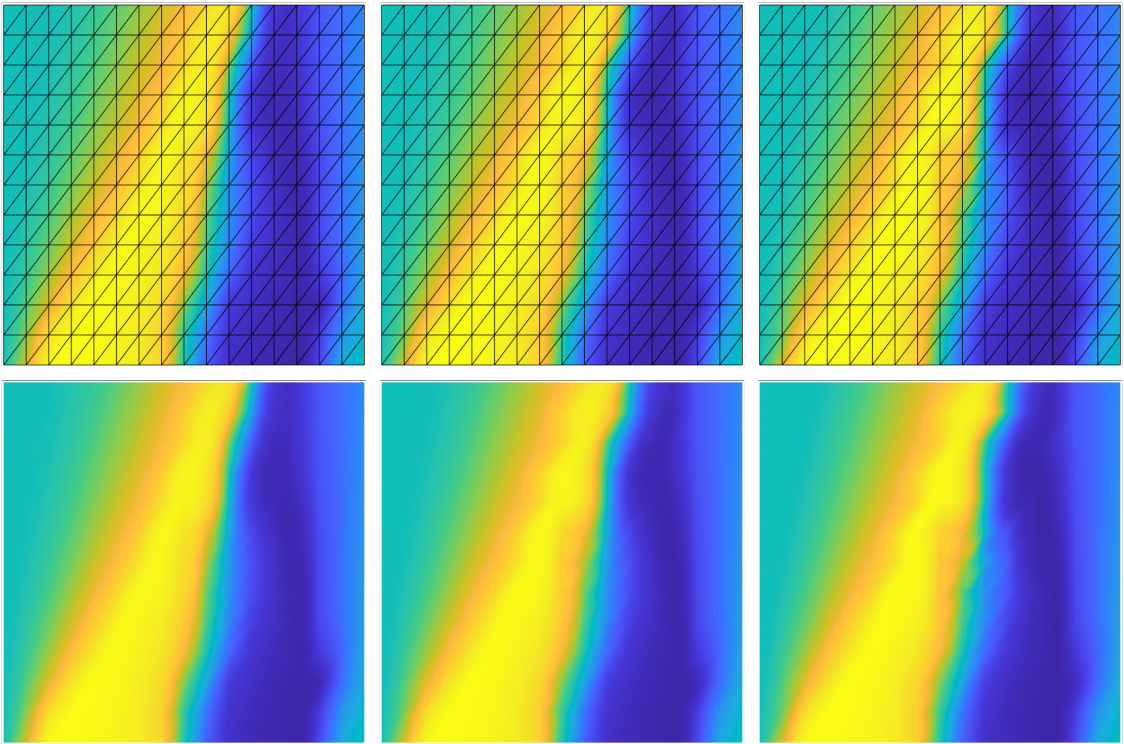

	\begin{figure}
		\centering
        \input{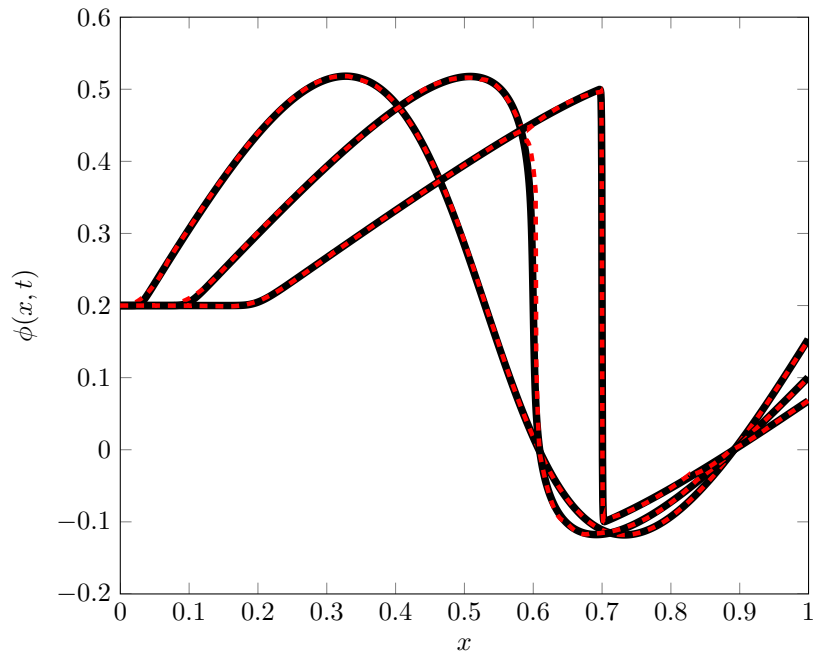}
		\caption{
			Temporal slices of the HOIST solution (pure $r$-adaptation)
			at times $t \in \{0.15, 0.5, 1.0\}$ (\ref{line:vburg_sptm_ist})
			and a reference solution computed with a highly refined
			second-order finite volume solution (\ref{line:vburg_sptm_fvmsoln})
			for the viscous Burgers' problem with straight shock.
		}
		\label{fig:vburg0_slice}
	\end{figure}
}

\subsubsection{One-dimensional steady shock: $r$- vs. $p$-adaptation study}
\label{subsubsec:vburg1}
Having demonstrated the ability of the proposed HOIST method ($r$-adaptation only)
to resolve thin viscous shocks in Section~\ref{subsubsec:vburg0}, we study how
to balance the two refinement mechanisms considered in this work: $r$- and
$p$-refinement. For this, we consider the steady viscous Burgers' equation in
(\ref{eqn:vburg}) with $\partial_t \phi = 0$ over the domain $\Omega = (-1, 1)$
with the following analytical solution
\begin{equation} \label{eqn:vburg1_soln}
	\phi(x) = \phi_R + \frac{\phi_L - \phi_R}{1 + \exp(x/\nu)},
\end{equation}
where $\phi_L = 1$ and $\phi_R = -1$ are the left and right states, respectively,
of the viscous shock.
To study the trade-off between $r$- and $p$-adaptation, we consider five families
of polynomial degree distributions parametrized by $\hat{p}\in\Nbb$ and apply the
HOIST method in pure $r$-adaptation mode. We repeat this study for
$\hat{p} = 1, 2, \dots, 9$ and $\nu = 0.1, 0.01$.
This study will provide insight into the optimal distribution of the elements
sizes for a given polynomial distribution and shock thickness, which help us
predict the behavior of the full $rp$-adaptive HOIST method. Because this
problem is one-dimensional and relatively simple, none of the robustness measures in
Section~\ref{subsec:slvrimprov} were used, except the enrichment degree $\Delta=2$.

The first family considers a mesh with a single element with polynomial degree
$\hat{p}_1 : K \mapsto \hat{p}$. This corresponds to a global polynomial approximation,
which is expected to poorly approximate the viscous shock in (\ref{eqn:vburg1_soln}),
especially because there are no free nodes to move for the $r$-adaptation procedure;
however, it is included as a control. The second family considers a mesh with three
elements with polynomial degree distribution
\begin{equation}
	\hat{p}_2 : K \mapsto
	\begin{cases}
		\hat{p}, & K = (-1/3, 1/3) \\
		1,       & \text{otherwise}.
	\end{cases}
\end{equation}
The remaining families consider a mesh with five elements with polynomial degree
distributions
\begin{equation}
	\begin{aligned}
		\hat{p}_3 &: K \mapsto
        	\begin{cases}
			\hat{p}, & K \in \{(-0.6, -0.2), (0.2, 0.6)\} \\
        	        1,       & \text{otherwise}
        	\end{cases}
		\\
		\hat{p}_4 &: K \mapsto
        	\begin{cases}
        	        \hat{p}, & K \in \{(-0.6, -0.2), (0.2, 0.6)\} \\
			3,       & K = (-0.2, 0.2) \\
        	        1,       & \text{otherwise}
        	\end{cases}
		\\
		\hat{p}_5 &: K \mapsto
        	\begin{cases}
        	        \hat{p}, & K \in \{(-0.6, -0.2), (0.2, 0.6)\} \\
			5,       & K = (-0.2, 0.2) \\
        	        1,       & \text{otherwise}.
        	\end{cases}
	\end{aligned}
\end{equation}

The different polynomial distributions lead to very different optimized meshes
(Figure~\ref{fig:vburg1_pdist}). For family $\hat{p}_2$, the central element is
contained within the shock thickness for lower degrees ($\hat{p} \leq 4$ for $\nu = 0.1$
and $\hat{p} \leq 6$ for $\nu = 0.01$), and expands somewhat beyond for higher
degrees. The extent to which the central element expands beyond the shock
thickness depends on both approximation degree and viscosity; for
$\nu=0.01$, expansion beyond the shock thickness is limited
even for $\hat{p}=9$. For the $\hat{p}_3$ family, which forces the central
element to use a linear approximation, the optimal element sizes
grow as $\hat{p}$ increases. However, because the central element
uses a linear approximation, it is forced to be very small, especially
for larger values of $\hat{p}$, because the linear approximation is only
valid in the neighborhood of $x=0$. The elements adjacent to the central
element can be much larger than the shock element in the $\hat{p}_2$
family because each one is only approximating half of the transition.
Owing to their similarity, $\hat{p}_4$ and $\hat{p}_5$ families show similar trends
as the $\hat{p}_3$ family. The main difference comes from the smallness requirement
on the central element. While the central element is forced to be minuscule
for the $\hat{p}_3$ family, especially for $\nu = 0.01$, it can be much
larger in the $\hat{p}_4$ and $\hat{p}_5$ families because of its cubic and quintic
approximation, respectively; however, it still remains within the shock.
For low values of $\hat{p}$, the central element is nearly the entire
shock width because the cubic/quintic approximation is needed to resolve
the shock transition and low-order approximations are forced to the regions
where the solution is nearly constant. Finally, we observe that the
low-order ($\hat{p} = 1,2$) members of the $\hat{p}_2$, $\hat{p}_3$, and $\hat{p}_4$
families  have three elements compressed in the shock thickness for $\nu = 0.01$,
which is consistent with the optimized mesh computed in Section~\ref{subsubsec:vburg0}.
However, a quintic approximation in the central element ($\hat{p}_5$ family)
relaxes this requirement by allowing the other elements to extend beyond
the shock thickness. This behavior is also produced by enriching the
approximation degree of the elements adjacent to the central element
($\hat{p}_3$ and $\hat{p}_4$ families).

Despite the substantially different meshes produced for the different
polynomial degree distributions, the $L^2(\Omega_0)$ for all families
decays at roughly the same rate with respect to the number of flow variable
degrees of freedom ($N_\ubm$) (Figure~\ref{fig:vburg1_err}), and have similar
error constants. The one exception is the global polynomial degree $\hat{p}_1$,
which, as expected, only decays for the larger viscosity where the shock
is a relatively slow transition. From this study, we conclude that there
is not a preferred approach to balance $r$- and $p$-adaptation, at least
in terms of error vs the number of degrees of freedom. However, the optimal mesh strongly depends on the
chosen polynomial distribution. This means that the balance between $r$-
and $p$-adaptation can be chosen based on other criteria, e.g., realizability
of the optimal mesh, computational cost, or overall convenience. In the
$rp$-adaptation setting, we expect two extremes. If $r$-adaptivity is
prioritized at lower polynomial degrees, more element compression is
expected near the shocks and boundary layers and, if the elements are
sufficiently small, $p$-adaptation will not necessarily be focused
within the shock or boundary layer. This situation would produce
distributions similar to the $\hat{p}_3$ family. On the other hand, if
$p$-adaptivity is prioritized, the polynomial degree in the shock
and boundary layer would be rapidly increased due to underresolution,
and the element compression would be more modest when $r$-adaptation
begins, especially for larger viscosities. This situation would produce
distributions similar to the $\hat{p}_2$ family. As we will see, the proposed
$rp$-adaptive method tends to produce $\hat{p}_3$ family solutions because
(1) it allows a majority of the work to be performed on the lower $p$ discretizations
with fewer degrees of freedom and (2) the element relaxation afforded
by larger polynomial degrees is less beneficial for very thin features.

\ifbool{fastcompile}{}{
	\begin{figure}[!htbp]
		\centering
		\begin{tikzpicture}
\begin{groupplot}[
  group style={
      group size=2 by 4,
      horizontal sep=0.65cm,
      vertical sep=0.125cm
  },
  width=0.5\textwidth,
  axis equal image,
  xlabel={$x$},
  ylabel={},
  xtick = {-7.0, -3.5, 0.0, 3.5, 7.0},
  xticklabels={-1.0, -0.5, 0.0, 0.5, 1.0},
  ytick = {1, 3, 5, 7, 9},
  xmin=-7, xmax=7,
  ymin=0, ymax=10
]

\nextgroupplot[xlabel={}, xtick=\empty, ylabel={$\hat{p}$}]
\addplot graphics [xmin=-7, xmax=7, ymin=0, ymax=10] {{_img/1DBurgers_1p1}.pdf};

\nextgroupplot[xlabel={}, xtick=\empty, ytick=\empty]
\addplot graphics [xmin=-7, xmax=7, ymin=0, ymax=10] {{_img/1DBurgers_1p1_nu0.01}.pdf};

\nextgroupplot[xlabel={}, xtick=\empty, ylabel={$\hat{p}$}]]
\addplot graphics [xmin=-7, xmax=7, ymin=0, ymax=10] {{_img/1DBurgers_1p1p1}.pdf};

\nextgroupplot[xlabel={}, xtick=\empty, ytick=\empty]
\addplot graphics [xmin=-7, xmax=7, ymin=0, ymax=10] {{_img/1DBurgers_1p1p1_nu0.01}.pdf};

\nextgroupplot[xlabel={},xtick=\empty, ylabel={$\hat{p}$}]]
\addplot graphics [xmin=-7, xmax=7, ymin=0, ymax=10] {{_img/1DBurgers_1p3p1}.pdf};

\nextgroupplot[xlabel={}, xtick=\empty, ytick=\empty]
\addplot graphics [xmin=-7, xmax=7, ymin=0, ymax=10] {{_img/1DBurgers_1p3p1_nu0.01}.pdf};

\nextgroupplot[xticklabels={-1.0, -0.5, 0.0, 0.5, 1.0}, ylabel={$\hat{p}$}]]
\addplot graphics [xmin=-7, xmax=7, ymin=0, ymax=10] {{_img/1DBurgers_1p5p1}.pdf};

\nextgroupplot[xticklabels={-1.0, -0.5, 0.0, 0.5, 1.0}, ylabel={}, ytick=\empty]
\addplot graphics [xmin=-7, xmax=7, ymin=0, ymax=10] {{_img/1DBurgers_1p5p1_nu0.01}.pdf};

\end{groupplot}
\end{tikzpicture}\\\vspace{1mm}
		\begin{tikzpicture}
    \begin{axis}[
        colorbar,
        colormap/viridis,
        colorbar horizontal,
        point meta min=0,
        point meta max=1,
        colorbar style={
            width=10cm,
            xtick={0, 0.125, 0.25, 0.375, 0.5, 0.625, 0.75, 0.875, 1},
            xticklabels={1, 2, 3, 4, 5, 6, 7, 8, 9},
        },
        hide axis,
    ]
        \addplot [draw=none] coordinates {(0,0)};
    \end{axis}
\end{tikzpicture}
		\caption{
			HOIST solution (pure $r$-adaptation) for the steady Burgers'
			problem for four families of polynomial distributions
			$\hat{p}_2$, $\hat{p}_3$, $\hat{p}_4$, $\hat{p}_5$ (\textit{top-to-bottom})
			for parameter $\hat{p} = 1, 2, \dots, 9$ and viscosities
			$\nu = 0.1$ (\textit{left}) and $\nu = 0.01$ (\textit{right}).
			The vertical lines delineate the shock thickness, defined
			as $x\in[-\delta, \delta]$, where
			$|\phi(x)| \leq 0.95\max\{|\phi_L|, |\phi_R|\}$.
			The color corresponds to the polynomial degree used for
			each element.
		}
		\label{fig:vburg1_pdist}
	\end{figure}
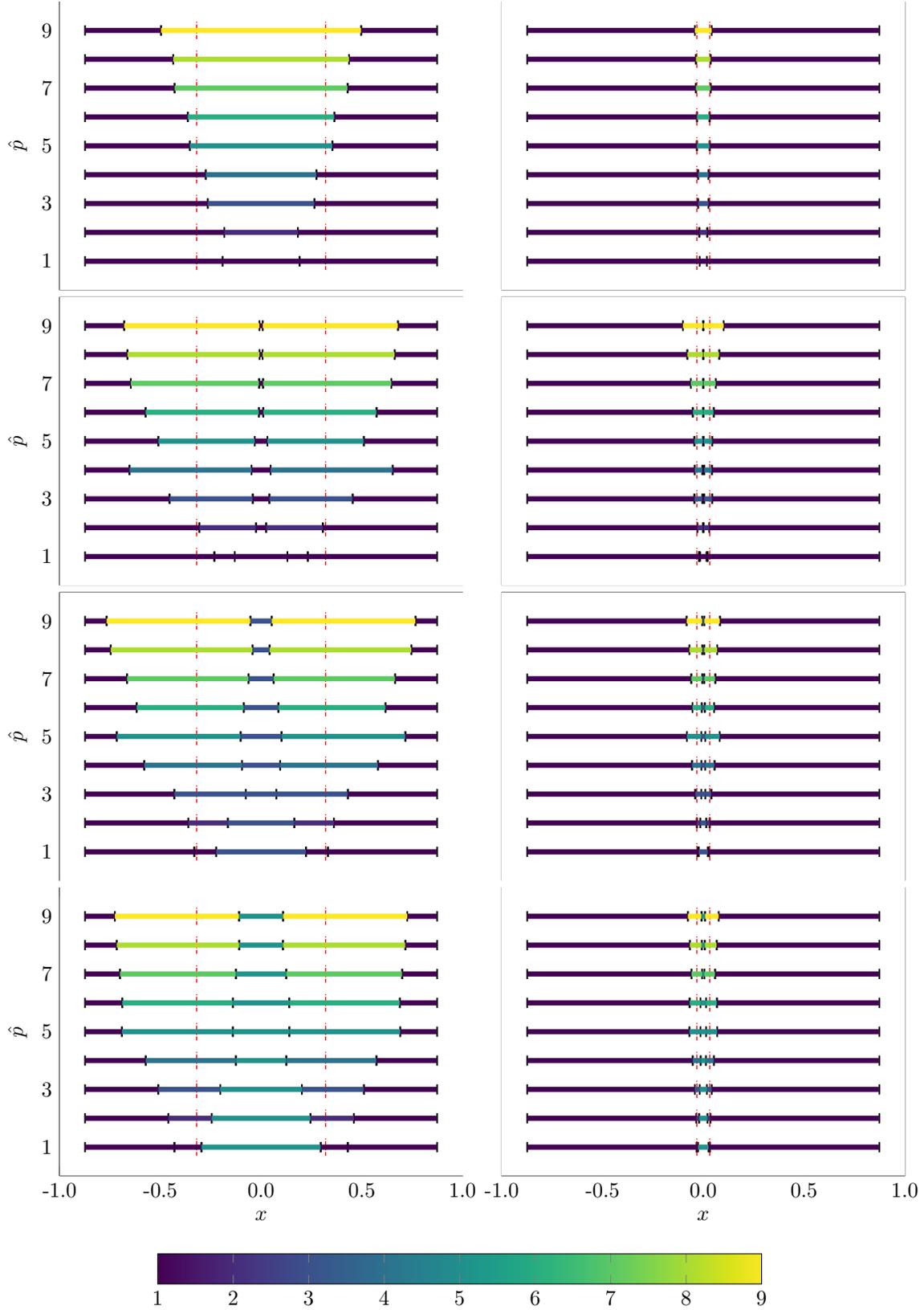
	
	\begin{figure}[!htbp]
		\centering
		\begin{tikzpicture}
\begin{groupplot}[
    group style={
    group size=2 by 1,
    horizontal sep=0.5cm, 
    vertical sep=0.5cm, 
    },
    width=0.5\textwidth,
    xlabel={$N_\ubm$},
    ylabel={$L_2(\Omega_0)$ error},
    legend style={at={(1.0,-0.2)},anchor=north, legend columns=-1},
    xtick = {0, 5, 10, 15, 20, 25, 30},
    xticklabels={0, 5, 10, 15, 20, 25, 30},
    ymin=1e-5, ymax=1e0,
]

\nextgroupplot[ymode=log]
\addplot [mark=none, solid, teal, line width=1pt] coordinates{
    (2, 0.24994696053560944)
    (3, 0.6227961825463061)
    (4, 0.47436375326325697)
    (5, 0.5937847811405941)
    (6, 0.16597912843142207)
    (7, 0.15187591669108585)
    (8, 0.04780711640434703)
    (9, 0.058919468671209174)
    (10, 0.042027055628262935)
   }; \label{line:p1}

\addplot  [mark=*, solid, red, line width=1pt] coordinates{
    (6, 0.14457349320455)
    (7, 0.12479080545655334)
    (8, 0.05407595845904571)
    (9, 0.05198724805475345)
    (10, 0.020770994629963523)
    (11, 0.018580143819406064)
    (12, 0.00825281948116569)
    (13, 0.00811352932190453)
    (14, 0.004213302411163727)
}; \label{line:p2}
\addplot [mark=diamond*, densely dashed, green!60!black
,line width=1pt]  coordinates{
     (10, 0.08101139967898886)
     (12, 0.030375813113598805)
     (14, 0.0052104824297294315)
     (16, 0.0020525089382057334)
     (18, 0.00306546151101574)
     (20, 0.0011666411757563025)
     (22, 0.00027124767988688715)
     (24, 0.0001283419748289417)
     (26, 0.00012130329519827549)
}; \label{line:p3}
\addplot [mark=square*, densely dashdotted, blue
        ,line width=1pt] coordinates{
    (12, 0.025314547073052127)
    (14, 0.015326383522546625)
    (16, 0.0062710128834093356)
    (18, 0.0009423291282421506)
    (20, 0.0005157138718284281)
    (22, 0.0005708337935673047)
    (24, 0.0002637348060543162)
    (26, 6.9083579712165e-5)
    (28, 3.0674372128932024e-5) 
}; \label{line:p4}
\addplot [mark=10-pointed star, dotted, purple
        ,line width=1pt] coordinates{
    (14, 0.008784408203865906)
    (16, 0.00521570941867164)
    (18, 0.0025298362056599983)
    (20, 0.0010101958684784241)
    (22, 0.00016009787943985698)
    (24, 0.0001642656641053826)
    (26, 0.00010925450239372656)
    (28, 5.4843164057404295e-5)
    (30, 1.855957382181568e-5)
}; \label{line:p5}

\nextgroupplot[ymode=log,  ylabel={}, ytick={}, yticklabels={}]]
\addplot [mark=none, solid, teal
,line width=1pt] coordinates{
    (3, 0.615883557759597)
    (4, 0.3733951404514678)
    (5, 0.4514817654032604)
    (6, 0.6955164568351679)
    (7, 0.7051496190618032)
    (8, 0.5400616155388636)
    (9, 0.5661140719158829)
    (10, 0.7136589509066146)
};
\addplot  [mark=*, solid, red
  ,line width=1pt] coordinates{
    (6, 0.15512889928557838)
    (7, 0.1295192051821777)
    (8, 0.04712310807570468)
    (9, 0.05628899422621676)
    (10, 0.03414952477481253)
    (11, 0.03556906751933864)
    (12, 0.01759775892137677)
    (13, 0.01885734807598146)
    (14, 0.011018448723992552)
};
\addplot [mark=diamond*, densely dashed, green!60!black
 ,line width=1pt]  coordinates{
    (10, 0.07285169889815842)
    (12, 0.009964438395546983)
    (14, 0.007394802989222622)
    (16, 0.002300050250119141)
    (18, 0.0013545246393748278)
    (20, 0.001596042849205843)
    (22, 0.00045751975556640115)
    (24, 0.00019180028409173663)
    (26, 0.0001240039426276915)
};
\addplot [mark=square*, densely dashdotted, blue
,line width=1pt] coordinates{
    (12, 0.03051945382754236)
    (14, 0.012925695205070323)
    (16, 0.008091505293871832)
    (18, 0.001355353625312277)
    (20, 0.0005172936054330627)
    (22, 0.0005637224919034511)
    (24, 0.0005729056014479059)
    (26, 0.00029635718620656907)
    (28, 5.067708711091364e-5)
};
\addplot [mark=10-pointed star, 
dotted, purple, line width=1pt] coordinates{
    (14, 0.026068034907123054)
    (16, 0.00815587776832746)
    (18, 0.007780198690730926)
    (20, 0.001716007137980186)
    (22, 0.0005615420772032045)
    (24, 0.0004328487902248574)
    (26, 9.910992771732198e-5)
    (28, 9.729824553572489e-5)
    (30, 3.451493898415422e-5)
};

\end{groupplot}
\end{tikzpicture}
		\caption{
			HOIST error (pure $r$-adaptation) for the steady Burgers'
			problem for the five families of polynomial distributions
			$\hat{p}_1$ (\ref{line:p1}),
			$\hat{p}_2$ (\ref{line:p2}),
			$\hat{p}_3$ (\ref{line:p3}),
			$\hat{p}_4$ (\ref{line:p4}),
			$\hat{p}_5$ (\ref{line:p5}) and viscosities
			$\nu = 0.1$ (\textit{left}) and $\nu = 0.01$ (\textit{right})
			as a function of the number of solution degrees of freedom.
		}
		\label{fig:vburg1_err}
	\end{figure}
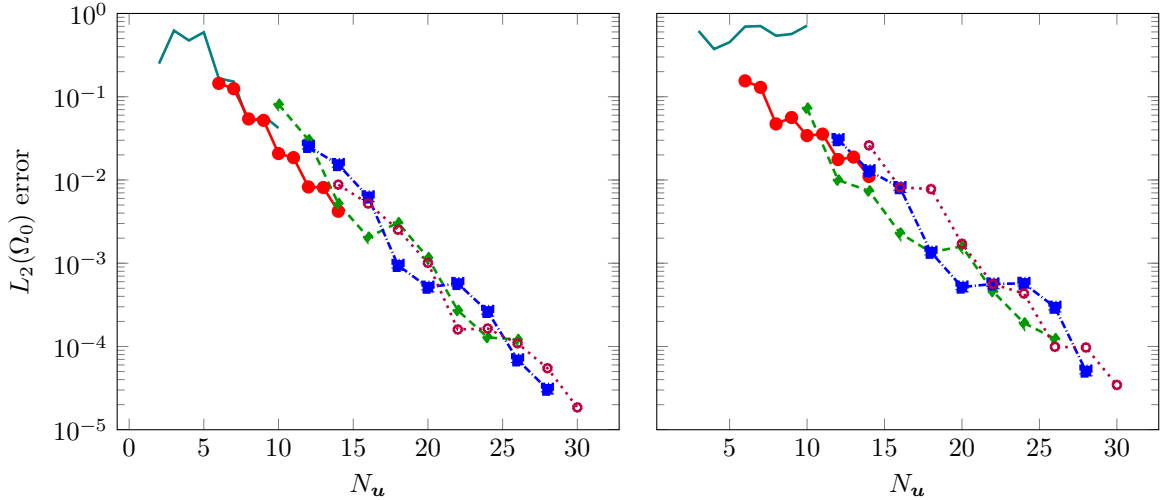
}

\subsubsection{Curved space-time shock: $rp$-adaptivity}
Next, we consider a modification of a problem found in \cite{masatsuka2009like},
defined by the spatial domain $\Omega = (-0.4, 1)$, temporal domain
$\Tcal = (0, 0.8)$, and smooth initial condition
\begin{equation}
\phi(x, 0) = \left \{ \begin{array}{rcl}
\phi_L \left(  -1 + \frac{2}{1+\exp(5x)} \right ) & \mbox{if} & x \leq 0 \\
&&\\
\phi_R (1-x) \left(  -1 + \frac{2}{1+\exp(5x)} \right ) & \mbox{if} & x > 0
\end{array}\right.
\end{equation}
where $\phi_L = 4$ and $\phi_R = 3$. The boundary conditions are constant
in time and consistent with the initial condition. This problem features
a smooth initial condition that steepens into a shock at a finite time,
which then accelerates until it leaves the domain. Despite the overall
simplicity of this problem, the solution possess two challenging
features in the space-time domain: (1) a singularity from the shock formation
and (2) a curved shock from the acceleration. We use this problem to
demonstrate the full $rp$-adaptive HOIST method.

\begin{figure}[!htbp]
	\centering
		\includegraphics[width=0.45\textwidth]{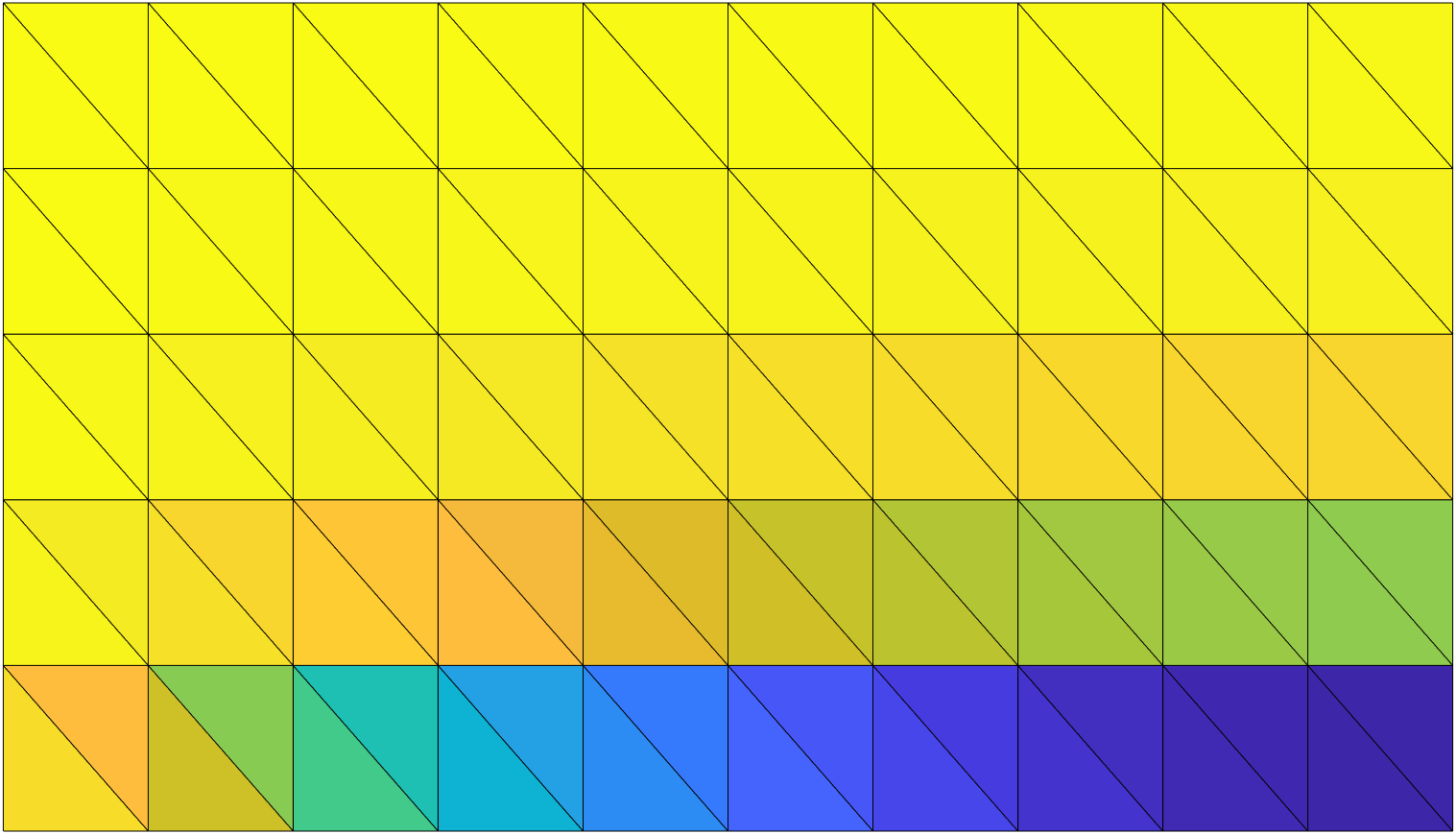}
		\colorbarMatlabParula{-1.42}{-0.20}{0.9}{2.0}{3.06}
		\caption{
			Mesh and solution (first-order finite volume with $\nu = 10^{-1}$)
			used to initialize the $rp$-adaptive HOIST method for the viscous
			Burgers' problem with curved shock.
		}
	\label{fig:vburg2_sptm_ic}
\end{figure}

We initially discretize the space-time domain using a structured mesh of $100$
straight-sided ($q=1$) elements with linear solution approximation ($p_0 : K \mapsto 1$).
The HOIST solver is initialized from a first-order finite volume solution of
(\ref{eqn:vburg}) with $\nu = 10^{-1}$ (Figure~\ref{fig:vburg2_sptm_ic}).
Viscosity continuation is performed in $20$ stages of equal increments beginning
with $\Xi_1 = 10^{-1}$ and ending with $\Xi_{20} = 10^{-3}$ with each continuation
stage allotted $10$ SQP iterations.
For this problem, we do not use
residual weighting (Section~\ref{subsubsec:resscale}),
step length modifications (Section~\ref{subsubsec:stepcon})
or step modifications (Section~\ref{subsubsec:stepmod}) as
they proved unnecessary, and we used the Hessian regularization
matrix of \cite{zahr_implicit_2020}. The enrichment degree is $\Delta = 2$.
After viscosity continuation is complete, the
mesh approximation is increased to $q = 2$ and additional SQP iterations are
performed for improved representation of the curved shock. A quadratic
mesh could have been used from the outset; however, we found that the
robustness of the SQP solver is improved and the mesh compression is
more severe by initializing with $q=1$. Finally, we perform three
rounds of $p$-refinement using the enriched residual error indicator
$\hat{s}_j^\mathrm{uwr}$.

\begin{figure}[!htbp]
	\centering
	\input{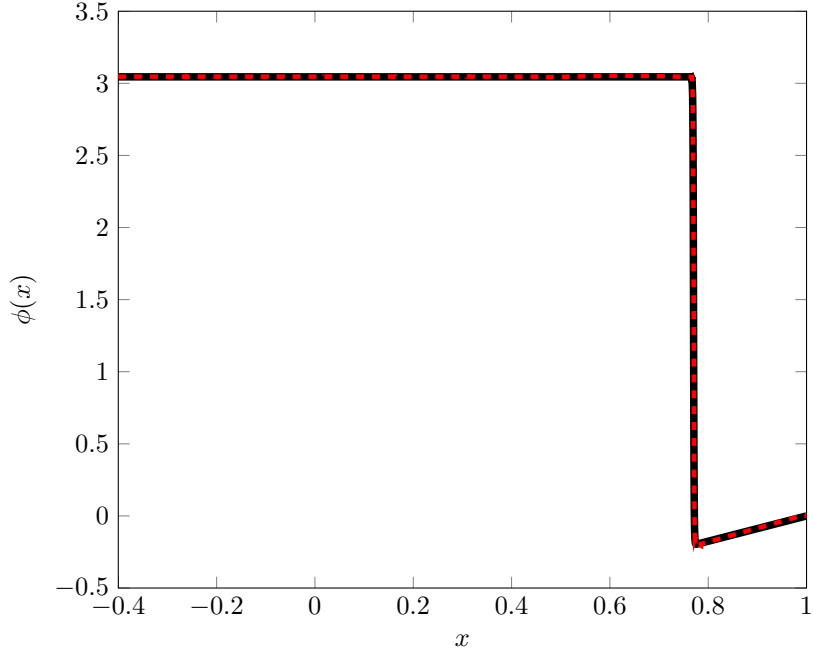}
	\caption{
		Temporal slice of the HOIST solution ($rp$-adaptation) at $t=0.8$
		(\ref{line:vburg_curshk_ist}) and a reference solution computed with a
		highly refined second-order finite volume method
		(\ref{line:vburg_curshk_fvm}) for the viscous Burgers' problem
		with curved shock.
	}
	\label{fig:vburg2_slice}
\end{figure}

Prior to $p$-adaptation, two elements are
compressed into the shock thickness and the enriched residual effectively
identifies the elements inside and near the shock as the main contributors
to the error. The $p$-refinement approach locally increases the polynomial
degree in these elements and the HOIST method adjusts the mesh appropriately.
In the end, the mesh has been effectively compressed near the singularity
and shock, and the polynomial degrees appropriately enhanced to reduce
the magnitude of the error indicator (Figure~\ref{fig:vburg2_sptm_physdom}).
The space-time solutions in the reference domain $\Omega_0$
clearly show the internal structure of the shock and the distribution
of the polynomial degree (Figure~\ref{fig:vburg2_sptm_refdom}). Both of
these figures demonstrate that the space-time solution is well-resolved by
the $rp$-adaptation approach. Temporal slices of the spatial solution field
further show the solution is well-resolved and oscillation-free, and the viscous
shock is represented as a continuous feature (Figure~\ref{fig:vburg2_slice}).
Furthermore, temporal slices of the HOIST solution match a reference solution
computed with a second-order finite volume method with superbee limiter on a grid
with 2000 cells that is integrated with fourth-order Runge--Kutta (explicit) with
fixed time step $\Delta t = 0.001$.

\begin{figure}[!htbp]
	\centering
	\input{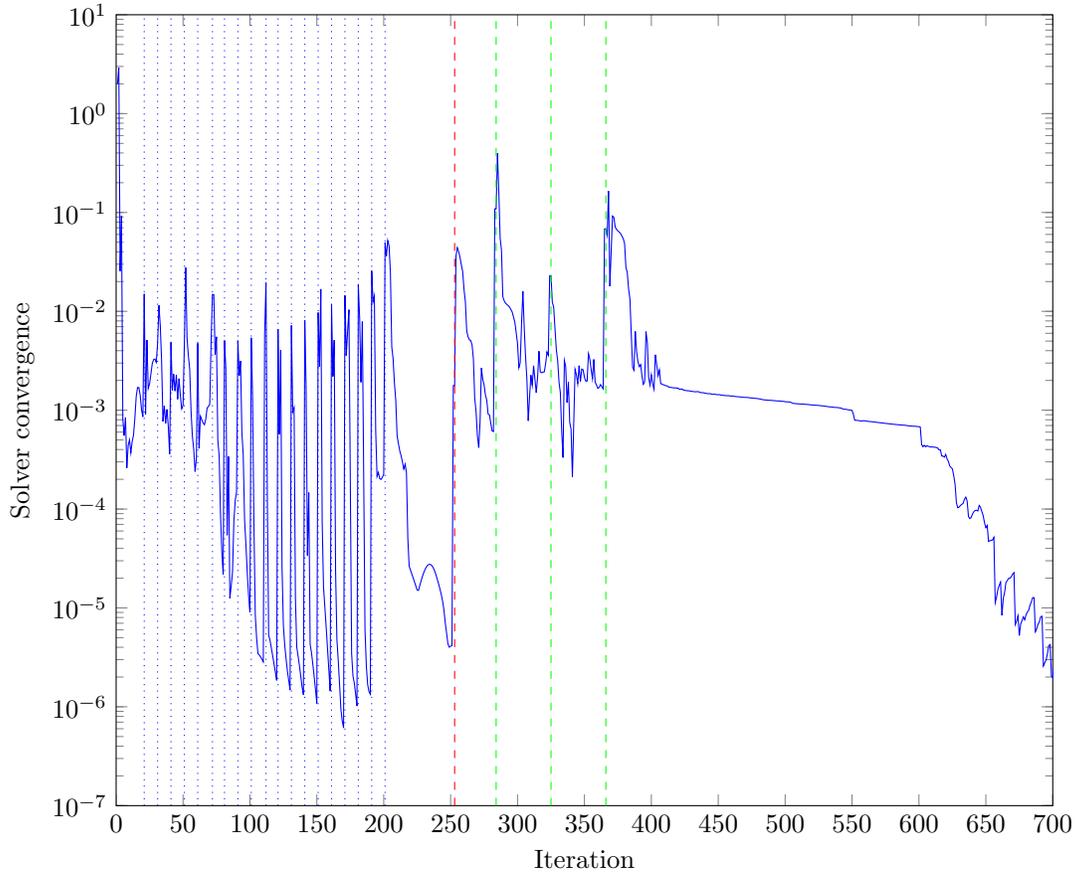}
	\caption{
		Convergence of the DG residual
		$\norm{\rbm^{(j)}(\ubm_k^{(j)}, \phibold(\ybm_k^{(j)}); \Xi_i)}$
		(\ref{line:vburg_curshk_c_nrm}) across all viscosity continuation stages
		$\{\Xi_i\}_{i=1}^{20}$ and $p$-adaptivity iterations $\{p_j\}_{j=0}^3$.
		Vertical lines indicate the end of a continuation state
		(\ref{line:padpat_c_nrm2}) (i.e., transition from $\Xi_i$ to $\Xi_{i+1}$)
		or the end of a $p$-adaptation iteration (\ref{line:padpat_c_nrm3})
		(i.e., transition from $p_j$ to $p_{j+1}$) for the viscous Burgers'
		problem with curved shock.
	}
	\label{fig:vburg2_solver}
\end{figure}

The behavior of the SQP solver across all continuation and $p$-adaptation stages
is shown in Figure~\ref{fig:vburg2_solver}. For a given value of the viscosity
$\Xi_i$ and polynomial distribution $p_j$, the HOIST method drives the residual
toward zero. However, the residual jumps up after the continuation parameter or
polynomial distribution is updated because the residual function changes. Notice
that 250 iterations are performed on the coarsest (i.e., least expensive)
discretization that uses a linear approximations in all elements. The
$p$-adaptation process is complete in another 150 iterations, at which
point the mesh is effectively frozen. We allow the SQP solver to run
an additional 300 iterations to obtain deep convergence, although a
nonlinear solve on a fixed mesh could have been used to obtain the
solution in only a few iterations.

\begin{figure}[!htbp]
	\centering
	\begin{tikzpicture}[scale=0.8]
\begin{groupplot}[
  group style={
      group size=3 by 4,
      horizontal sep=0.25cm,
      vertical sep=0.25cm
  },
  width=0.5\textwidth,
  axis equal image,
  xlabel={$x$},
  ylabel={$t$},
  xtick = {-0.4, 1.0},
  xticklabels={-0.4, 1.0},
  ytick = {0.0, 0.4, 0.8},
  xmin=-0.4, xmax=1.0,
  ymin=0, ymax=0.8
]

\nextgroupplot[xtick=\empty, xlabel={}, ytick=\empty, ylabel={}]
\addplot graphics [xmin=-0.4, xmax=1.0, ymin=0, ymax=0.8] {{_img/visBurgs_u_0th_phydom_mesh}.png};

\nextgroupplot[xtick=\empty, xlabel={}, ytick=\empty, ylabel={}]
\addplot graphics [xmin=-0.4, xmax=1.0, ymin=0, ymax=0.8] {{_img/visBurgs_error_0th_phydom}.png};

\nextgroupplot[xtick=\empty, xlabel={}, ytick=\empty, ylabel={}]
\addplot graphics [xmin=-0.4, xmax=1.0, ymin=0, ymax=0.8] {{_img/visBurgs_pdeg_0th_phydom}.png};

\nextgroupplot[xtick=\empty, xlabel={}, ytick=\empty, ylabel={}]
\addplot graphics [xmin=-0.4, xmax=1.0, ymin=0, ymax=0.8] {{_img/visBurgs_u_1st_phydom_mesh}.png};

\nextgroupplot[xtick=\empty, xlabel={}, ytick=\empty, ylabel={}]
\addplot graphics [xmin=-0.4, xmax=1.0, ymin=0, ymax=0.8] {{_img/visBurgs_error_1st_phydom}.png};

\nextgroupplot[xtick=\empty, xlabel={}, ytick=\empty, ylabel={}]
\addplot graphics [xmin=-0.4, xmax=1.0, ymin=0, ymax=0.8] {{_img/visBurgs_pdeg_1st_phydom}.png};

\nextgroupplot[xtick=\empty, xlabel={}, ytick=\empty, ylabel={}]
\addplot graphics [xmin=-0.4, xmax=1.0, ymin=0, ymax=0.8] {{_img/visBurgs_u_2nd_phydom_mesh}.png};

\nextgroupplot[xtick=\empty, xlabel={}, ytick=\empty, ylabel={}]
\addplot graphics [xmin=-0.4, xmax=1.0, ymin=0, ymax=0.8] {{_img/visBurgs_error_2nd_phydom}.png};

\nextgroupplot[xtick=\empty, xlabel={}, ytick=\empty, ylabel={}]
\addplot graphics [xmin=-0.4, xmax=1.0, ymin=0, ymax=0.8] {{_img/visBurgs_pdeg_2nd_phydom}.png};

\nextgroupplot[xtick=\empty, xlabel={}, ytick=\empty, ylabel={}]
\addplot graphics [xmin=-0.4, xmax=1.0, ymin=0, ymax=0.8] {{_img/visBurgs_u_3rd_phydom_mesh}.png};

\nextgroupplot[xtick=\empty, xlabel={}, ytick=\empty, ylabel={}]
\addplot graphics [xmin=-0.4, xmax=1.0, ymin=0, ymax=0.8] {{_img/visBurgs_error_3rd_phydom}.png};

\nextgroupplot[xtick=\empty, xlabel={}, ytick=\empty, ylabel={}]
\addplot graphics [xmin=-0.4, xmax=1.0, ymin=0, ymax=0.8] {{_img/visBurgs_pdeg_3rd_phydom}.png};

\end{groupplot}
\node[anchor=west] at ($(group c1r4.west) + (-0.85cm, -3)$) {\colorbarMatlabParulandscalse{-1.42}{-0.1}{1.0}{2.0}{3.06}{3.06}{5cm}};

\node[anchor=west] at ($(group c2r4.west) + (-0.28cm, -3.3)$) {\colorbarMatlabParulandscalse{2e-6}{2e-2}{4.1e-2}{4.1e-2}{4.1e-2}{4.1e-2}{4.8cm}};

\node[anchor=west] at ($(group c3r4.west) + (+0.1cm, -3)$) {\colorbarMatlabParulandscalse{1}{2}{3}{4}{4}{4}{5cm}};
\end{tikzpicture}
	\caption{
		HOIST solution ($rp$-adaptation) (\textit{left}), the enriched residual
		error indicator (\textit{middle}), and polynomial degree distribution
		(\textit{right}) \textit{in the physical domain} $\Gcal(\Omega_0)$
		to the viscous Burgers' problem with curved shock. These quantities
		are provided prior to $p$-adaptation and after each $p$-adaptation
		iteration (\textit{top-to-bottom}).
	}
	\label{fig:vburg2_sptm_physdom}
\end{figure}

\begin{figure}[!htbp]
	\centering
	\begin{tikzpicture}[scale=0.8]
\begin{groupplot}[
  group style={
      group size=3 by 4,
      horizontal sep=0.25cm,
      vertical sep=0.25cm
  },
  width=0.5\textwidth,
  axis equal image,
  xlabel={$x$},
  ylabel={$t$},
  xtick = {-0.4, 1.0},
  xticklabels={-0.4, 1.0},
  ytick = {0.0, 0.4, 0.8},
  xmin=-0.4, xmax=1.0,
  ymin=0, ymax=0.8
]

\nextgroupplot[xtick=\empty, xlabel={}, ytick=\empty, ylabel={}]
\addplot graphics [xmin=-0.4, xmax=1.0, ymin=0, ymax=0.8] {{_img/visBurgs_u_0th_refdom}.png};

\nextgroupplot[xtick=\empty, xlabel={}, ytick=\empty, ylabel={}]
\addplot graphics [xmin=-0.4, xmax=1.0, ymin=0, ymax=0.8] {{_img/visBurgs_error_0th_refdom}.png};

\nextgroupplot[xtick=\empty, xlabel={}, ytick=\empty, ylabel={}]
\addplot graphics [xmin=-0.4, xmax=1.0, ymin=0, ymax=0.8] {{_img/visBurgs_pdeg_0th_refdom}.png};

\nextgroupplot[xtick=\empty, xlabel={}, ytick=\empty, ylabel={}]
\addplot graphics [xmin=-0.4, xmax=1.0, ymin=0, ymax=0.8] {{_img/visBurgs_u_1st_refdom}.png};

\nextgroupplot[xtick=\empty, xlabel={}, ytick=\empty, ylabel={}]
\addplot graphics [xmin=-0.4, xmax=1.0, ymin=0, ymax=0.8] {{_img/visBurgs_error_1st_refdom}.png};

\nextgroupplot[xtick=\empty, xlabel={}, ytick=\empty, ylabel={}]
\addplot graphics [xmin=-0.4, xmax=1.0, ymin=0, ymax=0.8] {{_img/visBurgs_pdeg_1st_refdom}.png};

\nextgroupplot[xtick=\empty, xlabel={}, ytick=\empty, ylabel={}]
\addplot graphics [xmin=-0.4, xmax=1.0, ymin=0, ymax=0.8] {{_img/visBurgs_u_2nd_refdom}.png};

\nextgroupplot[xtick=\empty, xlabel={}, ytick=\empty, ylabel={}]
\addplot graphics [xmin=-0.4, xmax=1.0, ymin=0, ymax=0.8] {{_img/visBurgs_error_2nd_refdom}.png};

\nextgroupplot[xtick=\empty, xlabel={}, ytick=\empty, ylabel={}]
\addplot graphics [xmin=-0.4, xmax=1.0, ymin=0, ymax=0.8] {{_img/visBurgs_pdeg_2nd_refdom}.png};

\nextgroupplot[xtick=\empty, xlabel={}, ytick=\empty, ylabel={}]
\addplot graphics [xmin=-0.4, xmax=1.0, ymin=0, ymax=0.8] {{_img/visBurgs_u_3rd_refdom}.png};

\nextgroupplot[xtick=\empty, xlabel={}, ytick=\empty, ylabel={}]
\addplot graphics [xmin=-0.4, xmax=1.0, ymin=0, ymax=0.8] {{_img/visBurgs_error_3rd_refdom}.png};

\nextgroupplot[xtick=\empty, xlabel={}, ytick=\empty, ylabel={}]
\addplot graphics [xmin=-0.4, xmax=1.0, ymin=0, ymax=0.8] {{_img/visBurgs_pdeg_3rd_refdom}.png};

\end{groupplot}

\end{tikzpicture}
	\caption{
		HOIST solution ($rp$-adaptation) (\textit{left}), the enriched residual
		error indicator (\textit{middle}), and polynomial degree distribution
		(\textit{right}) \textit{in the reference domain} to the viscous Burgers'
		problem with curved shock. These quantities are provided prior to
                $p$-adaptation and after each $p$-adaptation iteration
		(\textit{top-to-bottom}). Colorbar in Figure~\ref{fig:vburg2_sptm_physdom}.
	}
	\label{fig:vburg2_sptm_refdom}
\end{figure}
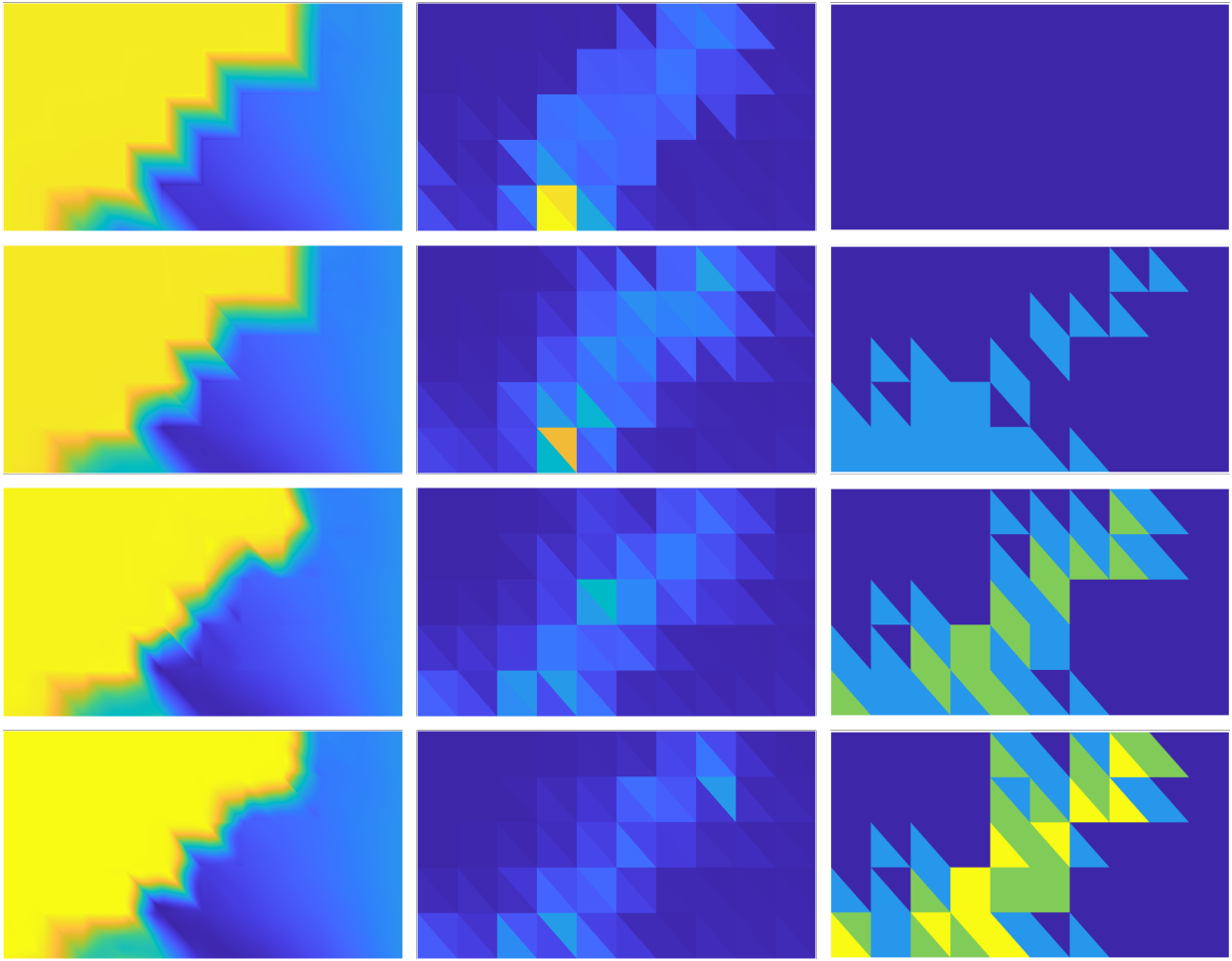

\subsection{Navier-Stokes equations}
\label{subsec:ns}
Next, we consider steady compressible, viscous flow through a domain
$\Omega \subset \Rbb^d$  with governing equations given by
\begin{equation}\label{eqn:nsii2d}
\begin{split}
\pder{}{x_j}\left(\rho v_j \right) & = 0 \\
\pder{}{x_j}\left(\rho v_i v_j  + P\delta_{ij}\right) & = \pder{}{x_j} \varpi_{ij}  \\
\pder{}{x_j}\left(\left(\rho E+P\right)v_j\right) &= \pder{}{x_j} \left(\varpi_{ij} v_i - q_j\right)
\end{split}
\end{equation}
for $i = 1, \cdots, d$. The density of the fluid $\func{\rho}{\Omega}{\Rbb_{>0}}$,
the fluid velocity $\func{v}{\Omega}{\Rbb^d}$, and the total energy of the fluid
$\func{\rho E}{\Omega}{\Rbb_{>0}}$ are implicitly defined
as the solution of (\ref{eqn:nsii2d}). We assume the fluid is an ideal gas
with thermal and caloric state equations
\begin{equation}\label{eqn:press}
P = \rho R T, \quad e = \frac{1}{\gamma - 1} \frac{P}{\rho},
\end{equation}
where $P : \Omega \rightarrow \Rbb_{>0}$ and $T : \Omega \rightarrow \Rbb_{>0}$
are the pressure and temperature of the fluid, $e : \Omega \rightarrow \Rbb_{>0}$
is the specific internal energy, $R \in \Rbb_{>0}$ is the specific gas constant,
and $\gamma > 1$ is the ratio of specific heats. The total energy is the sum
of the internal energy and the kinetic energy
\begin{equation}\label{eqn:energy}
\rho E = \rho e + \frac{1}{2}\rho v_i v_i.
\end{equation}
For a calorically perfect gas, the pressure is directly related to the conservative
variables as
\begin{equation}\label{eqn:press2}
P = \left(\gamma - 1 \right ) \left(\rho E - \rho v_i v_i/2\right).
\end{equation}
Assuming a Newtonian fluid with heat transfer that obeys Fourier's
law, the shear stress,
$\func{\varpi}{\Omega}{\Rbb^{d\times d}}$,
and heat flux,
$\func{q}{\Omega}{\Rbb^d}$, are related to
the primitive flow variables as
\begin{equation}
 \varpi_{ij} = \mu\left(\pder{v_i}{x_j}+\pder{v_j}{x_i}-\frac{2}{3}\pder{v_k}{x_k}\delta_{ij}\right), \qquad
 q_i = -\kappa \pder{T}{x_i},
\end{equation}
where $\func{\mu}{\Omega}{\Rbb_{>0}}$ is the
dynamic viscosity (independent of the fluid state by assumption), and
$\func{\kappa}{\Omega}{\Rbb_{>0}}$ is the
thermal conductivity. We use a constant viscosity model $\mu = \mathrm{Pr}\cdot\kappa/c_P$,
where $\mathrm{Pr}\in\Rbb_{>0}$ is the Prandl number and $c_P\in\Rbb_{>0}$ is the specific
heat capacity at constant pressure.

\subsubsection{Flow over a laminar flat plate: $rp$- vs. $h$-adaptation study}
\label{subsubsec:fp}
In this problem, we consider flow over a flat plate to demonstrate the ability of
the $rp$-adaptive method to resolve boundary layers, and compare it side-by-side
to a high-order, anisotropic $h$-adaptation method. For this study,
we consider a nondimensional version of the flat plate problem from \cite{economon2016su2},
nondimensionalized with respect to static inlet conditions (density, sound
speed, temperature) and the length of the flat plate.
The inlet Mach number is $M_\infty = 0.2$, the target Reynolds number
is $\mathrm{Re} = 1301233.166$, the ratio of specific heats is $\gamma = 1.4$,
and the Prandtl number is $\mathrm{Pr} = 0.72$.
The boundary conditions are summarized in Figure~\ref{fig:ns0_geom}.
The flow angle, stagnation temperature, and pressure are prescribed
at the inlet ($v/\| v \|_2 = (1, 0)$, $P_0 = 0.7345$, $T_0 = 1.008$),
the static pressure is prescribed at the outlet ($P_\mathrm{out} = 0.7143$),
and the flat plate is modeled as an adiabatic no-slip condition.
Additionally, a homogeneous normal viscous flux is imposed on
the inlet and outlet boundaries.

We discretize the domain using a structured, uniform $30\times 10$ grid 
of straight-sided ($q = 1$) quadrilateral elements with quadratic solution
approximation ($p_0 : K \mapsto 2$). The HOIST solver is initialized
from the DG solution on this coarse discretization (Figure~\ref{fig:ns0_soln_ic}).
To ensure the boundary layer is properly emphasized, we use the element
scaling in Section~\ref{subsubsec:resscale} with $\lambda = 100$. Reynolds number
continuation is performed in $15$ stages of equal increments beginning with
$\Xi_1 = 10^4$ and ending with $\Xi_{15} = 1301233.166$ with each continuation
stage allotted $20$ SQP iterations. We use step length modifications based
on the density and pressure with $\theta_\mathrm{l} = 0.1$ and $\theta_\mathrm{u} = 10$
(Section~\ref{subsubsec:stepcon}); PTC-based step modifications with
$\epsilon = 0.9$, $n_\mathrm{ptc} = 5$, and $\omega_\mathrm{mod}=5$
(Section~\ref{subsubsec:stepmod}); and enriched residual-based Hessian regularization
with $C_\mathrm{l} = 1$, $C_\mathrm{u} = 10$, $\eta_1 = 0.8$, $\eta_2 = 1$,
and $\hat\gamma = 10^4$
(Section~\ref{subsubsec:hessreg}). Finally, the enrichment degree is $\Delta = 2$ and
we perform four rounds of $p$-refinement using the dual-weighted residual error indicator
based on the drag on the flat plate.

\begin{figure}[!htbp]
	\centering
	\input{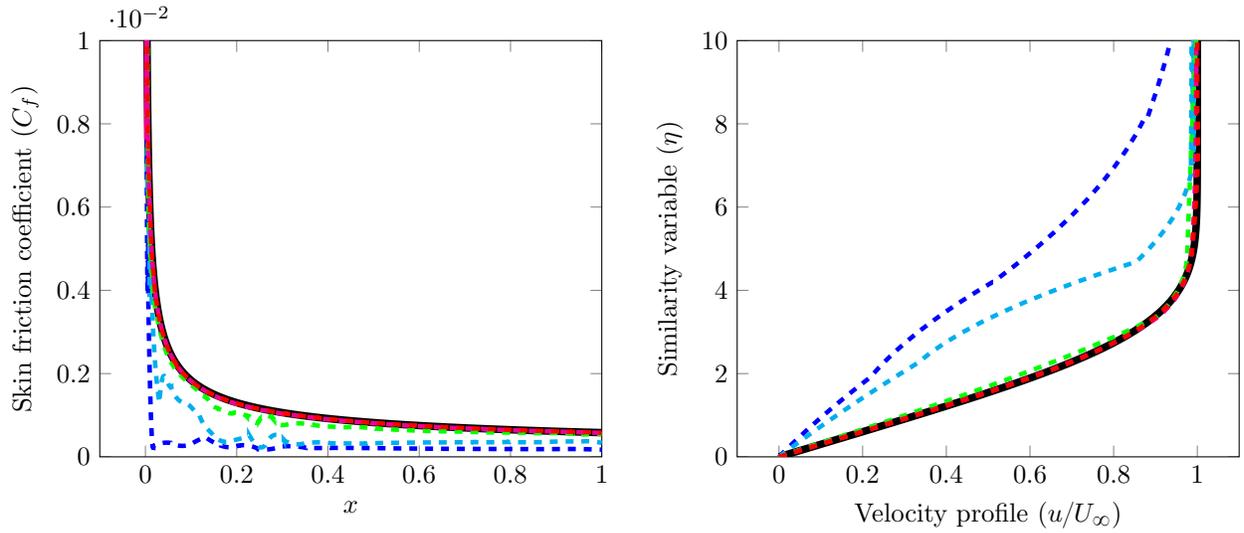}
	\caption{
		The skin friction coefficient (\textit{left}) and boundary layer
		profile (\textit{right}) computed from the Blasius solution
		(\ref{line:Cf_Blasius}) and the HOIST solution ($rp$-adaptation)
		prior to $p$-adaptation (\ref{line:fp 0th cycle}) and after the first
		(\ref{line:fp 1st cycle}), second (\ref{line:fp 2nd cycle}), third (\ref{line:fp 3rd cycle}), and fourth (\ref{line:fp 4th cycle})
		$p$-adaptation iterations for the flat plate problem.
	}
	\label{fig:ns0_blas}
\end{figure}

During the continuation procedure, mesh elements are compressed toward the
leading edge singularity and two layers of elements are squeezed near the
boundary layer. The $p$-refinement process concentrates the highest degrees
near the leading edge singularity and the beginning of the flat plate
where the boundary layer is thinnest. The boundary layer at the end of
the plate is well-resolved with a single element of moderate degree
(Figure~\ref{fig:ns0_soln_physdom}). From this figure, we also see
the $p$-refinement process effectively reduces the magnitude of the
error indicator across the domain. The internal structure of the
boundary layer, as well as the distribution of error and polynomial
degree, are more clear when visualized on the reference domain
$\Omega_0$ because the boundary layer is expanded
(Figure~\ref{fig:ns0_soln_refdom}).
From these figures, we see the proposed $rp$-adaptive method successfully
morphs a uniform mesh into one with substantial refinement, both in terms
of element size and polynomial degree, in the boundary layer and near the
leading edge singularity. Furthermore, after four $p$-refinement cycles,
the boundary layer profile and skin friction coefficient closely agree
with the Blasius solution (Figure~\ref{fig:ns0_blas}) (the profile agrees
after only three adaptation cycles).

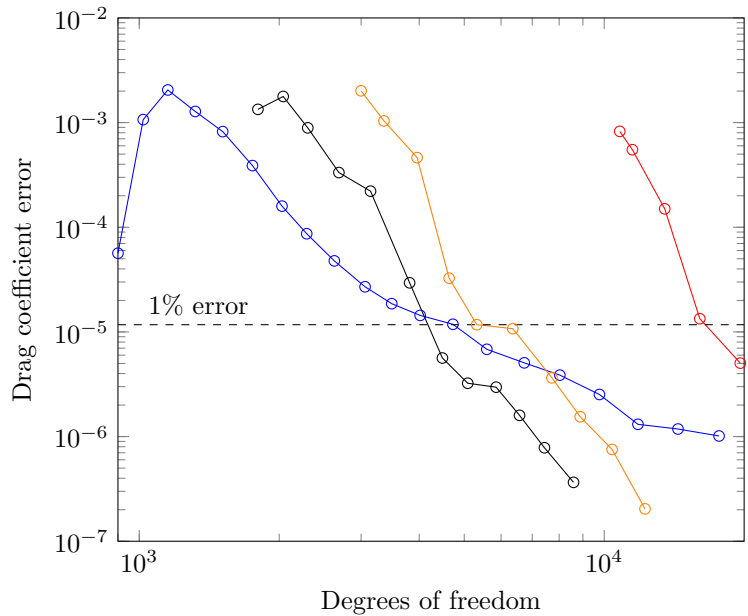
\begin{figure}[!htbp]
	\centering
	\begin{tikzpicture}
    \begin{loglogaxis}[width=0.6\textwidth, xmin=9.00000000e+02, xmax=2e4, xlabel={Degrees of freedom}, ymin=1e-7, ymax = 1e-2, ylabel={Drag coefficient error}, legend style={draw=none,}, legend pos=south west]

    \addplot [color=blue, mark=o]
    coordinates {
    ( 9.00000000e+02,  5.64838950e-05)
    ( 1.02000000e+03,  1.06830792e-03)
    ( 1.15200000e+03,  2.05060756e-03)
    ( 1.32000000e+03,  1.27755197e-03)
    ( 1.51200000e+03,  8.20149536e-04)
    ( 1.75200000e+03,  3.88319384e-04)
    ( 2.02800000e+03,  1.59180204e-04)
    ( 2.29200000e+03,  8.67435531e-05)
    ( 2.62800000e+03,  4.77257390e-05)
    ( 3.06000000e+03,  2.69957201e-05)
    ( 3.49200000e+03,  1.86201488e-05)
    ( 4.02000000e+03,  1.43852730e-05)
    ( 4.72800000e+03,  1.18304807e-05)
    ( 5.59200000e+03,  6.82281583e-06)
    ( 6.73200000e+03,  5.06565393e-06)
    ( 8.02800000e+03,  3.86794865e-06)
    ( 9.76800000e+03,  2.52518904e-06)
    ( 1.18200000e+04,  1.30991383e-06)
    ( 1.44120000e+04,  1.18111782e-06)
    ( 1.76640000e+04,  1.01356018e-06)
};\label{line:fp_p1};

    \addplot [color=black, mark=o]
    coordinates {
    ( 1.80000000e+03,  1.33766141e-03)
    ( 2.04000000e+03,  1.77732794e-03)
    ( 2.30400000e+03,  8.90707435e-04)
    ( 2.68800000e+03,  3.33218754e-04)
    ( 3.14400000e+03,  2.21251249e-04)
    ( 3.81600000e+03,  2.94550860e-05)
    ( 4.48800000e+03,  5.63372731e-06)
    ( 5.08800000e+03,  3.23402203e-06)
    ( 5.85600000e+03,  2.96578322e-06)
    ( 6.57600000e+03,  1.59593659e-06)
    ( 7.44000000e+03,  7.82234438e-07)
    ( 8.59200000e+03,  3.65293179e-07)};\label{line:fp_p2};

    \addplot [color=orange, mark=o]
    coordinates {
    ( 3.00000000e+03,  2.01359861e-03)
    ( 3.36000000e+03,  1.03793210e-03)
    ( 3.96000000e+03,  4.62606220e-04)
    ( 4.64000000e+03,  3.26132101e-05)
    ( 5.32000000e+03,  1.17257994e-05)
    ( 6.36000000e+03,  1.07426809e-05)
    ( 7.72000000e+03,  3.64207997e-06)
    ( 8.88000000e+03,  1.54673269e-06)
    ( 1.04000000e+04,  7.53198706e-07)
    ( 1.22400000e+04,  2.03789952e-07)};\label{line:fp_p3};

    \addplot [red, mark=o ] coordinates {
    (10800,   0.0008246922367959608)
    (11500,   0.0005503723834461753)
    (13496,   0.0001499810346767968)
    (16056,   1.3388521940094905e-5)
    (19624,   5.033876888964426e-6)
    };\label{line:fp_ist};
    
    \addplot[dashed, black] coordinates {
    (9e+02,  1e-2*1.171809938287728e-03)   (2e4, 1e-2*1.171809938287728e-03)
   };
   \node[anchor=south west, black] at (axis cs:1e3, 1e-2*1.171809938287728e-03) {$1\%$ error};
    \end{loglogaxis}
\end{tikzpicture}
    	\caption{
		Convergence of the error in the drag coefficient for the
		HOIST method ($rp$-adaptation) (\ref{line:fp_ist}) and the $h$-adaptive
		method with linear (\ref{line:fp_p1}), quadratic (\ref{line:fp_p2}),
		and cubic (\ref{line:fp_p3}) constant polynomial degree for the flat plate
		problem. Both flow and mesh degrees of freedom are counted
		the HOIST method ($N_\ubm^{(j)} + N_\xbm$), while only flow
		degrees of freedom are present for the $h$-adaptive method.
	}
	\label{fig:ns0_cderr}
\end{figure}

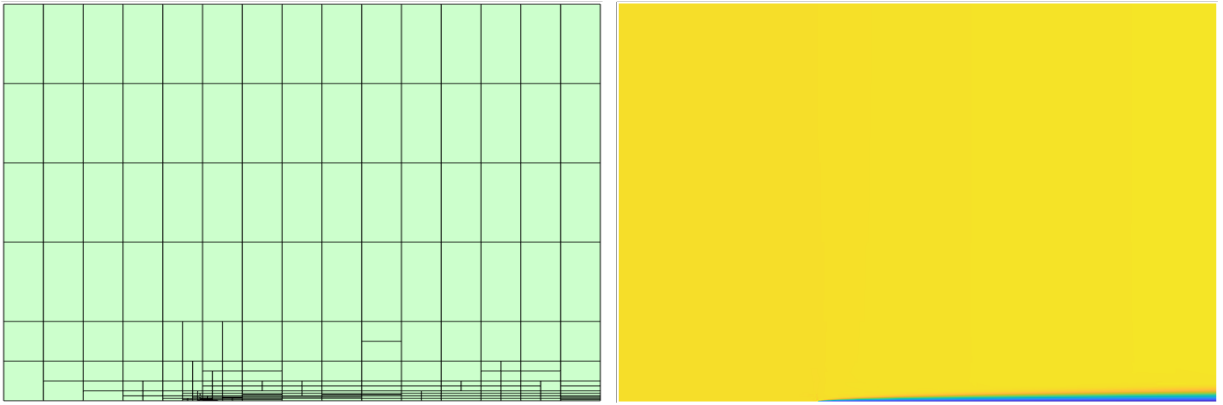
\begin{figure}[!htbp]
	\centering
	\begin{tikzpicture}[scale=0.8]
\begin{groupplot}[
  group style={
      group size=2 by 1,
      horizontal sep=0.25cm,
      vertical sep=0.25cm
  },
  width=0.7\textwidth,
  axis equal image,
  xmin=-0.5, xmax=1.0,
  ymin=0, ymax=1
]

\nextgroupplot[xtick=\empty, xlabel={}, ytick=\empty, ylabel={}]
\addplot graphics [xmin=-0.5, xmax=1.0, ymin=0, ymax=1] {{_img/fp_masa_msh}.png};

\nextgroupplot[xtick=\empty, xlabel={}, ytick=\empty, ylabel={}]
\addplot graphics [xmin=-0.5, xmax=1.0, ymin=0, ymax=1] {{_img/fp_masa_Mach_womsh}.png};

\end{groupplot}
\end{tikzpicture}
	\caption{
		The $h$-adapted mesh (\textit{left}) and corresponding solution
		(Mach number) (\textit{right}) for the flat plate problem
		using drag-based dual-weighted residual adaptation
		used for the comparative study with the
		HOIST method ($rp$-adaptation).
	}
	\label{fig:ns0_amr}
\end{figure}

Finally, we compare the proposed $rp$-adaptive method to an anisotropic 
$h$-adaptive method. The $h$-adaptive method uses $p = 1$, $2$, and $3$ discretizations,
with the initial uniform mesh comprised $15 \times 5$ elements.
For both methods, error adaptation is driven by the dual-weighted residual error estimate
based on the flat plate drag. Both methods reach a solution with sub-$1\%$ error;
however, the proposed $rp$-adaptive method requires more degrees of freedom than
the $h$-adaptive method to achieve such a solution for all polynomial
degrees considered (Figure~\ref{fig:ns0_cderr}).
A representative example of the $h$-adapted mesh
and solution produced by the $h$-adaptive method is provided in
Figure~\ref{fig:ns0_amr}. While this method demonstrated the flexibility
of the proposed method to automatically adapt a uniform mesh to boundary
layers, it is not necessarily competitive with other approaches. In the
next section, we show that it is highly competitive for flows involving
strong shocks.

\begin{figure}[!htbp]
        \centering
        \begin{tikzpicture}
\begin{axis}[
width=0.55\textwidth,
axis line style={gray},
axis x line*=bottom,
axis y line*=left,
xtick={-0.5, 0.0, 1.0},
ytick={0.0, 0.05, 0.1},
ymax=0.11000000000000001,
xmax=1.15,
xmin=-0.65,
ymin=-0.010000000000000002]
\addplot [opacity=0.6, fill=black!30!white, opacity=0.6, forget plot]
coordinates {
(-5.00000000e-01,  0.00000000e+00)
( 0.00000000e+00,  0.00000000e+00)
( 1.00000000e+00,  0.00000000e+00)
( 1.00000000e+00,  1.00000000e-01)
(-5.00000000e-01,  1.00000000e-01)
(-5.00000000e-01,  0.00000000e+00)};

\addplot [thick, color=green]
coordinates {
(-5.00000000e-01,  0.00000000e+00)
( -5.00000000e-01,  1.00000000e-01)};\label{line:fp:inlet}

\addplot [thick, color=black]
coordinates {
(-5.00000000e-01,  0.00000000e+00)
( 0.00000000e+00,  0.00000000e+00)};\label{line:fp:symm}

\addplot [thick, color=red]
coordinates {
( 0.00000000e+00,  0.00000000e+00)
( 1.00000000e+00, 0.00000000e+00)};\label{line:fp:noslipwall}

\addplot [thick, color=blue]
coordinates {
( 1.00000000e+00,  0.00000000e+00)
( 1.00000000e+00,  1.00000000e-01)};\label{line:fp:outflow}

\addplot [thick, color=blue]
coordinates {
( -5.00000000e-01, 1.00000000e-01)
( 1.00000000e+00,  1.00000000e-01)};

\end{axis}
\end{tikzpicture}
        \caption{
		Schematic of flat plate domain and boundary conditions. Boundary
		conditions: subsonic inlet (\ref{line:fp:inlet}), subsonic outlet
		(\ref{line:fp:outflow}), adiabatic no-slip wall (\ref{line:fp:noslipwall}),
		symmetry condition (\ref{line:fp:symm}).
        }
        \label{fig:ns0_geom}
\end{figure}
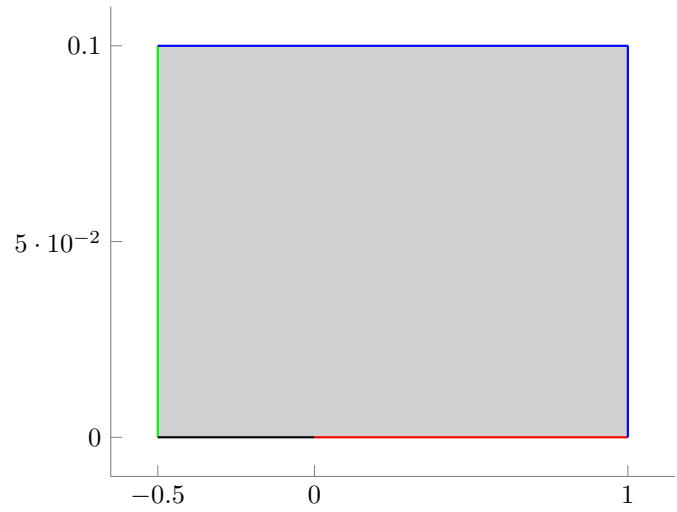

\begin{figure}[!htbp]
	\centering
	\begin{tikzpicture}
\begin{groupplot}[
  group style={
      group size=1 by 1,
      horizontal sep=1cm,
      vertical sep=0.5cm
  },
  width=0.5\textwidth,
  axis equal image,
  xlabel={$x_1$},
  ylabel={$x_2$},
  xtick = {-0.5, 0.0, 0.5, 1.0},
  xticklabels={-0.5, 0.0, 0.5, 1.0},
  ytick = {0.0, 0.5, 1.0},
  yticklabels={0.0, 0.1},
  xmin=-0.5, xmax=1.0,
  ymin=0, ymax=1
]

\nextgroupplot[xticklabels={-0.5, 0.0, 0.5, 1.0}, yticklabels={0.0, 0.05, 0.1}]
\addplot graphics [xmin=-0.5, xmax=1.0, ymin=0, ymax=1.0] {{_img/flatplateRe1e4_initialMesh}.png};

\end{groupplot}
\node[anchor=north] at ($(group c1r1.south) + (0.3cm, -0.5) $) {\colorbarMatlabParula{0}{0.05}{0.1}{0.15}{0.2}};
\end{tikzpicture}
	\caption{
		Mesh and solution (DG solution with $\mathrm{Re = 10^4}$) used to
		initialize the $rp$-adaptive HOIST method for the flat plate problem.
	}
	\label{fig:ns0_soln_ic}
\end{figure}
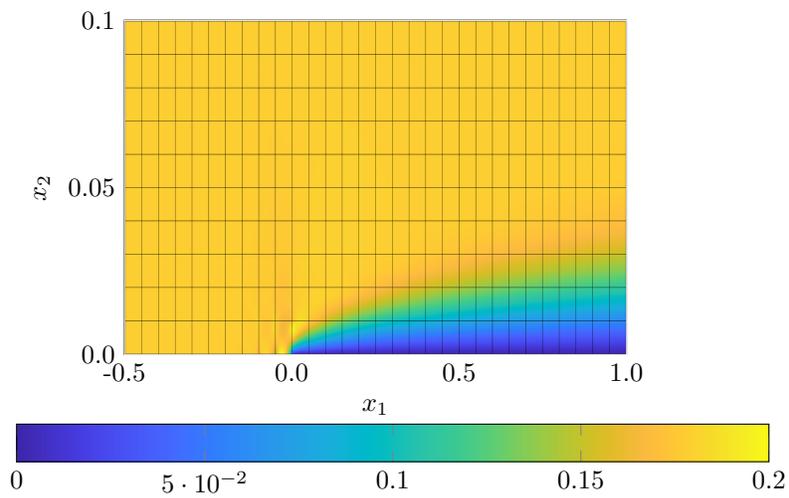

\begin{figure}[!htbp]
	\centering
	\begin{tikzpicture}[scale=0.8]
\begin{groupplot}[
  group style={
      group size=3 by 5,
      horizontal sep=0.25cm,
      vertical sep=0.25cm
  },
  width=0.5\textwidth,
  axis equal image,
  xmin=-0.5, xmax=1.0,
  ymin=0, ymax=1
]

\nextgroupplot[xtick=\empty, xlabel={}, ytick=\empty, ylabel={}]
\addplot graphics [xmin=-0.5, xmax=1.0, ymin=0, ymax=1] {{_img/fp_laminar_Mach_0th_phydom_mesh}.png};

\nextgroupplot[xtick=\empty, xlabel={}, ytick=\empty, ylabel={}]
\addplot graphics [xmin=-0.5, xmax=1.0, ymin=0, ymax=1] {{_img/fp_laminar_log10Error_0th_phydom}.png};

\nextgroupplot[xtick=\empty, xlabel={}, ytick=\empty, ylabel={}]
\addplot graphics [xmin=-0.5, xmax=1.0, ymin=0, ymax=1] {{_img/fp_laminar_pdeg_0th_phydom}.png};

\nextgroupplot[xtick=\empty, xlabel={}, ytick=\empty, ylabel={}]
\addplot graphics [xmin=-0.5, xmax=1.0, ymin=0, ymax=1] {{_img/fp_laminar_Mach_1st_phydom_mesh}.png};

\nextgroupplot[xtick=\empty, xlabel={}, ytick=\empty, ylabel={}]
\addplot graphics [xmin=-0.5, xmax=1.0, ymin=0, ymax=1] {{_img/fp_laminar_log10Error_1st_phydom}.png};

\nextgroupplot[xtick=\empty, xlabel={}, ytick=\empty, ylabel={}]
\addplot graphics [xmin=-0.5, xmax=1.0, ymin=0, ymax=1] {{_img/fp_laminar_pdeg_1st_phydom}.png};

\nextgroupplot[xtick=\empty, xlabel={}, ytick=\empty, ylabel={}]
\addplot graphics [xmin=-0.5, xmax=1.0, ymin=0, ymax=1] {{_img/fp_laminar_Mach_2nd_phydom_mesh}.png};

\nextgroupplot[xtick=\empty, xlabel={}, ytick=\empty, ylabel={}]
\addplot graphics [xmin=-0.5, xmax=1.0, ymin=0, ymax=1] {{_img/fp_laminar_log10Error_2nd_phydom}.png};

\nextgroupplot[xtick=\empty, xlabel={}, ytick=\empty, ylabel={}]
\addplot graphics [xmin=-0.5, xmax=1.0, ymin=0, ymax=1] {{_img/fp_laminar_pdeg_2nd_phydom}.png};

\nextgroupplot[xtick=\empty, xlabel={}, ytick=\empty, ylabel={}]
\addplot graphics [xmin=-0.5, xmax=1.0, ymin=0, ymax=1] {{_img/fp_laminar_Mach_3rd_phydom_mesh}.png};

\nextgroupplot[xtick=\empty, xlabel={}, ytick=\empty, ylabel={}]
\addplot graphics [xmin=-0.5, xmax=1.0, ymin=0, ymax=1] {{_img/fp_laminar_log10Error_3rd_phydom}.png};

\nextgroupplot[xtick=\empty, xlabel={}, ytick=\empty, ylabel={}]
\addplot graphics [xmin=-0.5, xmax=1.0, ymin=0, ymax=1] {{_img/fp_laminar_pdeg_3rd_phydom}.png};

\nextgroupplot[xtick=\empty, xlabel={}, ytick=\empty, ylabel={}]
\addplot graphics [xmin=-0.5, xmax=1.0, ymin=0, ymax=1] {{_img/fp_laminar_Mach_4th_phydom_mesh}.png};

\nextgroupplot[xtick=\empty, xlabel={}, ytick=\empty, ylabel={}]
\addplot graphics [xmin=-0.5, xmax=1.0, ymin=0, ymax=1] {{_img/fp_laminar_log10Error_4th_phydom}.png};

\nextgroupplot[xtick=\empty, xlabel={}, ytick=\empty, ylabel={}]
\addplot graphics [xmin=-0.5, xmax=1.0, ymin=0, ymax=1] {{_img/fp_laminar_pdeg_4th_phydom}.png};

\end{groupplot}
\node[anchor=west] at ($(group c1r5.west) + (-0.3cm, -3)$) {\colorbarMatlabParulandscalse{0}{0.05}{0.1}{0.15}{0.2}{0.2}{5cm}};

\node[anchor=west] at ($(group c2r5.west) + (-0.28cm, -3)$) {\colorbarMatlabParulandscalse{-17}{-14}{-11}{-8}{-5}{-2}{5cm}};

\node[anchor=west] at ($(group c3r5.west) + (-0.1cm, -3)$) {\colorbarMatlabParulandscalse{2}{3}{4}{5}{6}{6}{5cm}};
\end{tikzpicture}
	\caption{
		HOIST solution ($rp$-adaptation) (Mach number) (\textit{left}),
		the drag-based dual-weighted residual error estimate (\textit{middle}),
		and the polynomial degree distribution (\textit{right}) \textit{in the
		physical domain} $\Gcal(\Omega_0)$ to the flat plate problem.
		These quantities are provided prior to $p$-adaptation and after
		each $p$-adaptation iteration (\textit{top-to-bottom}).
	}
	\label{fig:ns0_soln_physdom}
\end{figure}

\begin{figure}[!htbp]
	\centering
	\begin{tikzpicture}[scale=0.8]
\begin{groupplot}[
  group style={
      group size=3 by 5,
      horizontal sep=0.25cm,
      vertical sep=0.25cm
  },
  width=0.5\textwidth,
  axis equal image,
  xmin=-0.5, xmax=1.0,
  ymin=0, ymax=1
]

\nextgroupplot[xtick=\empty, xlabel={}, ytick=\empty, ylabel={}]
\addplot graphics [xmin=-0.5, xmax=1.0, ymin=0, ymax=1] {{_img/fp_laminar_Mach_0th_refdom}.png};

\nextgroupplot[xtick=\empty, xlabel={}, ytick=\empty, ylabel={}]
\addplot graphics [xmin=-0.5, xmax=1.0, ymin=0, ymax=1] {{_img/fp_laminar_log10Error_0th_refdom}.png};

\nextgroupplot[xtick=\empty, xlabel={}, ytick=\empty, ylabel={}]
\addplot graphics [xmin=-0.5, xmax=1.0, ymin=0, ymax=1] {{_img/fp_laminar_pdeg_0th_refdom}.png};

\nextgroupplot[xtick=\empty, xlabel={}, ytick=\empty, ylabel={}]
\addplot graphics [xmin=-0.5, xmax=1.0, ymin=0, ymax=1] {{_img/fp_laminar_Mach_1st_refdom}.png};

\nextgroupplot[xtick=\empty, xlabel={}, ytick=\empty, ylabel={}]
\addplot graphics [xmin=-0.5, xmax=1.0, ymin=0, ymax=1] {{_img/fp_laminar_log10Error_1st_refdom}.png};

\nextgroupplot[xtick=\empty, xlabel={}, ytick=\empty, ylabel={}]
\addplot graphics [xmin=-0.5, xmax=1.0, ymin=0, ymax=1] {{_img/fp_laminar_pdeg_1st_refdom}.png};

\nextgroupplot[xtick=\empty, xlabel={}, ytick=\empty, ylabel={}]
\addplot graphics [xmin=-0.5, xmax=1.0, ymin=0, ymax=1] {{_img/fp_laminar_Mach_2nd_refdom}.png};

\nextgroupplot[xtick=\empty, xlabel={}, ytick=\empty, ylabel={}]
\addplot graphics [xmin=-0.5, xmax=1.0, ymin=0, ymax=1] {{_img/fp_laminar_log10Error_2nd_refdom}.png};

\nextgroupplot[xtick=\empty, xlabel={}, ytick=\empty, ylabel={}]
\addplot graphics [xmin=-0.5, xmax=1.0, ymin=0, ymax=1] {{_img/fp_laminar_pdeg_2nd_refdom}.png};

\nextgroupplot[xtick=\empty, xlabel={}, ytick=\empty, ylabel={}]
\addplot graphics [xmin=-0.5, xmax=1.0, ymin=0, ymax=1] {{_img/fp_laminar_Mach_3rd_refdom}.png};

\nextgroupplot[xtick=\empty, xlabel={}, ytick=\empty, ylabel={}]
\addplot graphics [xmin=-0.5, xmax=1.0, ymin=0, ymax=1] {{_img/fp_laminar_log10Error_3rd_refdom}.png};

\nextgroupplot[xtick=\empty, xlabel={}, ytick=\empty, ylabel={}]
\addplot graphics [xmin=-0.5, xmax=1.0, ymin=0, ymax=1] {{_img/fp_laminar_pdeg_3rd_refdom}.png};

\nextgroupplot[xtick=\empty, xlabel={}, ytick=\empty, ylabel={}]
\addplot graphics [xmin=-0.5, xmax=1.0, ymin=0, ymax=1] {{_img/fp_laminar_Mach_4th_refdom}.png};

\nextgroupplot[xtick=\empty, xlabel={}, ytick=\empty, ylabel={}]
\addplot graphics [xmin=-0.5, xmax=1.0, ymin=0, ymax=1] {{_img/fp_laminar_log10Error_4th_refdom}.png};

\nextgroupplot[xtick=\empty, xlabel={}, ytick=\empty, ylabel={}]
\addplot graphics [xmin=-0.5, xmax=1.0, ymin=0, ymax=1] {{_img/fp_laminar_pdeg_4th_refdom}.png};

\end{groupplot}

\end{tikzpicture}
	\caption{
                HOIST solution ($rp$-adaptation) (Mach number) (\textit{left}),
                the drag-based dual-weighted residual error estimate (\textit{middle}),
		and the polynomial degree distribution (\textit{right}) \textit{in the
                reference domain} to the flat plate problem. These quantities are
                provided prior to $p$-adaptation and after each $p$-adaptation
		iteration (\textit{top-to-bottom}). Colorbar in
		Figure~\ref{fig:ns0_soln_physdom}.
        }
	\label{fig:ns0_soln_refdom}
\end{figure}
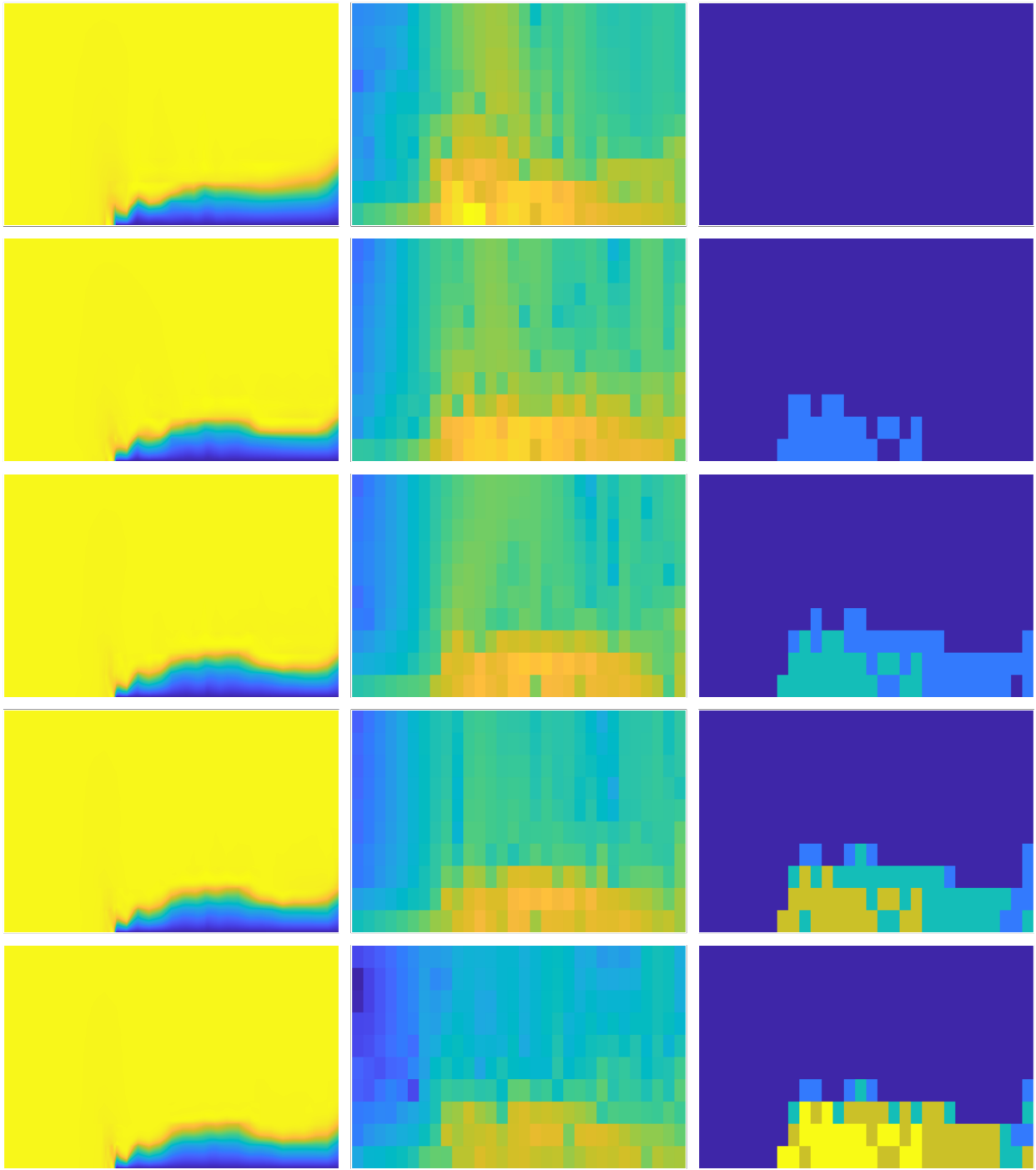

\subsubsection{Hypersonic flow over cylinder: $rp$- vs. $h$-adaptation study}
\label{subsubsec:bow}
In this problem, we consider hypersonic flow over a cylinder to demonstrate the
ability of the $rp$-adaptive method to accurately resolve a strong, thin bow
shock and thin boundary layer, and compare it side-by-side to a high-order,
anisotropic $h$-adaptive method. For this study, we consider
the problem in \cite{kercher2021moving} in a slightly modified computational
geometry with cylinder wall defined as the unit circle and the farfield
boundary defined as the ellipse $(x_1/6)^2 + (x_2/3)^2 = 1$. The inlet
Mach number is $M_\infty = 5$, the target Reynolds number based on the
radius of the cylinder is $\mathrm{Re} = 1000$,
the ratio of specific heats is $\gamma = 1.4$, the Prandtl number
is $\mathrm{Pr} = 0.72$, and the isothermal wall temperature is
$T_w/T_\infty = 2.5$. The geometry
and boundary conditions are illustrated in Figure~\ref{fig:ns1_geom}.
The entire primitive state is prescribed at the supersonic inlet and
the cylinder wall is modeled as an isothermal no-slip condition.
Additionally, a homogeneous normal viscous flux is imposed on
the inlet and outlet boundaries.

We discretize the domain using a structured, uniform mesh of $210$
curved ($q = 2$) quadrilateral elements with quadratic solution
approximation ($p_0 : K \mapsto 2$). The HOIST solver is initialized
from a shock capturing DG solution using PDE-based artificial viscosity
from \cite{barter2010shock} on this coarse discretization (Figure~\ref{fig:ns1_soln_ic}).
To ensure the boundary layer is properly emphasized, we use the element
scaling in Section~\ref{subsubsec:resscale} with $\lambda = 10$. Reynolds number
continuation is performed in $10$ stages of equal increments beginning with $\Xi_1 = 10^2$
and ending with $\Xi_{10} = 10^3$ with the first continuation stage allotted
$40$ SQP iterations and the remaining stages allowed $20$ SQP iterations.
We use step length modifications based on the density and pressure
with $\theta_\mathrm{l} = 0.1$ and $\theta_\mathrm{u} = 10$
(Section~\ref{subsubsec:stepcon}); PTC-based step modifications with
$n_\mathrm{ptc} = 5$, $\epsilon = 0.9$, and $\omega_\mathrm{mod}=10$
(Section~\ref{subsubsec:stepmod}); and enriched residual-based Hessian regularization
with $C_\mathrm{l} = 1$, $C_\mathrm{u} = 10$, $\eta_1 = \eta_2 = 1$,
and $\hat\gamma = 100$ (Section~\ref{subsubsec:hessreg}).
Finally, the enrichment degree $\Delta = 2$ is used and we perform
five rounds of $p$-refinement using the enriched residual error indicator.
After all SQP iterations are exhausted in the final $p$-adaptation iteration,
PTC is used to drive the DG residual to $\Ocal(10^{-9})$ in less than $20$ time
steps using the (now fixed) $rp$-adapted approximation space without any form of nonlinear stabilization.

During the continuation procedure, mesh elements are compressed toward the
bow shock and boundary layer with the compression being particularly tight
in the vicinity of the stagnation streamline. Interestingly, the $p$-refinement
process concentrates the highest polynomial degrees in the bow shock,
moderate degrees between the bow shock and boundary layer, and
lowest degrees in the boundary layer (Figure~\ref{fig:ns1_soln_physdom}).
The lower degrees in the boundary layer are likely due to the extremely
tight grid spacing (due to $r$-adaptation) making them sufficient to
resolve the feature, and from our choice to use the unweighted enriched
residual indicator. This figure also shows the $p$-refinement process
effectively reduces the magnitude of the error indicator across the domain.
The internal structure of the bow shock and boundary layer, as well as the
distribution of error and polynomial degree, are more clear when visualized
on the reference domain $\Omega_0$ because the thin features are expanded
(Figure~\ref{fig:ns1_soln_refdom}).
From these figures, we see the proposed $rp$-adaptive method successfully
morphs a structured mesh into one with substantial refinement, in terms
of element size and/or polynomial degree, near the shock and boundary
layer. Furthermore, after five $p$-refinement cycles, the heat flux
profile, a notoriously difficult quantity to predict, is well-resolved
(Figure~\ref{fig:ns1_hflx}); however, the combination of $r$- and
$p$-adaptivity was crucial as $r$-refinement alone and early
$p$-refinement iterations produce a highly oscillatory heat flux
profile with substantial inter-element jumps.
Finally, it is interesting to note that the heat flux computed from the
lifted viscous flux ($\sigma$) provides a smoother profile with smaller inter-element
jumps than direct differentiation (Figure~\ref{fig:ns1_hflux_liftvsnolift}); see
discussion in Remark~\ref{rem:liftgrad}. 

Next, we compare the proposed $rp$-adaptive method to a high-order
$h$-adaptive method with adaptation driven
by the dual-weighted residual error estimate based on the integrated heat flux
over the cylinder. The $h$-adaptive method uses $p = 1$, $2$, and $3$ discretizations and the same initial mesh as the $rp$-adaptive method.
Both methods achieve solutions with sub-$1\%$ error, and the
proposed $rp$-adaptive method requires fewer degrees of freedom than the $h$-adaptive
method to do so (Figure~\ref{fig:ns1_hflxerr}). This shows the HOIST method
with $rp$-adaptive method effectively resolves complex high-speed flow features
and achieves high accuracy per degree of freedom, making it a competitive approach
for such problems. A representative example of the $h$-adapted mesh
and solution produced by the $h$-adaptive method is provided in Figure~\ref{fig:ns1_amr}. The DWR method refines the shock only in the regions relevant to the heat flux evaluation.

The behavior of the SQP solver across all continuation and $p$-adaptation stages
is shown in Figure~\ref{fig:ns1_sqp}. For a given value of the viscosity
$\Xi_i$ and polynomial distribution $p_j$, the HOIST method drives the residual
toward zero. However, the residual jumps up after the continuation parameter or
polynomial distribution is updated because the residual function changes. Notice
that 260 iterations are performed on the coarsest (i.e., least expensive)
discretization that uses a quadratic approximations in all elements. The
$p$-adaptation process is complete in another 160 iterations, at which
point the mesh is effectively frozen. After the SQP iterations are complete,
PTC is used to deeply converge the residual on a fixed mesh to a tolerance
of $10^{-9}$.

The remainder of the this section studies the impact of the
various aspects of the complete $rp$-adaptive HOIST formulation and solver.

\ifbool{fastcompile}{}{

	\begin{figure}[!htbp]
	        \centering
            \begin{tikzpicture}
\begin{axis}[
width=0.65\textwidth,
axis equal image,
axis line style={gray},
axis x line*=bottom,
axis y line*=left,
xtick={-6, -1, 0, 1,6},
ytick={0, 3},
xlabel={$x_1$},
ymax=3.5,
xmax=6.5,
ylabel={$x_2$},
xmin=-6.5,
ymin=-0.5]
\addplot [blue, thick, solid,opacity=0.6, fill=black!30!white, forget plot]
coordinates {
(6.00000000e+00, 0.00000000e+00)
(5.99697925e+00, 9.51838005e-02)
(5.98792006e+00, 1.90271759e-01)
(5.97283154e+00, 2.85168130e-01)
(5.95172888e+00, 3.79777361e-01)
(5.92463333e+00, 4.74004188e-01)
(5.89157218e+00, 5.67753733e-01)
(5.85257872e+00, 6.60931598e-01)
(5.80769221e+00, 7.53443962e-01)
(5.75695784e+00, 8.45197671e-01)
(5.70042671e+00, 9.36100337e-01)
(5.63815572e+00, 1.02606043e+00)
(5.57020760e+00, 1.11498737e+00)
(5.49665074e+00, 1.20279161e+00)
(5.41755923e+00, 1.28938474e+00)
(5.33301269e+00, 1.37467957e+00)
(5.24309626e+00, 1.45859021e+00)
(5.14790048e+00, 1.54103217e+00)
(5.04752120e+00, 1.62192245e+00)
(4.94205949e+00, 1.70117959e+00)
(4.83162155e+00, 1.77872379e+00)
(4.71631857e+00, 1.85447696e+00)
(4.59626666e+00, 1.92836283e+00)
(4.47158670e+00, 2.00030700e+00)
(4.34240423e+00, 2.07023703e+00)
(4.20884933e+00, 2.13808251e+00)
(4.07105647e+00, 2.20377513e+00)
(3.92916440e+00, 2.26724872e+00)
(3.78331600e+00, 2.32843939e+00)
(3.63365812e+00, 2.38728552e+00)
(3.48034146e+00, 2.44372786e+00)
(3.32352038e+00, 2.49770956e+00)
(3.16335281e+00, 2.54917629e+00)
(3.00000000e+00, 2.59807621e+00)
(2.83362645e+00, 2.64436009e+00)
(2.66439968e+00, 2.68798132e+00)
(2.49249008e+00, 2.72889599e+00)
(2.31807075e+00, 2.76706288e+00)
(2.14131733e+00, 2.80244358e+00)
(1.96240778e+00, 2.83500246e+00)
(1.78152225e+00, 2.86470672e+00)
(1.59884288e+00, 2.89152648e+00)
(1.41455361e+00, 2.91543470e+00)
(1.22884001e+00, 2.93640734e+00)
(1.04188907e+00, 2.95442326e+00)
(8.53889030e-01, 2.96946433e+00)
(6.65029199e-01, 2.98151539e+00)
(4.75499741e-01, 2.99056433e+00)
(2.85491495e-01, 2.99660202e+00)
(9.51957830e-02, 2.99962238e+00)
(-9.51957830e-02, 2.99962238e+00)
(-2.85491495e-01, 2.99660202e+00)
(-4.75499741e-01, 2.99056433e+00)
(-6.65029199e-01, 2.98151539e+00)
(-8.53889030e-01, 2.96946433e+00)
(-1.04188907e+00, 2.95442326e+00)
(-1.22884001e+00, 2.93640734e+00)
(-1.41455361e+00, 2.91543470e+00)
(-1.59884288e+00, 2.89152648e+00)
(-1.78152225e+00, 2.86470672e+00)
(-1.96240778e+00, 2.83500246e+00)
(-2.14131733e+00, 2.80244358e+00)
(-2.31807075e+00, 2.76706288e+00)
(-2.49249008e+00, 2.72889599e+00)
(-2.66439968e+00, 2.68798132e+00)
(-2.83362645e+00, 2.64436009e+00)
(-3.00000000e+00, 2.59807621e+00)
(-3.16335281e+00, 2.54917629e+00)
(-3.32352038e+00, 2.49770956e+00)
(-3.48034146e+00, 2.44372786e+00)
(-3.63365812e+00, 2.38728552e+00)
(-3.78331600e+00, 2.32843939e+00)
(-3.92916440e+00, 2.26724872e+00)
(-4.07105647e+00, 2.20377513e+00)
(-4.20884933e+00, 2.13808251e+00)
(-4.34240423e+00, 2.07023703e+00)
(-4.47158670e+00, 2.00030700e+00)
(-4.59626666e+00, 1.92836283e+00)
(-4.71631857e+00, 1.85447696e+00)
(-4.83162155e+00, 1.77872379e+00)
(-4.94205949e+00, 1.70117959e+00)
(-5.04752120e+00, 1.62192245e+00)
(-5.14790048e+00, 1.54103217e+00)
(-5.24309626e+00, 1.45859021e+00)
(-5.33301269e+00, 1.37467957e+00)
(-5.41755923e+00, 1.28938474e+00)
(-5.49665074e+00, 1.20279161e+00)
(-5.57020760e+00, 1.11498737e+00)
(-5.63815572e+00, 1.02606043e+00)
(-5.70042671e+00, 9.36100337e-01)
(-5.75695784e+00, 8.45197671e-01)
(-5.80769221e+00, 7.53443962e-01)
(-5.85257872e+00, 6.60931598e-01)
(-5.89157218e+00, 5.67753733e-01)
(-5.92463333e+00, 4.74004188e-01)
(-5.95172888e+00, 3.79777361e-01)
(-5.97283154e+00, 2.85168130e-01)
(-5.98792006e+00, 1.90271759e-01)
(-5.99697925e+00, 9.51838005e-02)
(-6.00000000e+00, 3.67394040e-16)
};\label{line:cyl:inflow}

\addplot [red, thick, opacity=1.0, fill=white, solid, forget plot]
coordinates {
(1.00000000e+00, 0.00000000e+00)
(9.99496542e-01, 3.17279335e-02)
(9.97986676e-01, 6.34239197e-02)
(9.95471923e-01, 9.50560433e-02)
(9.91954813e-01, 1.26592454e-01)
(9.87438889e-01, 1.58001396e-01)
(9.81928697e-01, 1.89251244e-01)
(9.75429787e-01, 2.20310533e-01)
(9.67948701e-01, 2.51147987e-01)
(9.59492974e-01, 2.81732557e-01)
(9.50071118e-01, 3.12033446e-01)
(9.39692621e-01, 3.42020143e-01)
(9.28367933e-01, 3.71662456e-01)
(9.16108457e-01, 4.00930535e-01)
(9.02926538e-01, 4.29794912e-01)
(8.88835449e-01, 4.58226522e-01)
(8.73849377e-01, 4.86196736e-01)
(8.57983413e-01, 5.13677392e-01)
(8.41253533e-01, 5.40640817e-01)
(8.23676581e-01, 5.67059864e-01)
(8.05270258e-01, 5.92907929e-01)
(7.86053095e-01, 6.18158986e-01)
(7.66044443e-01, 6.42787610e-01)
(7.45264450e-01, 6.66769001e-01)
(7.23734038e-01, 6.90079011e-01)
(7.01474888e-01, 7.12694171e-01)
(6.78509412e-01, 7.34591709e-01)
(6.54860734e-01, 7.55749574e-01)
(6.30552667e-01, 7.76146464e-01)
(6.05609687e-01, 7.95761841e-01)
(5.80056910e-01, 8.14575952e-01)
(5.53920064e-01, 8.32569855e-01)
(5.27225468e-01, 8.49725430e-01)
(5.00000000e-01, 8.66025404e-01)
(4.72271075e-01, 8.81453363e-01)
(4.44066613e-01, 8.95993774e-01)
(4.15415013e-01, 9.09631995e-01)
(3.86345126e-01, 9.22354294e-01)
(3.56886222e-01, 9.34147860e-01)
(3.27067963e-01, 9.45000819e-01)
(2.96920375e-01, 9.54902241e-01)
(2.66473814e-01, 9.63842159e-01)
(2.35758936e-01, 9.71811568e-01)
(2.04806668e-01, 9.78802446e-01)
(1.73648178e-01, 9.84807753e-01)
(1.42314838e-01, 9.89821442e-01)
(1.10838200e-01, 9.93838464e-01)
(7.92499569e-02, 9.96854776e-01)
(4.75819158e-02, 9.98867339e-01)
(1.58659638e-02, 9.99874128e-01)
(-1.58659638e-02, 9.99874128e-01)
(-4.75819158e-02, 9.98867339e-01)
(-7.92499569e-02, 9.96854776e-01)
(-1.10838200e-01, 9.93838464e-01)
(-1.42314838e-01, 9.89821442e-01)
(-1.73648178e-01, 9.84807753e-01)
(-2.04806668e-01, 9.78802446e-01)
(-2.35758936e-01, 9.71811568e-01)
(-2.66473814e-01, 9.63842159e-01)
(-2.96920375e-01, 9.54902241e-01)
(-3.27067963e-01, 9.45000819e-01)
(-3.56886222e-01, 9.34147860e-01)
(-3.86345126e-01, 9.22354294e-01)
(-4.15415013e-01, 9.09631995e-01)
(-4.44066613e-01, 8.95993774e-01)
(-4.72271075e-01, 8.81453363e-01)
(-5.00000000e-01, 8.66025404e-01)
(-5.27225468e-01, 8.49725430e-01)
(-5.53920064e-01, 8.32569855e-01)
(-5.80056910e-01, 8.14575952e-01)
(-6.05609687e-01, 7.95761841e-01)
(-6.30552667e-01, 7.76146464e-01)
(-6.54860734e-01, 7.55749574e-01)
(-6.78509412e-01, 7.34591709e-01)
(-7.01474888e-01, 7.12694171e-01)
(-7.23734038e-01, 6.90079011e-01)
(-7.45264450e-01, 6.66769001e-01)
(-7.66044443e-01, 6.42787610e-01)
(-7.86053095e-01, 6.18158986e-01)
(-8.05270258e-01, 5.92907929e-01)
(-8.23676581e-01, 5.67059864e-01)
(-8.41253533e-01, 5.40640817e-01)
(-8.57983413e-01, 5.13677392e-01)
(-8.73849377e-01, 4.86196736e-01)
(-8.88835449e-01, 4.58226522e-01)
(-9.02926538e-01, 4.29794912e-01)
(-9.16108457e-01, 4.00930535e-01)
(-9.28367933e-01, 3.71662456e-01)
(-9.39692621e-01, 3.42020143e-01)
(-9.50071118e-01, 3.12033446e-01)
(-9.59492974e-01, 2.81732557e-01)
(-9.67948701e-01, 2.51147987e-01)
(-9.75429787e-01, 2.20310533e-01)
(-9.81928697e-01, 1.89251244e-01)
(-9.87438889e-01, 1.58001396e-01)
(-9.91954813e-01, 1.26592454e-01)
(-9.95471923e-01, 9.50560433e-02)
(-9.97986676e-01, 6.34239197e-02)
(-9.99496542e-01, 3.17279335e-02)
(-1.00000000e+00, 1.22464680e-16)
};\label{line:cyl:wall}

\addplot [black, thick, solid, forget plot]
coordinates {
(6.00000000e+00, 0.00000000e+00)
(1.00000000e+00, 0.00000000e+00)
};\label{line:cyl:outflow}

\addplot [black, thick, solid, forget plot]
coordinates {
(-6.00000000e+00, 0.00000000e+00)
(-1.00000000e+00, 0.00000000e+00)
};
\end{axis}
\end{tikzpicture}
	        \caption{
	                Schematic of cylinder domain and boundary conditions. Boundary
	                conditions: supersonic inlet (\ref{line:cyl:inflow}), supersonic
			outlet (\ref{line:cyl:outflow}), isothermal no-slip wall
			(\ref{line:cyl:wall}).
	        }
	        \label{fig:ns1_geom}
	\end{figure}
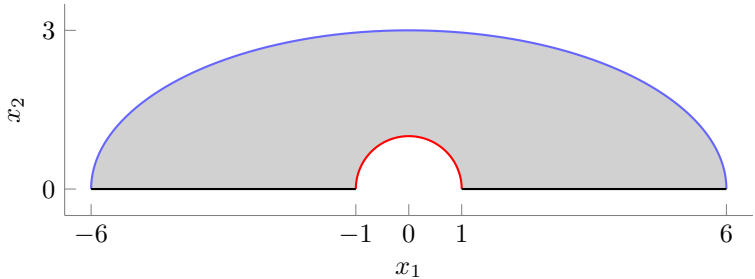

	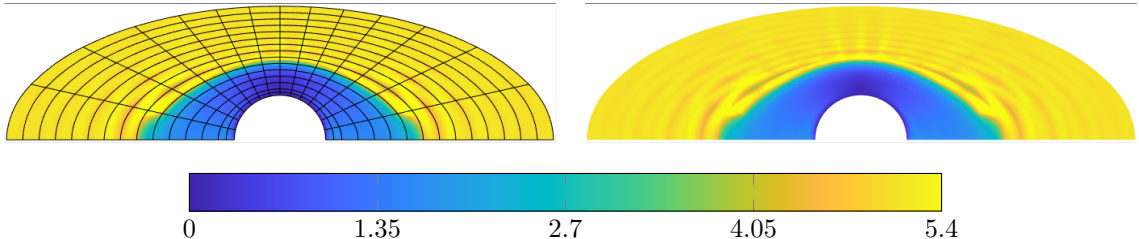
\begin{figure}[!htbp]
		\centering
		\begin{tikzpicture}[scale=0.8]
\begin{groupplot}[
  group style={
      group size=2 by 1,
      horizontal sep=0.5cm,
      vertical sep=0.5cm
  },
  width=0.65\textwidth,
  axis equal image,
  xlabel={$x_1$},
  ylabel={$x_2$},
  xtick = {-6.0, 0.0, 6.0},
  xticklabels={-6.0, 0.0, 6.0},
  ytick = {0.0, 3.0},
  xmin=-6.0, xmax=6.0,
  ymin=0, ymax=3.0
]

\nextgroupplot[xtick=\empty, xlabel={}, ytick=\empty, ylabel={}]
\addplot graphics [xmin=-6.0, xmax=6.0, ymin=0, ymax=3.0] {{_img/cyl_M5Re1e3_shkcap_M_msh}.png};

\nextgroupplot[xtick=\empty, xlabel={}, ytick=\empty, ylabel={}]
\addplot graphics [xmin=-6.0, xmax=6.0, ymin=0, ymax=3.0] {{_img/cyl_M5Re1e3_shkcap_M_womsh}.png};

\end{groupplot}
\node[anchor=north] at ($(group c1r1.south)!0.5!(group c2r1.south)$) {\colorbarMatlabParula{0}{1.35}{2.7}{4.05}{5.4}};
\end{tikzpicture}
   		\caption{
			Mesh and solution (shock capturing for $\mathrm{Re} = 10^2$) used
			to initialize the $rp$-adaptive HOIST method for the cylinder
			problem.
		}
		\label{fig:ns1_soln_ic}
	\end{figure}

	\begin{figure}[!htbp]
		\centering
	    	\begin{tikzpicture}[scale=0.8]
\begin{groupplot}[
  group style={
      group size=3 by 6,
      horizontal sep=0.25cm,
      vertical sep=0.25cm
  },
  width=0.5\textwidth,
  axis equal image,
  xtick = {-6.0, 0.0, 6.0},
  xticklabels={-6.0, 0.0, 6.0},
  ytick = \empty,
  xmin=-6.0, xmax=6.0,
  ymin=0, ymax=3.0
]

\nextgroupplot[xtick=\empty, xlabel={}]
\addplot graphics [xmin=-6.0, xmax=6.0, ymin=0, ymax=3.0] {{_img/cyl_M5Re1e3_0th_Mach_phydom_mesh}.png};

\nextgroupplot[xtick=\empty, xlabel={}, ytick=\empty, ylabel={}]
\addplot graphics [xmin=-6.0, xmax=6.0, ymin=0, ymax=3.0] {{_img/cyl_M5Re1e3_Error_0th_phydom}.png};

\nextgroupplot[xtick=\empty, xlabel={}, ytick=\empty, ylabel={}]
\addplot graphics [xmin=-6.0, xmax=6.0, ymin=0, ymax=3.0] {{_img/cyl_M5Re1e3_0th_pdeg_phydom}.png};

\nextgroupplot[xtick=\empty, xlabel={}]
\addplot graphics [xmin=-6.0, xmax=6.0, ymin=0, ymax=3.0] {{_img/cyl_M5Re1e3_1st_Mach_phydom_mesh}.png};

\nextgroupplot[xtick=\empty, xlabel={}, ytick=\empty, ylabel={}]
\addplot graphics [xmin=-6.0, xmax=6.0, ymin=0, ymax=3.0] {{_img/cyl_M5Re1e3_Error_1st_phydom}.png};

\nextgroupplot[xtick=\empty, xlabel={}, ytick=\empty, ylabel={}]
\addplot graphics [xmin=-6.0, xmax=6.0, ymin=0, ymax=3.0] {{_img/cyl_M5Re1e3_1st_pdeg_phydom}.png};

\nextgroupplot[xtick=\empty, xlabel={}]
\addplot graphics [xmin=-6.0, xmax=6.0, ymin=0, ymax=3.0] {{_img/cyl_M5Re1e3_2nd_Mach_phydom_mesh}.png};

\nextgroupplot[xtick=\empty, xlabel={}, ytick=\empty, ylabel={}]
\addplot graphics [xmin=-6.0, xmax=6.0, ymin=0, ymax=3.0] {{_img/cyl_M5Re1e3_Error_2nd_phydom}.png};

\nextgroupplot[xtick=\empty, xlabel={}, ytick=\empty, ylabel={}]
\addplot graphics [xmin=-6.0, xmax=6.0, ymin=0, ymax=3.0] {{_img/cyl_M5Re1e3_2nd_pdeg_phydom}.png};

\nextgroupplot[xtick=\empty, xlabel={}]
\addplot graphics [xmin=-6.0, xmax=6.0, ymin=0, ymax=3.0] {{_img/cyl_M5Re1e3_3rd_Mach_phydom_mesh}.png};

\nextgroupplot[xtick=\empty, xlabel={}, ytick=\empty, ylabel={}]
\addplot graphics [xmin=-6.0, xmax=6.0, ymin=0, ymax=3.0] {{_img/cyl_M5Re1e3_Error_3rd_phydom}.png};

\nextgroupplot[xtick=\empty, xlabel={}, ytick=\empty, ylabel={}]
\addplot graphics [xmin=-6.0, xmax=6.0, ymin=0, ymax=3.0] {{_img/cyl_M5Re1e3_3rd_pdeg_phydom}.png};

\nextgroupplot[xtick=\empty, xlabel={}]
\addplot graphics [xmin=-6.0, xmax=6.0, ymin=0, ymax=3.0] {{_img/cyl_M5Re1e3_4th_Mach_phydom_mesh}.png};

\nextgroupplot[xtick=\empty, xlabel={}, ytick=\empty, ylabel={}]
\addplot graphics [xmin=-6.0, xmax=6.0, ymin=0, ymax=3.0] {{_img/cyl_M5Re1e3_Error_4th_phydom}.png};

\nextgroupplot[xtick=\empty, xlabel={}, ytick=\empty, ylabel={}]
\addplot graphics [xmin=-6.0, xmax=6.0, ymin=0, ymax=3.0] {{_img/cyl_M5Re1e3_4th_pdeg_phydom}.png};

\nextgroupplot[xtick=\empty, xlabel={}]
\addplot graphics [xmin=-6.0, xmax=6.0, ymin=0, ymax=3.0] {{_img/cyl_M5Re1e3_5th_Mach_phydom_mesh}.png};

\nextgroupplot[xtick=\empty, xlabel={}, ytick=\empty, ylabel={}]
\addplot graphics [xmin=-6.0, xmax=6.0, ymin=0, ymax=3.0] {{_img/cyl_M5Re1e3_Error_5th_phydom}.png};

\nextgroupplot[xtick=\empty, xlabel={}, ytick=\empty, ylabel={}]
\addplot graphics [xmin=-6.0, xmax=6.0, ymin=0, ymax=3.0] {{_img/cyl_M5Re1e3_5th_pdeg_phydom}.png};

\end{groupplot}
\node[anchor=west] at ($(group c1r6.west) + (-0.25cm, -2)$) {\colorbarMatlabParulandscalse{0}{1.25}{2.5}{3.75}{5.10}{5.10}{5cm}};

\node[anchor=west] at ($(group c2r6.west) + (-0.25cm, -2)$) {\colorbarMatlabParulandscalse{4.5e-5}{0.26}{0.26}{0.26}{0.51}{0.51}{4.5cm}};

\node[anchor=west] at ($(group c3r6.west) + (-0.15cm, -2)$) {\colorbarMatlabParulandscalse{2}{3}{4}{5}{6}{7}{4.8cm}};
\end{tikzpicture}
	    	\caption{
			HOIST solution ($rp$-adaptation) (Mach number) (\textit{left}),
			the enriched residual error estimate (\textit{middle}),
			and the polynomial degree distribution
			(\textit{right}) in \textit{the physical domain} $\Gcal(\Omega_0)$
			to the cylinder
			problem. These quantities are provided prior to $p$-adaptation
			and after each $p$-adaptation iteration (\textit{top-to-bottom}).
		}
		\label{fig:ns1_soln_physdom}
	\end{figure}
	
	\begin{figure}[!htbp]
		\centering
	    	\begin{tikzpicture}[scale=0.8]
\begin{groupplot}[
  group style={
      group size=3 by 7,
      horizontal sep=0.25cm,
      vertical sep=0.25cm
  },
  width=0.5\textwidth,
  axis equal image,
  xtick = {-6.0, 0.0, 6.0},
  xticklabels={-6.0, 0.0, 6.0},
  ytick = \empty,
  xmin=-6.0, xmax=6.0,
  ymin=0, ymax=3.0
]

\nextgroupplot[xtick=\empty, xlabel={}]
\addplot graphics [xmin=-6.0, xmax=6.0, ymin=0, ymax=3.0] {{_img/cyl_M5Re1e3_0th_Mach_refdom}.png};

\nextgroupplot[xtick=\empty, xlabel={}, ytick=\empty, ylabel={}]
\addplot graphics [xmin=-6.0, xmax=6.0, ymin=0, ymax=3.0] {{_img/cyl_M5Re1e3_Error_0th_refdom}.png};

\nextgroupplot[xtick=\empty, xlabel={}, ytick=\empty, ylabel={}]
\addplot graphics [xmin=-6.0, xmax=6.0, ymin=0, ymax=3.0] {{_img/cyl_M5Re1e3_0th_pdeg_refdom}.png};

\nextgroupplot[xtick=\empty, xlabel={}]
\addplot graphics [xmin=-6.0, xmax=6.0, ymin=0, ymax=3.0] {{_img/cyl_M5Re1e3_1st_Mach_refdom}.png};

\nextgroupplot[xtick=\empty, xlabel={}, ytick=\empty, ylabel={}]
\addplot graphics [xmin=-6.0, xmax=6.0, ymin=0, ymax=3.0] {{_img/cyl_M5Re1e3_Error_1st_refdom}.png};

\nextgroupplot[xtick=\empty, xlabel={}, ytick=\empty, ylabel={}]
\addplot graphics [xmin=-6.0, xmax=6.0, ymin=0, ymax=3.0] {{_img/cyl_M5Re1e3_1st_pdeg_refdom}.png};

\nextgroupplot[xtick=\empty, xlabel={}]
\addplot graphics [xmin=-6.0, xmax=6.0, ymin=0, ymax=3.0] {{_img/cyl_M5Re1e3_2nd_Mach_refdom}.png};

\nextgroupplot[xtick=\empty, xlabel={}, ytick=\empty, ylabel={}]
\addplot graphics [xmin=-6.0, xmax=6.0, ymin=0, ymax=3.0] {{_img/cyl_M5Re1e3_Error_2nd_refdom}.png};

\nextgroupplot[xtick=\empty, xlabel={}, ytick=\empty, ylabel={}]
\addplot graphics [xmin=-6.0, xmax=6.0, ymin=0, ymax=3.0] {{_img/cyl_M5Re1e3_2nd_pdeg_refdom}.png};

\nextgroupplot[xtick=\empty, xlabel={}]
\addplot graphics [xmin=-6.0, xmax=6.0, ymin=0, ymax=3.0] {{_img/cyl_M5Re1e3_3rd_Mach_refdom}.png};

\nextgroupplot[xtick=\empty, xlabel={}, ytick=\empty, ylabel={}]
\addplot graphics [xmin=-6.0, xmax=6.0, ymin=0, ymax=3.0] {{_img/cyl_M5Re1e3_Error_3rd_refdom}.png};

\nextgroupplot[xtick=\empty, xlabel={}, ytick=\empty, ylabel={}]
\addplot graphics [xmin=-6.0, xmax=6.0, ymin=0, ymax=3.0] {{_img/cyl_M5Re1e3_3rd_pdeg_refdom}.png};

\nextgroupplot[xtick=\empty, xlabel={}]
\addplot graphics [xmin=-6.0, xmax=6.0, ymin=0, ymax=3.0] {{_img/cyl_M5Re1e3_4th_Mach_refdom}.png};

\nextgroupplot[xtick=\empty, xlabel={}, ytick=\empty, ylabel={}]
\addplot graphics [xmin=-6.0, xmax=6.0, ymin=0, ymax=3.0] {{_img/cyl_M5Re1e3_Error_4th_refdom}.png};

\nextgroupplot[xtick=\empty, xlabel={}, ytick=\empty, ylabel={}]
\addplot graphics [xmin=-6.0, xmax=6.0, ymin=0, ymax=3.0] {{_img/cyl_M5Re1e3_4th_pdeg_refdom}.png};

\nextgroupplot[xtick=\empty, xlabel={}]
\addplot graphics [xmin=-6.0, xmax=6.0, ymin=0, ymax=3.0] {{_img/cyl_M5Re1e3_5th_Mach_refdom}.png};

\nextgroupplot[xtick=\empty, xlabel={}, ytick=\empty, ylabel={}]
\addplot graphics [xmin=-6.0, xmax=6.0, ymin=0, ymax=3.0] {{_img/cyl_M5Re1e3_Error_5th_refdom}.png};

\nextgroupplot[xtick=\empty, xlabel={}, ytick=\empty, ylabel={}]
\addplot graphics [xmin=-6.0, xmax=6.0, ymin=0, ymax=3.0] {{_img/cyl_M5Re1e3_5th_pdeg_refdom}.png};

\end{groupplot}

\end{tikzpicture}
	    	\vspace{0.1cm}
	        \caption{
	                HOIST solution ($rp$-adaptation) (Mach number) (\textit{left}),
	                the enriched residual error estimate (\textit{middle}),
			and the polynomial degree distribution
	                (\textit{right}) in \textit{the reference domain} to the cylinder
	                problem. These quantities are provided prior to $p$-adaptation
	                and after each $p$-adaptation iteration (\textit{top-to-bottom}).
			Colorbar in Figure~\ref{fig:ns1_soln_physdom}.
	        }
		\label{fig:ns1_soln_refdom}
	\end{figure}
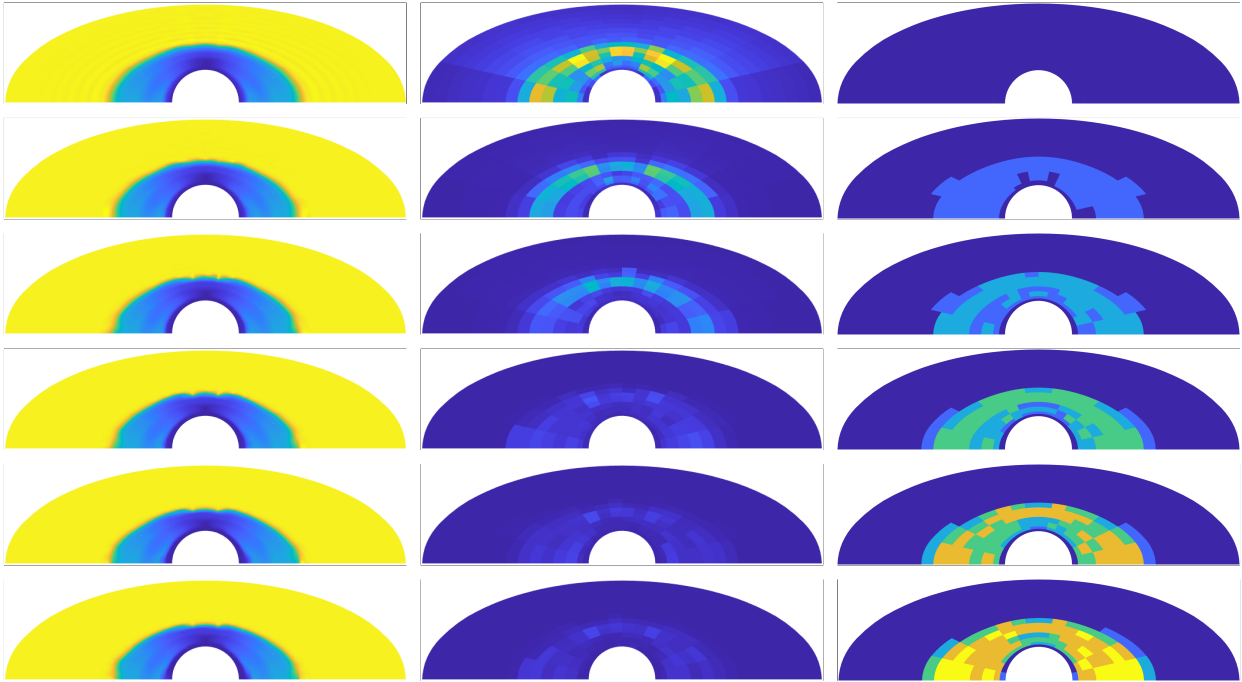
	
	\begin{figure}[!htbp]
		\centering
	    	\input{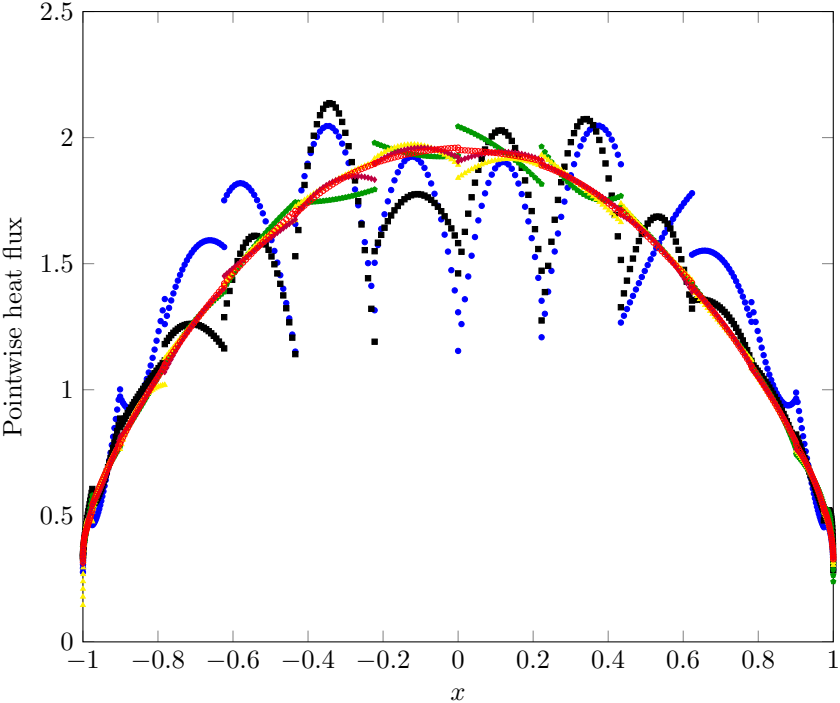}
		\caption{
			Pointwise heat flux over the cylinder surface computed
			from the HOIST solution ($rp$-adaptation) prior to $p$-adaptation
			(\ref{dots:0th}) and after the first (\ref{dots:1st}), second
			(\ref{dots:2nd}), third (\ref{dots:3rd}), fourth (\ref{dots:4th}),
			and fifth (\ref{dots:5th}) $p$-adaptation iteration.
		}
		\label{fig:ns1_hflx}
	\end{figure}

	\begin{figure}
		\centering
		\input{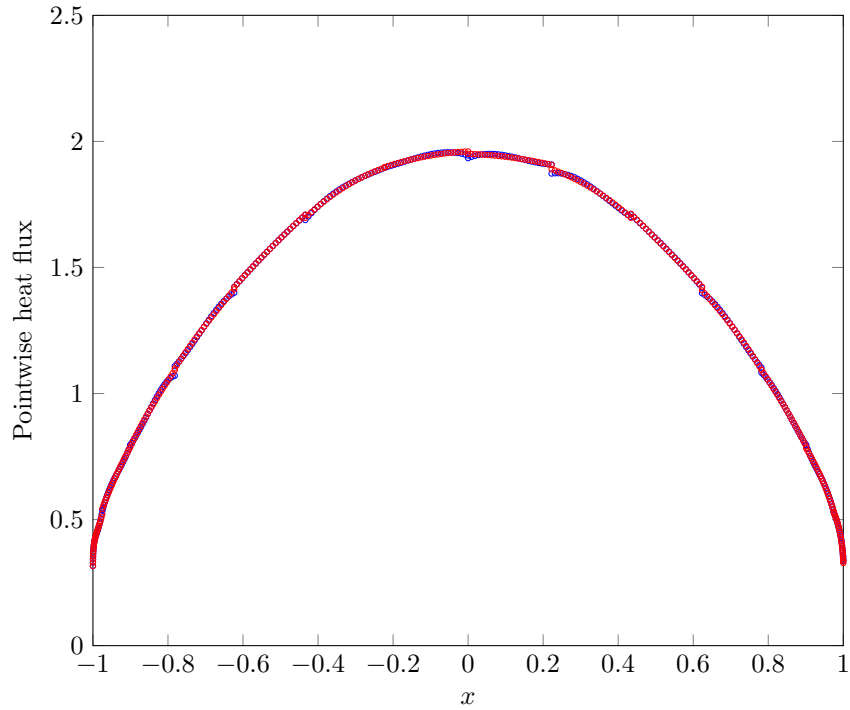}
		\caption{
			Pointwise heat flux profile over the cylinder surface
			computed from the HOIST solution ($rp$-adaptation) based
			on lifted viscous flux (\ref{dots:liftedvflux}) and
			via direct differentiation (\ref{dots:woliftvflux})
			(Remark~\ref{rem:liftgrad}) after fifth $p$-adaptation
			iteration.
		}
		\label{fig:ns1_hflux_liftvsnolift}
	\end{figure}
 
	\begin{figure}[!htbp]
		\centering
	    	\begin{tikzpicture}
    \begin{semilogxaxis}[width=0.6\textwidth, xmin=2.52000000e+03, xmax=1.05036000e+05, xlabel={Degrees of freedom}, ymin=3.69, ymax = 4.08, ytick={3.7, 3.75, 3.8, 3.85, 3.9, 3.95, 4.0, 4.05}, ylabel={Integrated heat flux}]

    \addplot [color=blue, mark=o]
    coordinates {
    ( 2.52000000e+03,  4.94775855e+00)
    ( 3.07200000e+03,  7.30356385e+00)
    ( 3.50400000e+03,  8.85062882e+00)
    ( 4.08000000e+03,  9.05792143e+00)
    ( 4.63200000e+03,  8.22245296e+00)
    ( 5.52000000e+03,  8.10595209e+00)
    ( 6.67200000e+03,  6.99568723e+00)
    ( 8.01600000e+03,  6.82798131e+00)
    ( 9.52800000e+03,  6.67013021e+00)
    ( 1.18560000e+04,  6.33065343e+00)
    ( 1.43760000e+04,  5.98501702e+00)
    ( 1.77000000e+04,  5.89900527e+00)
    ( 2.18040000e+04,  5.51930612e+00)
    ( 2.69520000e+04,  4.53356303e+00)
    ( 3.32640000e+04,  4.14777381e+00)
    ( 4.08600000e+04,  3.99264132e+00)
    ( 4.94160000e+04,  3.97403859e+00)
    ( 5.98560000e+04,  3.91648685e+00)
    ( 7.27200000e+04,  3.91725677e+00)
    ( 8.73360000e+04,  3.89880405e+00)
    ( 1.05036000e+05,  3.88715428e+00)};\label{dots:cyl_p1};

    \addplot [color=black, mark=o]
    coordinates {
    ( 5.04000000e+03,  5.49403741e+00)
    ( 5.80800000e+03,  5.46774744e+00)
    ( 6.62400000e+03,  4.32556079e+00)
    ( 7.68000000e+03,  4.01311988e+00)
    ( 9.26400000e+03,  3.96576280e+00)
    ( 1.10400000e+04,  3.91492959e+00)
    ( 1.40160000e+04,  3.89461159e+00)
    ( 1.65840000e+04,  3.88133941e+00)
    ( 2.04000000e+04,  3.88053426e+00)
    ( 2.53680000e+04,  3.87787642e+00)
    ( 3.11520000e+04,  3.86983132e+00)};\label{dots:cyl_p2};

    \addplot [color=orange, mark=o]
    coordinates {
    ( 8.40000000e+03,  5.66332031e+00)
    ( 9.36000000e+03,  4.46428409e+00)
    ( 1.08800000e+04,  3.95490535e+00)
    ( 1.36000000e+04,  3.89747110e+00)
    ( 1.65600000e+04,  3.87806406e+00)
    ( 2.15200000e+04,  3.87003772e+00)};\label{dots:cyl_p3};

    \addplot [red, mark=o ] coordinates {
    (5040,   3.9552544166084536)
    (6240,   3.820031577872347)
    (7436,   3.862276212770823)
    (8900,   3.8545650916726797)
    (10456,  3.8556580162052994)
    (12132, 3.8573185213358285)};\label{dots:cyl_ist};

\addplot[dotted, black] coordinates {
    (2.52000000e+03,  3.903554357065188)   (1.05036000e+05, 3.903554357065188)

    (2.52000000e+03, 3.826256250984689)   (1.05036000e+05, 3.826256250984689)
    };
   \addplot[dashed, black] coordinates {
    (2.52000000e+03,  3.8649053040249384)   (1.05036000e+05, 3.8649053040249384)
   };
    
   \node[anchor=south west, black] at (axis cs:2.55000000e+03, 3.91) {$1\%$ error};
    \end{semilogxaxis}
\end{tikzpicture}
	    	\caption{
			Convergence of the integrated heat flux for the HOIST method
			($rp$-adaptation) (\ref{dots:cyl_ist}) and the $h$-adaptive method
			with linear (\ref{dots:cyl_p1}),
			quadratic (\ref{dots:cyl_p2}), and cubic (\ref{dots:cyl_p3})
			constant polynomial degree for the cylinder problem. Both
			flow and mesh degrees of freedom are counted for the HOIST
			method ($N_\ubm^{(j)} + N_\xbm$), while only flow degrees
			of freedom are present for the $h$-adaptive method.
		}
		\label{fig:ns1_hflxerr}
	\end{figure}
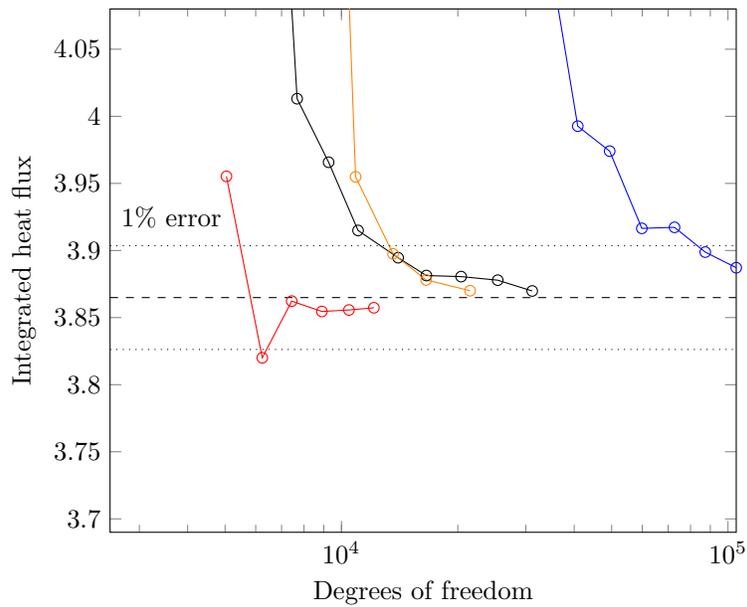
	
       \begin{figure}[!htbp]
        	\centering
		\begin{tikzpicture}[scale=0.8]
\begin{groupplot}[
  group style={
      group size=2 by 1,
      horizontal sep=1cm,
      vertical sep=0.5cm
  },
  width=0.7\textwidth,
  axis equal image,
  xmin=-6.0, xmax=6.0,
  ymin=0, ymax=3.0
]

\nextgroupplot[xtick=\empty, xlabel={}, ytick=\empty, ylabel={}]
\addplot graphics [xmin=-6.0, xmax=6.0, ymin=0, ymax=3.0] {{_img/cyl_masa_msh}.png};

\nextgroupplot[xtick=\empty, xlabel={}, ytick=\empty, ylabel={}]
\addplot graphics [xmin=-6.0, xmax=6.0, ymin=0, ymax=3.0] {{_img/cyl_masa_Mach_womsh}.png};

\end{groupplot}
\end{tikzpicture}
        	\caption{
        		The $h$-adapted mesh (\textit{left}) and corresponding solution
        		(Mach number) (\textit{right}) for the cylinder problem
        		using integrated heat flux-based dual-weighted residual
			adaptation used for the comparative study with the
        		HOIST method ($rp$-adaptation).
        	}
        	\label{fig:ns1_amr}
        \end{figure}
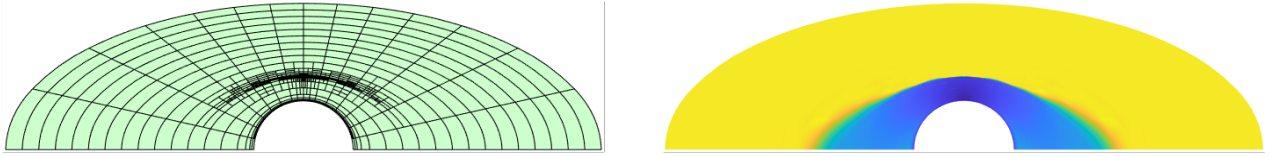

	\begin{figure}[!htbp]
		\centering
		\input{_tikz/cyl_M5Re1e3_c_nrm.tikz}
		\caption{
			Convergence of the DG residual
			$\norm{\rbm^{(j)}(\ubm_k^{(j)},\phibold(\ybm_k^{(j)}); \Xi_i)}$
			(\ref{line:c_nrm}) and enriched DG residual
			$\norm{\Rbm_\rho^{(j)}(\ubm_k^{(j)},\phibold(\ybm_k^{(j)}); \Xi_i)}$
			(\ref{line:Rerr_nrm}) across all Reynolds number continuation stage
			$\{\Xi_i\}_{i=1}^{10}$ and $p$-adaptivity iterations
			$\{p_j\}_{j=0}^5$. Vertical lines indicate the end of a
			continuation stage (\ref{line:recon_c_nrm}) (i.e., transition from
			$\Xi_i$ to $\Xi_{i+1}$) or the end of a $p$-adaptation iteration
			(\ref{line:padapt_c_nrm}) (i.e., transition from $p_j$ to
			$p_{j+1}$) for the cylinder problem.
	 	}
		\label{fig:ns1_sqp}
	\end{figure}

}

\subsubsection{Hypersonic flow over cylinder: importance of Reynolds continuation}
\label{subsubsec:bow_recont}
In this section, we consider the hypersonic cylinder problem in Section~\ref{subsubsec:bow}
and demonstrate the importance of Reynolds number continuation with all other
solver and formulation parameters identical to those in Section~\ref{subsubsec:bow}.
Without Reynolds number continuation, i.e., HOIST is directly applied to the target
Reynolds number $\Xi_1 = \mathrm{Re} = 1000$, carbuncle
phenomena \cite{peery1988blunt,robinet2000shock} arise,
which are particularly detrimental in an $r$-adaptation setting because the mesh tracks
the forming carbuncle (Figure~\ref{fig:ns1_soln_carb}).

\ifbool{fastcompile}{}{
	\begin{figure}[!htbp]
    		\centering
    		\begin{tikzpicture}[scale=0.8]
\begin{groupplot}[
  group style={
      group size=2 by 1,
      horizontal sep=1cm,
      vertical sep=0.5cm
  },
  width=0.7\textwidth,
  axis equal image,
  xlabel={$x_1$},
  ylabel={$x_2$},
  xtick = {-6.0, 0.0, 6.0},
  xticklabels={-6.0, 0.0, 6.0},
  ytick = {0.0, 3.0},
  xmin=-6.0, xmax=6.0,
  ymin=0, ymax=3.0
]
\nextgroupplot[xlabel={}, ylabel={}, xtick=\empty, ytick=\empty]
\addplot graphics [xmin=-6.0, xmax=6.0, ymin=0, ymax=3.0] {{_img/cyl_M5Re1e3_ReConfinal_M_msh}.png};

\nextgroupplot[xlabel={}, ylabel={}, xtick=\empty, ytick=\empty]
\addplot graphics [xmin=-6.0, xmax=6.0, ymin=0, ymax=3.0] {{_img/cyl_M5Re1e3_noReCon_M_msh}.png};

\end{groupplot}
\end{tikzpicture}
    		\caption{
			HOIST solution (Mach number) ($r$-adaptation) with the
			Reynolds number continuation strategy of Section~\ref{subsec:cont}
			(\textit{left}) and without continuation
			($\Xi_1 = \mathrm{Re} = 1000$) (\textit{right})
			for the cylinder problem.
			Colorbar in Figure~\ref{fig:ns1_soln_ic}.
		}
    		\label{fig:ns1_soln_carb}
	\end{figure}
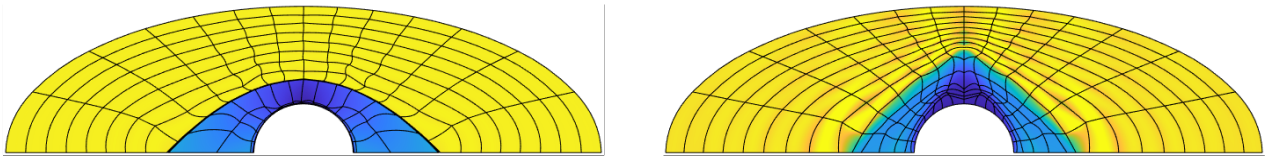

}

\subsubsection{Hypersonic flow over cylinder: impact of enrichment degree}
\label{subsubsec:bow_enrch}
Next, we consider the hypersonic cylinder problem in Section~\ref{subsubsec:bow}
and demonstrate the impact of the enrichment degree $\Delta$. In
particular, we compare the choice $\Delta = 1$ to the $\Delta = 2$
choice used throughout this work with all other solver parameters
identical to those in Section~\ref{subsubsec:bow}. We find that $\Delta = 1$ leads to
substantially less mesh compression in the boundary layer
(Figure~\ref{fig:ns1_soln_delta1}) and worse prediction of the
heat flux profile (Figure~\ref{fig:ns1_hflux_delta1}).

\ifbool{fastcompile}{}{
	\begin{figure}[!htbp]
		\centering
    		\begin{tikzpicture}[scale=0.8]
\begin{groupplot}[
  group style={
      group size=2 by 2,
      horizontal sep=1cm,
      vertical sep=0.5cm
  },
  width=0.7\textwidth,
  axis equal image,
  xlabel={$x_1$},
  ylabel={$x_2$},
  xtick = {-6.0, 0.0, 6.0},
  xticklabels={-6.0, 0.0, 6.0},
  ytick = {0.0, 3.0},
]

\nextgroupplot[xlabel={}, ylabel={}, xtick=\empty, ytick=\empty]
\addplot graphics [xmin=-3.0, xmax=3.0, ymin=0, ymax=2.0] {{_img/cyl_M5Re1e3_final_T_msh}.png};

\draw[red, thick] (axis cs:-0.02,0.998) rectangle (axis cs:0.02,1.003);

\nextgroupplot[xlabel={}, ylabel={}, xtick=\empty, ytick=\empty]
\addplot graphics [xmin=-3.0, xmax=3.0, ymin=0, ymax=2.0] {{_img/cyl_M5Re1e3_enrespp1_T_msh}.png};

\draw[red, thick] (axis cs:-0.02,0.998) rectangle (axis cs:0.02,1.003);

\nextgroupplot[xlabel={}, ylabel={}, xtick=\empty, ytick=\empty]
\addplot graphics [xmin=-0.02, xmax=0.02, ymin=0.998, ymax=1.003] {{_img/cyl_M5Re1e3_final_T_zoomedmsh0}.png};

\nextgroupplot[xlabel={}, ylabel={}, xtick=\empty, ytick=\empty]
\addplot graphics [xmin=-0.02, xmax=0.02, ymin=0.998, ymax=1.003] {{_img/cyl_M5Re1e3_enrespp1_T_zoomedmsh0}.png};

\end{groupplot}
\node[anchor=north] at ($(group c1r2.south)!0.5!(group c2r2.south)$) {\colorbarMatlabParulandscalse{1}{2}{3}{4}{5}{6}{10cm}};
\end{tikzpicture}
    		\caption{
			HOIST solution (temperature) ($rp$-adaptation) based on
			test space enrichment $\Delta = 2$ (\textit{left}) and
			$\Delta = 1$ (\textit{right}) for the cylinder problem. 
		}
		\label{fig:ns1_soln_delta1}
	\end{figure}
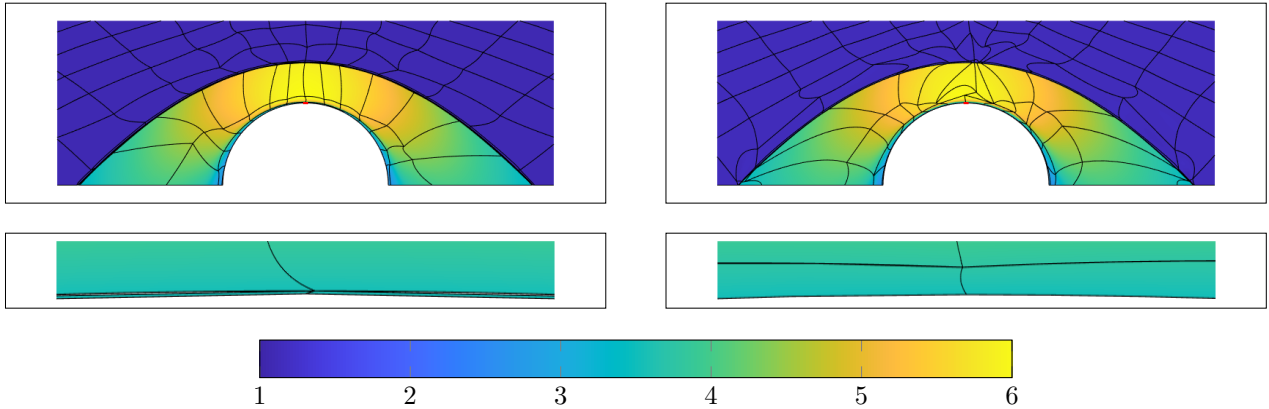

	\begin{figure}
		\centering
		\input{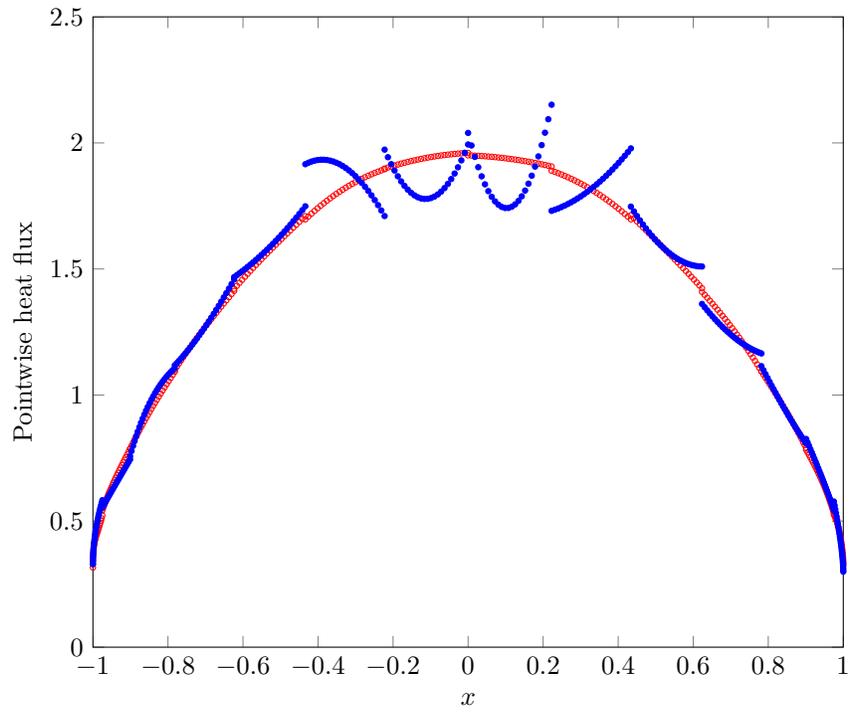}
		\caption{
			Pointwise heat flux profile over the cylinder surface
			computed from the HOIST solution ($rp$-adaptation) based
			on test space enrichment $\Delta = 2$ (\ref{line:enrespp2}) and
			$\Delta = 1$ (\ref{line:enrespp1}).
		}
		\label{fig:ns1_hflux_delta1}
	\end{figure}
}

\subsubsection{Hypersonic flow over cylinder: impact of boundary residual weighting}
\label{subsubsec:bow_bndwght}
Next, we consider the hypersonic cylinder problem in Section~\ref{subsubsec:bow}
and demonstrate the impact of boundary residual weighting $\lambda$. In
particular, we compare the choice $\lambda = 1$ (unweighted) to the
$\lambda = 10$ choice used in Section~\ref{subsubsec:bow}. All other solver parameters
are identical to those in Section~\ref{subsubsec:bow}. We find that the unweighted
scenario does not sufficiently compress the mesh to resolve the
boundary layer (Figure~\ref{fig:ns1_soln_bndwght}), which leads to
inaccurate heat flux profiles (Figure~\ref{fig:ns1_hflux_bndwght}).

\ifbool{fastcompile}{}{
	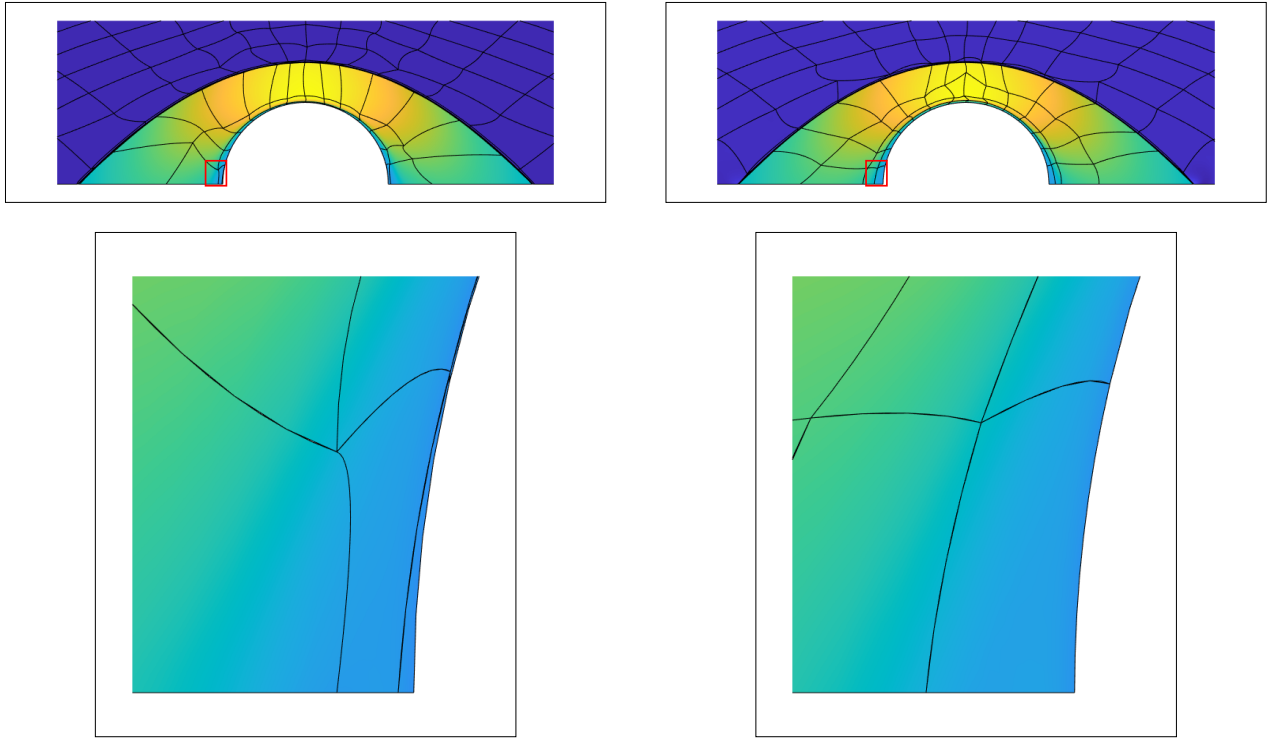
\begin{figure}[!htbp]
    		\centering
    		\begin{tikzpicture}[scale=0.8]
\begin{groupplot}[
  group style={
      group size=2 by 2,
      horizontal sep=1cm,
      vertical sep=0.5cm
  },
  width=0.7\textwidth,
  axis equal image,
  xlabel={$x_1$},
  ylabel={$x_2$},
  xtick = {-6.0, 0.0, 6.0},
  xticklabels={-6.0, 0.0, 6.0},
  ytick = {0.0, 3.0},
]

\nextgroupplot[xlabel={}, ylabel={}, xtick=\empty, ytick=\empty]
\addplot graphics [xmin=-3.0, xmax=3.0, ymin=0, ymax=2.0] {{_img/cyl_M5Re1e3_final_T_msh}.png};

\draw[red, thick] (axis cs:-1.2,0) rectangle (axis cs:-0.95,0.3);

\nextgroupplot[xlabel={}, ylabel={}, xtick=\empty, ytick=\empty]
\addplot graphics [xmin=-3.0, xmax=3.0, ymin=0, ymax=2.0] {{_img/cyl_M5Re1e3_nowelem_T_msh}.png};

\draw[red, thick] (axis cs:-1.2,0) rectangle (axis cs:-0.95,0.3);

\nextgroupplot[xlabel={}, ylabel={}, xtick=\empty, ytick=\empty]
\addplot graphics [xmin=-1.2, xmax=-0.95, ymin=0, ymax=0.3] {{_img/cyl_M5Re1e3_final_T_zoomedmsh}.png};

\nextgroupplot[xlabel={}, ylabel={}, xtick=\empty, ytick=\empty]
\addplot graphics [xmin=-1.2, xmax=-0.95, ymin=0, ymax=0.3] {{_img/cyl_M5Re1e3_nowelem_T_zoomedmsh}.png};

\end{groupplot}

\end{tikzpicture}
   		 \caption{
			HOIST solution (temperature) ($rp$-adaptation) based on
			boundary residual scaling $\lambda = 10$ (\textit{left})
			and $\lambda = 1$ (unweighted) (\textit{right}).
			Colorbar in Figure~\ref{fig:ns1_soln_delta1}. 
		 }
		 \label{fig:ns1_soln_bndwght}
	\end{figure}

	\begin{figure}
		\centering
		\input{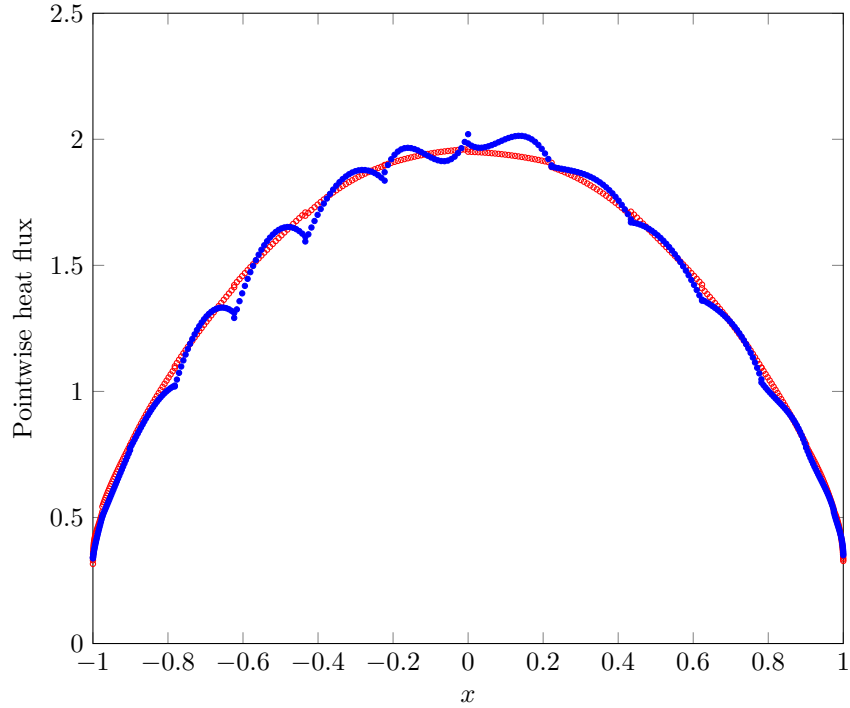}
		\caption{
			Pointwise heat flux profile over the cylinder surface
                        computed from the HOIST solution ($rp$-adaptation) based
			boundary residual scaling $\lambda = 10$ (\ref{line:welem})
			and $\lambda = 1$ (unweighted) (\ref{line:nowelem}).
		}
		\label{fig:ns1_hflux_bndwght}
	\end{figure}
}

\subsubsection{Hypersonic flow over cylinder: $p$-adaptation error indicator}
\label{subsubsec:bow_snsr}
Finally, we consider the hypersonic cylinder problem in Section~\ref{subsubsec:bow}
and demonstrate the performance of the HOIST method ($rp$-adaptation)
with different error indicators driving $p$-adaptation. All of the
solver parameters are identical to those in Section~\ref{subsubsec:bow}. The
dual-weighted residual error indicator based on the integrated
heat flux concentrates the highest polynomial degrees near the
stagnation line, whereas the feature-based sensor based on the
magnitude of the velocity gradient only enriches the polynomial
degree in the shock layer (Figure~\ref{fig:ns1_pdeg_errindic}).

The distribution of polynomial degree produced by the dual-weighted
residual sensor is not ideal because the lowest polynomial degree
(two) is used in the lateral regions (minor contributions
to the integrated heat flux), which means $r$-adaptation alone
is used to resolve the bow shock in these regions. Regardless,
the integrated heat flux reaches sub-$1\%$ error using this
trial space, which is not surprising because the
polynomial distribution is tailored to accurately predict
this quantity of interest. However, because the shock is
not adequately resolved, the PTC solver fails to deeply
converge the solution without additional stabilization.
The dual-weighted residual method is known to be more effective than residual-based
methods for approximating quantities of interest~\cite{Becker_2001_DWR,fidkowski2011review}.
However, in the present formulation, mesh optimization plays a role of improving not
just accuracy by providing better approximation spaces but also stability by providing
shock-aligned meshes. The former is desired, but the latter is required; the dual-weighted
residual method by construction accounts for the former but not the latter, and hence is
not best suited for the present formulation.

Furthermore, the distribution of polynomial degree produced by the
feature-based sensor is suboptimal despite highly focusing the $p$-refinement
in the bow shock because (1) the boundary layer is not being refined
and (2) the bow shock will always be targeted for refinement, even
once it has been adequately resolved due to the use of a feature-based
(not error-based) sensor. Unlike the trial space produced from the
dual-weighted residual indicator, the integrated heat flux is inaccurate
using this feature-adapted trial space due to lack of boundary layer
refinement. However, deep convergence using PTC is obtained
because the shock is adequately resolved.

The enriched residual sensor used in Section~\ref{subsubsec:bow} combines the
best aspects of the other two sensors. It is an error-based sensor
so refinement will be placed where resolution is needed, regardless
of the steepness of the features. This ensures sufficient refinement
will be placed in both the shock and boundary layer; the former enables
deep convergence using PTC without stabilization and the latter leads
to accurate heat flux calculations.

\ifbool{fastcompile}{}{
	\begin{figure}[!htbp]
		\centering
		\begin{tikzpicture}[scale=0.8]
\begin{groupplot}[
  group style={
      group size=2 by 2,
      horizontal sep=1cm,
      vertical sep=0.5cm
  },
  width=0.7\textwidth,
  axis equal image,
  xlabel={$x_1$},
  ylabel={$x_2$},
  xtick = {-6.0, 0.0, 6.0},
  xticklabels={-6.0, 0.0, 6.0},
  ytick = {0.0, 3.0},
  xmin=-6.0, xmax=6.0,
  ymin=0, ymax=3.0
]

\nextgroupplot[xlabel={}, ylabel={}, xtick=\empty, ytick=\empty]
\addplot graphics [xmin=-6.0, xmax=6.0, ymin=0, ymax=3.0] {{_img/cyl_dwr_5th_pdeg_phydom}.png};

\nextgroupplot[xlabel={}, ylabel={}, xtick=\empty, ytick=\empty]
\addplot graphics [xmin=-6.0, xmax=6.0, ymin=0, ymax=3.0] {{_img/cyl_dwr_5th_pdeg_refdom}.png};

\nextgroupplot[xlabel={}, ylabel={}, xtick=\empty, ytick=\empty]
\addplot graphics [xmin=-6.0, xmax=6.0, ymin=0, ymax=3.0] {{_img/cyl_velgrad_5th_pdeg_phydom}.png};

\nextgroupplot[xlabel={}, ylabel={}, xtick=\empty, ytick=\empty]
\addplot graphics [xmin=-6.0, xmax=6.0, ymin=0, ymax=3.0] {{_img/cyl_velgrad_5th_pdeg_refdom}.png};

\end{groupplot}
\node[anchor=north] at ($(group c1r2.south)!0.5!(group c2r2.south)$) {\colorbarMatlabParulandscalse{2}{3}{4}{5}{6}{7}{10cm}};
\end{tikzpicture}
		\caption{
			Distribution of the polynomial degree for the HOIST method
			($rp$-adaptation) after five $p$-adaption iterations using
			the dual-weighted residual sensor $\hat{s}^\mathrm{dwr}$
			based on the integrated heat flux quantity (\textit{top})
			and the magnitude of the velocity gradient $\hat{s}^\mathrm{fbs}$
			(\textit{bottom}). The sensors are provided on the
			physical domain (\textit{left}) and reference domain
			(\textit{right}) for clarity.
		}
		\label{fig:ns1_pdeg_errindic}
	\end{figure}
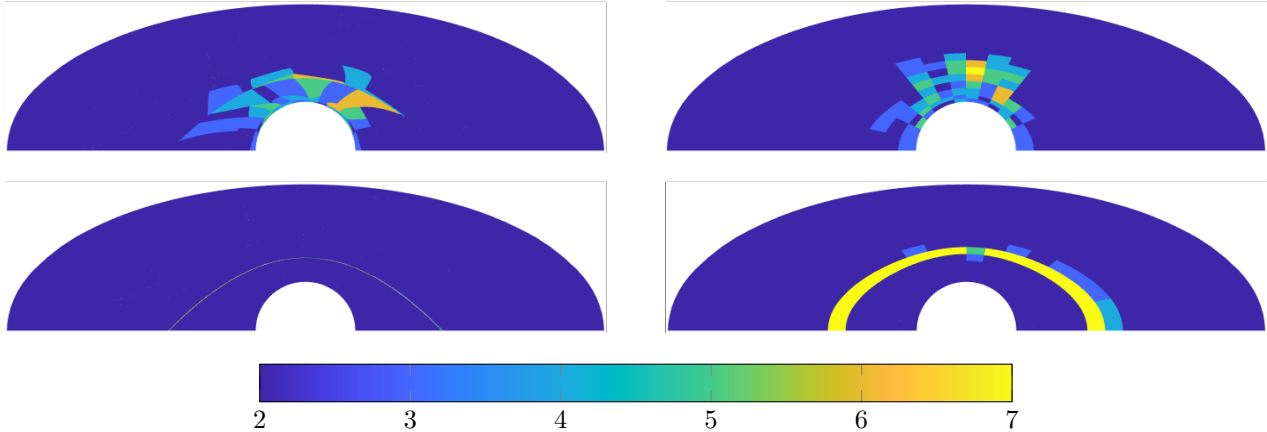
}

\section{Conclusion}
\label{sec:conclude}
In this work, we extended the implicit shock tracking method
of \cite{zahr2018shktrk,zahr2020implicit,huang2022robust} from a method
that automatically aligns element faces with discontinuities in the solution
field to one that resolves strong, rapid transitions (e.g., viscous shock waves
and boundary layers) using aggressive $rp$-adaptivity. While both approaches share
the same optimization formulation, several key innovations were necessary to resolve
the viscous features. First, an appropriate discretization for a second-order
conservation law is required (we use IPDG in this work). Second, an objective function
corresponding to a test space enriched by two polynomial degrees led to higher quality
solutions and meshes than the standard approach of enriching by a single polynomial
degree. Third, the incorporation of viscosity continuation was critical for extreme
yet robust compression of elements into transitions and to avoid carbuncles.
Fourth, $p$-adaptivity was used to locally increase the polynomial degree in the
regions with the most error, which limited the amount of mesh compression necessary.
Fifth, an artificial elemental weighting was applied to the residual of the objective
function to ensure the boundary layer was prioritized equally to the shock in the
objective function. Finally, a host of other measures were incorporated to improve
robustness of the solver, including constraints and modifications to the SQP step
and strategies to regularize the objective function Levenberg--Marquardt Hessian
approximation. A collection of numerical experiments show the new HOIST method
effectively resolves both shock waves and boundary layers in viscous shock-dominated
flows and compares favorably in terms of accuracy per degree of freedom to aggressive
$h$-adaptation with high-order methods. Furthermore, the $rp$-adaptive approach
effectively predicts the heat flux profile produced by hypersonic flow over a cylinder.

Relevant future research includes the development parallel preconditioners for the
linearized optimality system to effectively leverage high-performance computing, and
application of the proposed approach to realistic hypersonic flows. Another promising
avenue of research is the development of an implicit shock tracking method for viscous
flows that approximates strong shock waves as discontinuities. Such an approach could
lead to highly accurate approximations for high Reynolds flow without the need for
extreme mesh compression or $p$-adaptivity, and potentially use coarser grids than those
used in this work.

\section*{Acknowledgments}
This work is supported by AFOSR award numbers FA9550-20-1-0236,
FA9550-22-1-0002, FA9550-22-1-0004, ONR award number
N00014-22-1-2299, and NSF award number CBET-2338843.
The content of this publication does not necessarily reflect the position
or policy of any of these supporters, and no official endorsement should
be inferred.

\bibliographystyle{plain}
\bibliography{biblio}

\end{document}